\pdfoutput=1

\documentclass[preprint,10pt]{elsarticle}

\usepackage{fullpage}
\usepackage[colorlinks=true]{hyperref}
\usepackage{float}

\usepackage{amsmath,amssymb,amsfonts,amsthm}
\theoremstyle{definition}

\theoremstyle{lemma}
\newtheorem{lemma}{Lemma}
\newtheorem*{remark}{Remark}
\theoremstyle{theorem}
\newtheorem{theorem}{Theorem}
\theoremstyle{assumption}

\usepackage[titletoc,toc,title]{appendix}
\usepackage{array} 
\usepackage{listings}
\usepackage{mathtools}
\usepackage{pdfpages}
\usepackage[textsize=footnotesize,color=green]{todonotes}
\usepackage{bm}
\usepackage{bbm}

\usepackage{tikz}
\usetikzlibrary{shapes,arrows}
\usepackage[normalem]{ulem}
\usepackage{hhline}

\usepackage{graphicx}
\usepackage{subfig}
\usepackage{color}

\usepackage{pgfplots}
\usepackage{pgfplotstable}
\definecolor{markercolor}{RGB}{124.9, 255, 160.65}
\pgfplotsset{
compat=1.3,
width=10cm,
tick label style={font=\small},
label style={font=\small},
legend style={font=\small}
}

\usetikzlibrary{calc}
\usetikzlibrary{intersections} 

\newcommand{\logLogSlopeTriangleFlip}[5]
{

    \pgfplotsextra
    {
        \pgfkeysgetvalue{/pgfplots/xmin}{\xmin}
        \pgfkeysgetvalue{/pgfplots/xmax}{\xmax}
        \pgfkeysgetvalue{/pgfplots/ymin}{\ymin}
        \pgfkeysgetvalue{/pgfplots/ymax}{\ymax}

        \pgfmathsetmacro{\xBrel}{#1-#2}
        \pgfmathsetmacro{\yBrel}{#3}
        \pgfmathsetmacro{\xCrel}{#1}

        \pgfmathsetmacro{\lnxB}{\xmin*(1-(#1-#2))+\xmax*(#1-#2)} 
        \pgfmathsetmacro{\lnxA}{\xmin*(1-#1)+\xmax*#1} 
        \pgfmathsetmacro{\lnyA}{\ymin*(1-#3)+\ymax*#3} 
        \pgfmathsetmacro{\lnyC}{\lnyA+#4*(\lnxA-\lnxB)}
        \pgfmathsetmacro{\yCrel}{\lnyC-\ymin)/(\ymax-\ymin)} 

	\pgfmathsetmacro{\xArel}{\xBrel}
        \pgfmathsetmacro{\yArel}{\yCrel}

        \coordinate (A) at (rel axis cs:\xArel,\yArel);
        \coordinate (B) at (rel axis cs:\xBrel,\yBrel);
        \coordinate (C) at (rel axis cs:\xCrel,\yCrel);

        \draw[#5]   (A)-- node[pos=0.5,anchor=east] {#4}
                    (B)-- 
                    (C)-- node[pos=0.5,anchor=south] {}
                    cycle;
    }
}

\usepackage{stmaryrd}
\newcommand{\myfrac}[2]{\displaystyle \frac{#1}{#2}}

\newcommand{\LRp}[1]{\left( #1 \right)}

\newcommand{\LRc}[1]{\left\{ #1 \right\}}

\newcommand{\jump}[1] {\ensuremath{\llbracket#1\rrbracket}}
\newcommand{\avg}[1] {\ensuremath{\LRc{\!\{#1\}\!}}}

\newcommand{\eval}[2][\right]{\relax
  \ifx#1\right\relax \left.\fi#2#1\rvert}

\newcolumntype{C}[1]{>{\centering\let\newline\\\arraybackslash\hspace{0pt}}m{#1}}

\makeatletter
\renewcommand\d[1]{\mspace{6mu}\mathrm{d}#1\@ifnextchar\d{\mspace{-3mu}}{}}
\makeatother

\graphicspath{{./figs/}}

\begin{document}


\begin{frontmatter}
\title{A weight-adjusted discontinuous Galerkin method for the poroelastic wave equation: penalty fluxes and micro-heterogeneities}

\author{Khemraj Shukla$^\mathit{1}$\corref{cor}}
\ead{rajexplo@gmail.com}

\author{Jesse Chan$^\mathit{1}$ \corref{}}
\ead{jesse.chan@rice.edu}

\author{Maarten V. de Hoop$^\mathit{1}$ \corref{}}
\ead{mdehoop@rice.edu}

\author{Priyank Jaiswal$^\mathit{2}$ \corref{}}
\ead{priyank.jaiswal@okstate.edu}

\cortext[cor]{corresponding author}

\address{$^\mathit{1}$Department of Computational and Applied Mathematics, Rice University, 6100 Main St, Houston, TX, 77005\\
$^\mathit{2}$Boone Pickens School of Geology, Oklahoma State University, Stillwater, OK, 74078}

\begin{abstract}
We introduce a high-order weight-adjusted discontinuous Galerkin (WADG) scheme for the numerical solution of three-dimensional (3D) wave propagation problems in anisotropic porous media. We use a coupled first-order symmetric stress-velocity formulation \cite{lemoine2013, lemoine2016}. Careful attention is directed at (a) the derivation of an energy-stable penalty-based numerical flux, which offers high-order accuracy in presence of material discontinuities, and (b) proper treatment of micro-heterogeneities (sub-element variations) in the numerical scheme. The use of a penalty-based numerical flux avoids the diagonalization of Jacobian matrices into polarized wave constituents necessary when solving element-wise Riemann problems.  Micro-heterogeneities are accurately and stably incorporated in the numerical scheme using easily-invertible weight-adjusted mass matrices \cite{chan2018}. The convergence of the proposed numerical scheme is proven and verified by using convergence studies against analytical plane wave solutions. The proposed method is also compared against an existing implementation using the spectral element method to solve the poroelastic wave equation \cite{morency2008}.
\end{abstract}
\end{frontmatter}


\section{Introduction}
Theories of poroelasticity deal with the physics of elastic wave propagation in porous media. These physics are applicable where pore-filling materials are of interest, such as oil and gas exploration, gas hydrate detection, hydrogeology, seismic monitoring of CO$_2$ storage, and medical imaging. The most popular theory of poroelastic wave propagation was pioneered by Maurice Biot \cite{biot1956I,biot1962}. Biot obtained a dynamical system of equations describing wave propagation in a porous medium, saturated with a single phase fluid with an ability to flow through the pore networks. Biot's system of equations is divided in two categories, based on the material's characteristic frequency $(\omega_c)$ \cite{biot1962, carcione1995, ozdenvar1997}, which describes the nature of fluid flow (laminar or non-laminar) induced by wave motion. Biot's theory of poroelasticity predicts two compressional waves, the fast P wave and slow P wave, and one shear wave. 

{Biot's theory of poroelasticity predicts two compressional waves (a fast P wave and slow P wave) and one shear wave. In the fast P wave, the fluid and matrix are locked together and move in-phase. This locking arises through the action of viscous or inertial forces which result from wave-induced fluid flow.  Viscous forces arise due to friction between the layers of the viscous fluid (laminar flow), whereas inertial forces result from the momentum of the fluid flow (non-laminar flow).  The slow P wave results from relative motion between the fluid in the pores and the matrix.  Biot's system of equations is divided in two categories based on the material's characteristic frequency $(\omega_c)$ \cite{biot1962, carcione1995, ozdenvar1997}, which describes the nature of fluid flow (laminar or non-laminar) induced by wave motion.  Due to the relative motion between the solid and the fluid, a pore boundary layer is formed.  The pore boundary layer is related to the viscous-skin depth, which is the distance where the fluid particle velocity attains a certain percentage of its maximum value. }

A major cause of energy dissipation in porous media is a dissipative (or drag) ``memory'' force resulting from the relative motion between the solid skeleton and pore fluids. The dissipative force is incorporated in the equations of motion through a time dependent viscodynamic operator $\psi(t)$. In the low-frequency regime, the dissipative force is linearly proportional to the relative velocity between solid and fluid and $\psi(t)$ is compactly supported. This occurs because the viscous skin depth in the low-frequency regime is greater than the pore size, and thus does not form any pore boundary layer. However, in the high frequency regime, the dominance of inertial forces over viscous forces causes the formation of a pore boundary viscous layer. Thus, in principle, the relaxation mechanism of the medium due to drag force is expressed as a memory kernel and $\psi(t)$ is incorporated into the equation of motion through a convolution operator \cite{carcione1996}. Our current work focuses on poroelasticity in the low-frequency regime $(\omega_c)$.  In subsequent work, we will present numerical methods for the broad-band Biot system based on the Johnson-Koplik-Dashen (JKD) model of dynamic permeability \cite{johnson1987}, which will unify both the high and low frequency regime.

Numerical simulations of the poroelastic wave equation  in the low-frequency regime have been previously explored in the literature. The first numerical simulations were performed using the reflectivity method for flat layers \cite{stern1985,yamamoto1988} and cylindrical structures~\cite{rosenbaum1974}. Numerical simulation of the poroelastic wave equation using direct or grid based methods, such as finite difference (FD) and pseudo-spectral (PS) methods, dates back to the 1970s \cite{carcione2010}. Most of the methods presented in \cite{carcione2010} regard pseudospectral \cite{carcione1995,carcione1996}, staggered pseudospectral \cite{ozdenvar1997} and finite-difference methods \cite{garg1974,dai1995} and are based on 2D structured meshes. Santos and Ore{\~n}a \cite{santos1986} used a finite-element method to solve the poroelastic wave equation (second-order form) using quadrilateral meshes for spatial discretization. Recent work on the numerical solution of orthotropic poroelasticity is reported by Lemoine et al. \cite{lemoine2013}, using a finite-volume method on structured meshes.

In the present study, we introduce a high-order numerical scheme based on the discontinuous Galerkin method to solve the 3D poroelastic wave equation on unstructured tetrahedral meshes. High order methods provide one avenue towards improving fidelity in numerical simulations while maintaining reasonable computational costs, and methods which can accommodate unstructured meshes are desirable for problems with complex geometries. Among such methods, high order discontinuous Galerkin (DG) methods are particularly well-suited to the solution of time-dependent hyperbolic problems on modern computing architectures \cite{hesthaven2007, klockner2009}. The accuracy of high order methods can be attributed in part to their low numerical dissipation and dispersion compared to low order schemes \cite{ainsworth2004}. This accuracy has made them advantageous for the simulation of electro-magnetic and elastic wave propagation \cite{hesthaven2007, wilcox2010}.

Since the time-domain wave propagation is described by a hyperbolic system of partial differential equations, an explicit time integration can efficiently be applied. The finite-element methods, when coupled with an explicit time integrator, require the inversion of a global mass matrix, unless special techniques, such as diagonal mass lumping, are applied. Spectral element methods avoid the inversion of the global mass matrix for hexahedral elements by choosing nodal basis functions, which are discretely orthogonal with respect to an under-integrated $L^2$ inner product and result in a diagonal mass matrix. In contrast, high order DG methods produce block diagonal mass matrices, which are locally invertible. High order DG methods are often used for seismic simulation (elastic approximation) through the use of simplicial meshes \cite{kaser2007,de2007,ye2016}. 

DG methods impose inter-element continuity of approximate solutions between elements weakly through a numerical flux, of which several choices are common. de la Puente et al. \cite{de2008} solved the poroelastic wave equation using a local space-time DG method with a Rusanov type flux, which tends to be dissipative. In another study, Ward et al. \cite{ward2017}, derived an upwind flux by solving the exact Riemann problem on inter-element boundaries for the 2D isotropic strain-velocity form of the poroelastic wave equations. Recently, Zhan et al. \cite{zhan2019} also solved the poroelastic wave equation in 3D (formulated in stress-velocity form) using an upwind flux by solving the exact Riemann problem. The solution of Riemann problem requires diagonalization Jacobian matrices into polarized waves constituents, which is a computationally intensive process for the poroelastic system and does not extend naturally to anisotropic materials. Zhan et al. \cite{zhan2019} addresses the complexity of the diagonalization by reducing the rank of Jacobian matrices by using the method of generalized wave impedances.  
Ye at al. \cite{ye2016} completely avoid the process of diagonalization for the coupled acoustic-elastic wave equation by using a penalty flux based on natural boundary conditions.  In this study, we use a similar approach and derive an energy-stable penalty flux for the poroelastic wave equations. 

A poroelastic material is defined by the static acoustic properties of fluid, solid and frame exclusively, e.g. fifteen physical properties are required to define an orthotropic medium. The spatial scale of variability of these properties can range from macro (piece-wise constant approximations) to micro-heterogeneities (sub-element variations) in the medium. Most high order methods for wave equations on simplical meshes assume that material coefficients are constant over each element.  However, if the media is such that material gradients are non-zero in interior of the an element, piecewise constant approximations can yield inaccurate solutions \cite{chan2017I}. This can be circumvented by incorporating sub-element heterogeneities into weighted mass matrices, which recovers a highly accurate and energy-stable DG method \cite{mercerat2015}.  On tetrahedral meshes, this approach requires precomputation and storage of inverses for each local mass matrix, which increases the storage cost and  data movement at high orders of approximations. 
Chan et al. \cite{chan2018} circumvented these storage costs for elastic wave propagation by approximating weighted mass matrices with easily invertible ``weight-adjusted" approximations. We extend the same approach to the poroelastic wave equations, where matrix-valued weight functions arise after symmetrization of the system. 

In brief, the novelties of our approach are the following:
\begin{enumerate}
\item We obtain a consistent DG weak formulation for the poroelastic wave equations without the diagonalization into polarized wave constituents required for upwind fluxes. 
\item We introduce upwind-like dissipation through simple penalty terms in the numerical flux.
\item We implement an efficient ``weight-adjusted" approximate mass matrix to address micro-heterogeneities at the sub-element level present in poroelastic materials.
\item We prove stability and high order accuracy of the proposed DG method.
 \end{enumerate}    

The method is presented for tetrahedral meshes but extends naturally to quadrilateral and hexahedral meshes as well.  The outline of the paper is as follows: Section \ref{sec:intro} will discuss the system of equations describing the poroelastic wave equation along with interface conditions.  Section \ref{sec:ESdG} presents an energy stable formulation with simple penalty fluxes for the symmetric hyperbolic form of the poroelastic wave equations.  Section \ref{sec:ESdG} also discusses issues pertaining to storage and inversion of local mass matrices for material coefficients with micro-heterogeneities (sub-element variations).  Section \ref{sec:WADGForm} incorporates weight-adjusted approximations of weighted $L^2$ inner products and mass matrices with matrix-valued weights into DG discretizations of the poroelastic wave equations.  Finally, numerical results in Section \ref{sec:results} demonstrate the accuracy of this method for several problems in linear poroelasticity.

\section{System of equations describing porelastic waves}
\label{sec:intro}

Biot's theory describes wave propagation in a saturated porous medium i.e., a medium made up of a solid matrix (the skeleton or frame) fully saturated with single phase fluid.  Biot also assumed that an infinitesimal transformation relates the reference and current states of deformation; thus the displacement, strains, and particle velocity are small. The concept of infinitesimal deformation allows the use of the Lagrangian and Hamilton's principle to derive the equations governing the propagation of waves in such a medium.  In the following two subsections, we will review the constitutive equations, and equations of motions for a poroelastic medium. Readers are advised to refer to Biot's original papers \cite{biot1956I,biot1956II,biot1962} and \cite{carcione2014} for further detail.  

\subsection{Constitutive equations}
The constitutive equations for an inhomogeneous and anisotropic poroelastic medium are expressed as \cite{biot1956I,biot1956II, biot1962} as  
\begin{align}
\tau_{ij}&=c^{u}_{ijkl}\epsilon_{kl}-M\alpha_i \zeta \label{eq 1},\\
p&=M(\zeta-\alpha_k \epsilon_{kk}^{(m)}), \label{eq 2}
\end{align}
where $c^{u}_{ijkl}$ is the undrained stiffness constant, which can be expressed as $c^{u}_{ijkl}=c_{ijkl} + M \alpha_i \alpha_j$.  Here, $c_{ijkl}$ represents the elements of the stiffness tensor of the solid matrix, with $i,~j,~k,~l=1...6$ denoting indices of stress and strain tensors. The Biot effective stress coefficient $\alpha$ and Biot's parameter $M$ are derived experimentally from the bulk moduli of the solid and fluid~\cite{carcione2014}. Finally, $\smash[t]{\epsilon_{kl}^{(m)}}$ signifies the strain with $\zeta$ denoting the variation of fluid content.

\subsection{Dynamical equations and Darcy's law}

The dynamic equations describing the wave propagation in anisotropic  heterogeneous porous media are given by \cite{biot1956I,biot1956II,biot1962}
\begin{align}
\rho \frac{\partial\bm{v}}{\partial t} + \rho_f \frac{\partial\bm{q}}{\partial t}=\nabla \cdot \bm{\tau}_s, \label{eq 3} 
\end{align}
where $\rho=(1-\phi)\rho_s + \phi \rho_f$ is the bulk density, and $\rho_s$ and $\rho_f$ are the solid and fluid density, respectively. The variables $\bm{v}$ and $\bm{q}$ are the solid and fluid (relative to solid) particle velocities, respectively, and $\bm{\tau}_s$ is total stress (total) tensor. 

The generalized dynamic Darcy's law governing the fluid flow  in an anisotropic porous media, is expressed as
\begin{align}
\rho_f \frac{\partial{\bm{v}}}{{\partial t}}  + {\psi}_i(t) * \myfrac{\partial \bm{q}}{\partial t} =-\nabla p, \label{eq 4}
\end{align}  
where ${\psi}_i,~i=1...3$ are viscodynamic operators. Here, $``*"$ represents the convolution operation in time. In the low frequency range, i.e., for frequencies lower than $\omega_c=\min\left(\myfrac{\eta \phi }{\rho_f T_i \kappa_i}\right)$, $\psi_i$  can be expressed as 
 \begin{align}
 \psi_i(t)=m_i \delta(t) + \left(\eta/\kappa_i\right) H(t), \label{eq 5}
 \end{align}
 where $m_i=T_i \rho_f /\phi$, with $T_i$ being the tortuosity, $\eta$ the fluid viscosity, and $\kappa_i$s principal components of the global permeability tensor, while $\delta(t)$ is Dirac's function and $H(t)$ the Heaviside step function.
Substituting (\ref{eq 5}) into (\ref{eq 4}), we get 
\begin{align}
 \rho_f \frac{\partial{\bm{v}}}{{\partial t}} + {m_i} \frac{\partial \bm{q}}{\partial {t}} + \myfrac{\eta}{\kappa_i}\bm{q} =-\nabla p, \label{eq 6}
\end{align} 
The convolution operation in (\ref{eq 4}) represents the dissipation of energy due to fluid-flow induced by the wave motion.  As this work focuses on the low-frequency regime, the convolution is approximated by viscosity-dependent damping terms, shown in (\ref{eq 6}).

\subsection{The system of equations in matrix form}
To simplify notation, we introduce a matrix form of the system of equations while accounting for anisotropy and heterogeneity. Combining (\ref{eq 1}), (\ref{eq 2}), (\ref{eq 4}) and (\ref{eq 6}), we get
\begin{align} \label{eq 7}
\frac{\partial \bm{Q}}{\partial t} + \bm{A}\frac{\partial \bm{Q}}{\partial x_1} + \bm{B}\frac{\partial \bm{Q}}{\partial x_2} + \bm{C}\frac{\partial \bm{Q}}{\partial x_3}=\bm{D Q} + \bm{f},
\end{align}
where \begin{align*}
\bm{Q}=\left[\begin{array}{ccccccccccccc}
\tau_{11} ,& \tau_{22}, & \tau_{33}, & \tau_{23}, & \tau_{13}, & \tau_{12}, & p, & v_1, & v_2, & v_3, & q_1, & q_2, & q_3
\end{array}\right]^T,~
\bm{A}=\left[\begin{array}{c|c}
\bm{0}_{7 \times 7} & \bm{A}_{11}\\
\hline 
\bm{A}_{21} & \bm{0}_{6 \times 6}
\end{array}\right],
\end{align*}
\begin{align*}
\bm{B}=\left[\begin{array}{c|c}
\bm{0}_{7 \times 7} & \bm{B}_{11}\\
\hline 
\bm{B}_{21} & \bm{0}_{6 \times 6}
\end{array}\right],~\bm{C}=\left[\begin{array}{c|c}
\bm{0}_{7 \times 7} & \bm{C}_{11}\\
\hline 
\bm{C}_{21} & \bm{0}_{6 \times 6}
\end{array}\right],~\text{and}~\bm{D}=\left[\begin{array}{c|c}
\bm{0}_{7 \times 7} & \bm{0}_{6 \times 7} \\ 
\hline 
\bm{0}_{6 \times 7} & \bm{D}_{22}
\end{array}\right], ~\text{with}
\end{align*}
\begin{align*}
\bm{A}_{11}=-\left[\begin{array}{cccccc}
c_{11}^u & 0 & 0 & \alpha_1 M & 0 & 0\\
c_{12}^u & 0 & 0 & \alpha_2 M & 0 & 0\\
c_{13}^u & 0 & 0 & \alpha_3 M & 0 & 0\\
0 & 0 &0 & 0 & 0 & 0\\
0 & 0 & c_{55} & 0 & 0 & 0\\
0 & c_{66} & 0 & 0 & 0 & 0\\
-M \alpha_1 & 0 & 0 & - M & 0 & 0
\end{array}\right],~\bm{A}_{21}=-\left[\begin{array}{ccccccc}
\myfrac{m_1}{\beta_1} & 0 & 0 & 0 & 0 & 0 & \myfrac{\rho_f}{\beta_1}\\
0 & 0 & 0 & 0 & 0 & \myfrac{m_2}{\beta_2} & 0 \\
0 & 0 & 0 & 0 & \myfrac{m_3}{\beta_3} & 0 & 0\\
-\myfrac{\rho_f}{\beta_1} & 0 & 0 & 0 & 0 & 0 & -\myfrac{\rho}{\beta_1}\\
0 & 0 & 0 & 0 & 0 & -\myfrac{\rho_f}{\beta_2} & 0 \\
0 & 0 & 0 & 0 & -\myfrac{\rho_f}{\beta_3}&0&0     
\end{array}\right],
\end{align*}
\begin{align*}
\bm{B}_{11}=-\left[\begin{array}{cccccc}
0 & c_{12}^u & 0 & 0 & \alpha_1 M & 0 \\
0 & c_{22}^u & 0 & 0 & \alpha_2 M & 0 \\
0 & c_{33}^u & 0 & 0 & \alpha_3 M & 0 \\
0 & 0 & c_{44} & 0 & 0 &0 \\
0 & 0  & 0 & 0 & 0 & 0\\
c_{66} & 0 & 0 & 0 & 0 & 0\\
0 & - M\alpha_2 & 0 & 0 & -M & 0
\end{array}\right],~\bm{B}_{21}=-\left[\begin{array}{ccccccc}
0 & 0 & 0 & 0 & 0 & \myfrac{m_1}{\beta_1} & 0 \\
0 & \myfrac{m_2}{\beta_2} & 0 & 0 & 0 & 0 & \myfrac{\rho_f}{\beta_2} \\
0 & 0 & 0 & \myfrac{m_3}{\beta_3} & 0 & 0 & 0\\
0 & 0 & 0 & 0 & 0 & -\myfrac{\rho_f}{\beta_1} & 0\\
0 & -\myfrac{\rho_f}{\beta_2} & 0 & 0 & 0 & 0 & 0\\
0 & 0 & 0 & -\myfrac{\rho_f}{\beta_3} & 0 & 0 & 0 
\end{array}\right],
\end{align*}
\begin{align*}
\bm{C}_{11}=-\left[\begin{array}{cccccc}
0 & 0 & c_{13}^u & 0 & 0 & \alpha_1 M \\
0 & 0 & c_{23}^u & 0 & 0 & \alpha_2 M \\
0 & 0 & c_{33}^u & 0 & 0 & \alpha_3 M \\
0 & 0 & c_{44} & 0 & 0 & 0 \\
c_{55} & 0 & 0 & 0 & 0 & 0 \\
0 & 0 & 0 & 0 & 0 & 0 \\
0 & 0 & -M \alpha_3 & 0 & 0 & -M 
\end{array}\right],~\bm{C}_{21}=-\left[\begin{array}{ccccccc}
0 & 0 & 0 & 0 & \myfrac{m_1}{\beta_1} & 0 & 0 \\
0 & 0 & 0 & \myfrac{m_2}{\beta_2} & 0  & 0 & 0 \\
0 & 0 & \myfrac{m_3}{\beta_3} & 0 & 0 & 0 & \myfrac{\rho_f}{\beta_3} \\
0 & 0 & 0 & 0 & -\myfrac{\rho_f}{\beta_1}& 0 & 0 \\
0 & 0 & 0 & -\myfrac{\rho_f}{\beta_2} & 0 & 0 & 0 \\
0 & 0 & -\myfrac{\rho_f}{\beta_3} & 0 & 0  & 0 & -\myfrac{\rho}{\beta_3} 
\end{array}\right],~\text{and}
\end{align*}
\begin{align*}
\bm{D}_{22}=\left[\begin{array}{cccccc}
0 & 0 & 0 & \myfrac{\rho_f \eta}{\beta_1 \kappa_1} & 0 & 0\\
0 & 0 & 0 & 0 & \myfrac{\rho_f \eta}{\beta_2 \kappa_2} & 0\\
0 & 0 & 0 & 0 & 0 & \myfrac{\rho_f \eta}{\beta_3 \kappa_3}\\
0 & 0 & 0  & -\myfrac{\rho \eta}{\beta_1 \kappa_1} & 0 & 0\\
0 & 0 & 0 & 0 & -\myfrac{\rho\eta}{\beta_2 \kappa_2} & 0\\
0 & 0 & 0 & 0 & 0 & -\myfrac{\rho \eta}{\beta_3 \kappa_3}
\end{array}\right],
\end{align*}
where $\beta_i=\rho m_i - \rho_f^2$ and $\bm{f}$ is the vector of forcing terms.

\subsection{Symmetric form of system of poroelastic equations}
To ensure the stability and convergence of the numerical scheme, we utilize a symmetric form of (\ref{eq 7}). The symmetric form of (\ref{eq 7}) for velocity and stress is expressed as 
\begin{align}
{\bm {Q}_s} \myfrac{\partial \bm{\tau}}{\partial t}=\sum_{i=1}^{d}\bm{A}_i \myfrac{\partial \bm{v}} {\partial \bm{x}_i}, \label{eq8}\\
{\bm{Q}_v} \myfrac{\partial \bm{v}}{\partial t}=\sum_{i=1}^{d}\bm{A}_i^T \myfrac{\partial \bm{\tau}} {\partial \bm{x}_i} + \bm{D} \bm{v} + \bm{f}, \label{eq9}
\end{align}

where \begin{align*}
\bm{Q}_s=\left[\begin{array}{c|c}
\bm{S}& \bm{S \alpha} \\
\hline
\bm{\alpha^T} \bm{S} & \myfrac{1}{M} + \bm{\alpha}^T S \bm{\alpha}\\
\end{array}\right],~~
\bm{Q}_v=\left[\begin{array}{cccccc}
\rho&0&0&\rho_f&0&0\\
0&\rho&0&0&\rho_f&0\\
0&0&\rho&0&0&\rho_f\\
\rho_f&0&0&m_1&0&0\\
0&\rho_f&0&0&m_2&0\\
0&0&\rho_f&0&0&m_3\\
\end{array}\right],
\end{align*} with
\begin{align*}
\bm{S}&=\left[\begin{array}{cccccc}
\myfrac{c_{11}c_{33}-c_{13}^2 }{c_0} & \myfrac{c_{13}^2 - c_{12}c_{33}}{c_0}& \myfrac{- c_{13}}{c_1} & 0 & 0 & 0 \\
\myfrac{c_{13}^2 - c_{12}c_{33}}{c_0} & \myfrac{ c_{11}c_{33}- c_{13}^2}{c_0} & \myfrac{-c_{13}}{c_1} & 0 & 0 &0 \\
-\myfrac{c_{13}}{c_1} & -\myfrac{c_{13}}{c_1} & \myfrac{c_{11}+ c_{12}}{c_1} & 0 & 0 &0\\
0 & 0 & 0 & \myfrac{1}{c_{55}} & 0 & 0\\
0 & 0 & 0 & 0 & \myfrac{1}{c_{55}} & 0 \\
0 & 0 & 0 & 0 & 0&  \myfrac{2}{c_{11}-c_{12}} \\
\end{array}
\right],~~~\bm{\alpha}=\left[\begin{array}{c}
\alpha_1\\
\alpha_2\\
\alpha_3\\
0\\
0\\
0\\
\end{array}\right],
\end{align*} and 
$\bm{\tau}=[\tau_{11},~\tau_{22},~\tau_{33},~\tau_{23},~\tau_{13},~\tau_{12},~p]^T$, $\bm{v}=[v_1,~v_2,~v_3,~q_1,~q_2,~q_3]^T$, with
\begin{align*}
\bm{A}_1&=\left[ \begin{array}{cccccc}
1&0&0&0&0&0\\
0&0&0&0&0&0\\
0&0&0&0&0&0\\
0&0&0&0&0&0\\
0&0&1&0&0&0\\
0&1&0&0&0&0\\
0&0&0&-1&0&0\\
\end{array}\right], ~~
\bm{A}_2=\left[ \begin{array}{cccccc}
0&0&0&0&0&0\\
0&1&0&0&0&0\\
0&0&0&0&0&0\\
0&0&1&0&0&0\\
0&0&0&0&0&0\\
1&0&0&0&0&0\\
0&0&0&0&-1&0\\
\end{array}\right],~~\\
\bm{A}_3&=\left[ \begin{array}{cccccc}
0&0&0&0&0&0\\
0&0&0&0&0&0\\
0&0&1&0&0&0\\
0&1&0&0&0&0\\
1&0&0&0&0&0\\
0&0&0&0&0&0\\
0&0&0&0&0&-1\\
\end{array}\right],~
\bm{D}=\left[ \begin{array}{cccccc}
0&0&0&0&0&0\\
0&0&0&0&0&0\\
0&0&0&0&0&0\\
0&0&0&0&0&0\\
0&0&0&-\myfrac{\eta}{\kappa_1}&0&0\\
0&0&0&0&-\myfrac{\eta}{\kappa_2}&0\\
0&0&0&0&0&-\myfrac{\eta}{\kappa_3}\\
\end{array}\right].
\end{align*}

Explicit expressions for $\bm{Q}_s^{-1}$ and $\bm{Q}_{v}^{-1}$ are given by
\[
\bm{Q}_s^{-1}=\left[\begin{array}{ccccccc}
c_{11} + M\alpha_1^2 & c_{12} + M\alpha_1\alpha_2 & c_{13} + M\alpha_1\alpha_3 & 0 & 0 & 0 & -M\alpha_1 \\
c_{12} + M\alpha_1 \alpha_2 & c_{11} + M\alpha_2^2 & c_{13} + M\alpha_2\alpha_3 & 0 & 0 & 0 & -M\alpha_2\\
c_{13}+M\alpha_1 \alpha_3  & c_{13} + M\alpha_2 \alpha_3 & c_{33} + M\alpha_3^2  & 0 & 0 & 0 & -M\alpha_3 \\
0 & 0 & 0 & c_{55}& 0 & 0& 0\\
0 & 0 & 0 & 0 & c_{55}& 0 & 0\\
0 & 0 & 0 & 0 & 0 & \myfrac{c_{11} - c_{12}}{2} & 0\\
-M\alpha_1 & -M\alpha_2 & -M\alpha_3 & 0 & 0 & 0 & M \\
\end{array}
\right],
\]
\[\bm{Q}_v^{-1}=\left[\begin{array}{ccccccc}
\myfrac{m_1}{\rho m_1 -\rho_f^2} & 0 & 0 & \myfrac{-\rho_f}{\rho m_1 -\rho_f^2} & 0 & 0\\
0 & \myfrac{m_2}{\rho m_2  -\rho_f^2} & 0 & 0 & \myfrac{-\rho_f}{\rho m_2 -\rho_f^2} & 0\\
0 & 0 & \myfrac{m_3}{\rho m_3 - \rho_f^2} & 0 & 0 & \myfrac{-\rho_f}{\rho m_3 -\rho_f^2}\\
\myfrac{-\rho_f}{ \rho m_1 -\rho_f^2}&0&0 & \myfrac{\rho}{\rho m_1 -\rho_f^2}&0&0 \\
0 & \myfrac{-\rho_f}{\rho m_2-\rho_f^2}& 0 & 0 & \myfrac{\rho}{\rho m_2 -\rho_f^2} & 0 \\
0&0& \myfrac{-\rho_f}{\rho m_3 -\rho_f^2} & 0 & 0 & \myfrac{\rho}{\rho m_3 - \rho_f^2}\\
\end{array}
\right].
\]
The matrices $\bm{Q}_s$ and $\bm{Q}_v$ are Hessians of the potential and kinetic energy, respectively, computed with respect to state variables in $\bm{Q}$. Since the energy of the poroelastic system is quadratic positive-definite \cite{carcione1996}, the Hessians ($\bm{Q_s}$, $\bm{Q}_v$) are symmetric positive-definite. We also assume the Hessians are bounded as follows
\[
0 < s_{\text{min}} \le \bm{u}^T \bm{Q}_s(\bm{x})\bm{u} \le s_{\text{max}}  < \infty
\] 
\[
0 < \tilde{s}_{\text{min}} \le \bm{u}^T \bm{Q}_s^{-1}(\bm{x})\bm{u} \le \tilde{s}_{\text{max}}  < \infty
\]
\[
0 < v_{\text{min}} \le \bm{u}^T \bm{Q}_v(\bm{x})\bm{u} \le v_{\text{max}}  < \infty
\] 
\[
0 < \tilde{v}_{\text{min}} \le \bm{u}^T \bm{Q}_v^{-1}(\bm{x})\bm{u} \le \tilde{v}_{\text{max}}  < \infty
\]
for $\forall$ $\bm{x} \in \mathbb{R}^{d}$ and $\forall$ $\bm{u} \in \mathbb{R}^{N_d}$.
Note that in (\ref{eq9}), the dissipation matrix $\bm{D}$ is negative-definite.  

Using the expressions for $\bm{Q}_s$ and $\bm{Q}_v$, the kinetic $(K)$ and potential energy $(V)$ of the poroelastic system are expressed as
\begin{align}
K=\myfrac{1}{2}\bm{\tau}^T \bm{Q}_s \bm{\tau},~~ \label{eq 10}
V=\myfrac{1}{2}\bm{v}^T \bm{Q}_v\bm{v}.
\end{align}

\section{An energy stable discontinuous Galerkin formulation for poroelastic wave propgation}
\label{sec:ESdG}
Energy stable discontinuous Galerkin methods for elastic wave propagation have been constructed based on symmetric formulations of the elastodynamic equations \cite{chan2017I}, and it is straightforward to extend such discontinuous Galerkin formulations to the symmetric poroelastic system. We assume that the domain $\Omega$ is exactly triangulated by a mesh $\Omega_h$ which consists of elements $D^k$ which are images of a reference element $\hat{D}$ under the local affine mapping.
\[
\bm{x}^{k}=\Phi^k \widehat{\bm{x}},
\] 
where $\bm{x}^k=\{x^k, y^k\}$ for $d=2$ and $\bm{x}^k=\{x^k, y^k, z^k\}$ for $d=3$ denote the physical coordinates on $D^k$ and $\hat{\bm{x}}=\{\hat{x}, \hat{y}\}$ for $d=2$ and $\widehat{\bm{x}}=\{\widehat{x}, \widehat{y}, \widehat{z}\}$ for $d=3$ denote coordinates on the reference element. We denote the determinant of the Jacobian of $\Phi^k$ as $J$. 

Solutions over each element $D^k$ are approximated from a local approximation space $V_h(D^k)$, which is defined as composition of mapping $\Phi^k$ and reference approximation space $V_h(\widehat{D})$
\[
V_h(D^k)=\Phi^k \circ V_h(\widehat{D}).
\] 

Subsequently, the global approximation space $V_h(\Omega_h)$ is defined as 
\begin{align*}
V_h(\Omega_h)=\bigoplus_{k} V_h(D^k). 
\end{align*}
In this work, we will take $V_h(\widehat{D})=P^N(\widehat{D})$, with $P^N(\widehat{D})$ being the space of polynomials of total degree $N$ on the reference simplex. In two dimensions, $P^N$ on a triangle is 
\[
P^N(\widehat{D})=\{\widehat{x}^i\widehat{y}^j, 0 \le i +j \le N\},
\] 
and in three dimensions, $P^N$ on a tetrahedron is
\[
P^N(\widehat{D})=\{\widehat{x}^i\widehat{y}^j\widehat{x}^k, 0 \le i +j+k \le N\}.
\] 

The $L^2$ inner product and norm over $D^k$ is represented as
\begin{align*}
\LRp{\bm{g}, \bm{h}}=\int_{D^k} \bm{g} \cdot \bm{h}~d{\bm{x}} = \int_{\hat{D}} \bm{g} \cdot \bm{h} J~d\hat{\bm{x}}, \qquad ||\bm{g}||^2_{L^2{\Omega}} = (\bm{g}, \bm{g})_{L^2(D^k)},
\end{align*}
where $\bm{g}$ and $\bm{h}$ are real vector-valued functions.  Global $L^2$ inner products and squared norms are defined as the sum of local $L^2$ inner products and squared norms over each elements. The $L^2$ inner product and norm over the boundary $\partial D^k$ of an element are similarly defined as
\[
\left< \bm{u}, \bm{v}  \right>_{L^2(\partial D^k)}=\int_{\partial D^k} \bm{u} \cdot \bm{v}~d\bm{x} = \sum_{f \in \partial D^k} \int_{\hat{f}} \bm{u} \cdot \bm{v} J^f~d\hat{\bm{x}}, \qquad ||\bm{u}||^2_{L^2(\partial D^k)}=\left<\bm{u}, \bm{u}\right>, 
\]
where $J^f$ is the Jacobian of the mapping from a reference face $\hat{f}$ to a physical face $f$ of an element.

Let $f$ be a face of an element $D^k$ with neighboring element $D^{k,+}$ and unit outward normal $\bm{n}$. Let $u$ be a function with discontinuities across element interfaces. We define the interior value $u^-$ and exterior value $u^+$ on face $f$ of $D^k$
\[
u^- = u|_{f \cap \partial D^k}, \qquad u^+ = u|_{f\cap \partial D^{k,+}}.
\]
The jump and average of a scalar function $u \in V_h(\Omega_h)$ over $f$ are then defined as
\[
\jump{u}=u^+ - u^-, \qquad \avg{u}=\myfrac{u^+ +  u^-}{2}.
\]
Jumps and averages of vector-valued functions $\bm{u} \in \mathbb{R}^{m}$ and and matrix-valued functions $\tilde{\bm{S}} \in \mathbb{R}^{m \times n}$  are defined component-wise.
\[
\left( \jump{\bm{u}}\right)_i = \jump{\bm{u}_i}, \qquad~1 \le i \le m\qquad~\left(\jump{\tilde{\bm{S}}}\right)_{ij}=\jump{\tilde{\bm{S}}}
\]

We can now specify a DG formulation for poroelastic wave equation (\ref{eq 7}). 
The symmetric hyperbolic system in (\ref{eq 7}) readily admits a DG formulation based on a penalty flux \cite{chan2017II, chan2017I}. For the symmetric first order poroelastic wave equation, the DG formulation in strong form  is given as
\begin{equation}
\begin{aligned}
\label{eq17}
\sum_{D^k \in \Omega_h} \left( \bm{Q}_s \myfrac{\partial \bm{\tau}}{\partial t} , \bm{h} \right)_{L^2(D^k)}=&\sum_{{D^k \in \Omega_h}} \left(\left(\sum_{i=1}^{d}\bm{A}_i \myfrac{\partial \bm{v}}{\partial \bm{x}_i} , \bm{h}   \right)_{L^2(D^k)} + \left \langle \myfrac{1}{2} \bm{A_n}\jump{\bm{v}} + \myfrac{\alpha_{\bm{\tau}}}{2} \bm{A_n}\bm{A}_n^T \jump{\bm{\tau}}, \bm{h}  \right \rangle_{L^2(\partial D^k)}\right) \\
\sum_{D^k \in \Omega_h} \left( \bm{Q}_v \myfrac{\partial \bm{v}}{\partial t} , \bm{g} \right)_{L^2(D^k)}=&\sum_{{D^k \in \Omega_h}} \Biggl(\left(\sum_{i=1}^{d}{\bm{A}_i}^T \myfrac{\partial \bm{\tau}}{\partial \bm{x}_i} + \bm{f} ,\bm{g}   \right)_{L^2(D^k)} \\  & + \left \langle \myfrac{1}{2} \bm{A_n}^T \jump{\bm{\tau}} + \myfrac{\alpha_{\bm{v}}}{2} \bm{A}_n^T \bm{A_n} \jump{\bm{v}}, \bm{g}  \right \rangle_{L^2(\partial D^k)}  +  \left(\bm{D}\bm{v}, \bm{g}\right)\Biggr),
\end{aligned}
\end{equation}

for all $\bm{h},~\bm{g} \in V_h(\Omega_h)$.  Here, $\bm{A}_n$ is the normal matrix defined on a face $f$ of an element 
\[
\bm{A}_n=\sum_{i=1}^d n_i \bm{A}_i=\left[\begin{array}{cccccc}
n_x &  0 & 0 & 0 & 0 & 0\\
0 & n_y & 0 & 0 & 0 & 0\\
0 & 0 & n_z & 0 & 0 & 0\\
0 & n_z & n_y & 0 & 0 & 0\\
n_z & 0 & n_x & 0 & 0 & 0\\
n_y & n_x & 0 & 0 & 0 & 0\\
0 & 0 & 0 & -n_x & -n_y & -n_z\\
\end{array}\right].
\]

The factors $\alpha_{\bm{\tau}}$ and $\alpha_{\bm{v}}$ are penalty parameters and defined on element interfaces. We assume that $\alpha_{\bm{\tau}},~\alpha_{\bm{v}} \ge 0$ and are piecewise constant over each shared face between two elements. These penalty constants can be taken to be zero, which results in a non-dissipative central flux, while $\alpha_{\bm{\tau}},~\alpha_{\bm{v}} > 0$ results in energy dissipation similar to the upwind flux \cite{hesthaven2007}. The stability of DG formulations are independent of the magnitude of these penalty parameters. However, a naive choice of these parameters will result in a stiffer semi-discrete system of ODEs and necessitates a smaller time under explicit time integration schemes.  In this work, we take $\alpha_\tau, \alpha_v = O(1)$ unless stated otherwise.

\begin{remark}
For several problems in Section~\ref{sec:num}, $\bm{f}$ is taken as a scaled point source perturbation or Dirac delta function $\bm{f(x)}=\bm{\beta}(\bm{x})\Delta(\bm{x - x_0})$ (with $\bm{\beta(x)} \in \mathbb{R}$).  In this setting, $\bm{f}$ is not $L^2$ integrable and thus $\left(\bm{f}, \bm{g}\right)_{L^2(D^k)}$ may not be well-defined.  In such cases, we commit a variation crime and evaluate its contribution as
\[
\sum_{k} \left(\bm{\beta}(\bm{x})\delta(\bm{x}-\bm{x}_0), \bm{g}\right)_{L^2(\left(D^k\right)}=\int_\Omega \bm{g}\cdot \bm{\beta}\delta(\bm{x - x_0}) d\bm{x} = \bm{g(\bm{x}_0)} \cdot \bm{\beta}({\bm{x}_0})
\]
\end{remark}

\subsection{Physical interpretation of central flux terms}
Apart from providing the global solution of the system, the central flux terms corresponding to the velocity $(\bm{A}_n \jump{\bm{v}})$ and the stress $(\bm{A}_n^T \jump{\bm{\tau} })$ also enforce natural boundary conditions at the interface between two poroelastic media. These conditions are related to the phenomena of reflection, refraction, and diffraction of waves in the presence of inhomogeneities and interfaces.  The component-wise expressions of the central flux terms $\bm{A}_n^T \jump{\bm{\tau} }$ and $\bm{A}_n \jump{\bm{v}}$ correspond to
\begin{alignat}{2}
\bm{\tau_s}^+ \cdot \bm{n}-\bm{\tau_s}^- \cdot \bm{n}&=0, \label{eq11}\\
\bm{v}^+ - \bm{v}^-&=0, \label{eq12}, \\
\bm{q}^+ \cdot \bm{n} -\bm{q}^- \cdot \bm{n}&=0 \label{eq13}, \\
{p^+ -p^-}&=0, \label{eq14}
\end{alignat}
where the ``$\pm$" convention is determined by the outward interface normal, $\bm{n}$.  In this work, the outer normal vector is assumed point in the direction of the ``$+$'' side of the interface.

The physical significance of the above conditions are:
\begin{enumerate}
\item (\ref{eq11}) describes the continuity of traction across the interface.
\item (\ref{eq12}) states that the medium across the interface stays intact.
\item (\ref{eq13}) implies that all fluid entering the interface should exit the other side. 
\item (\ref{eq14}) assumes perfect hydraulic conductivity across the interface, providing continuity of pore-fluid pressure.
\end{enumerate}
These interface conditions are also consistent with the open-pore boundary conditions of~\cite{deresiewicz1963, gurevich1999} between two different poroelastic media.

\subsection{Boundary conditions}
In many applications, the boundary condition for the top surface of a domain is assumed to be a free surface (stress free), with the remaining surfaces taken to be absorbing boundaries.
We impose boundary conditions on the DG formulation by choosing appropriate exterior values which result in modified boundary numerical fluxes.  Boundary conditions on the normal component of the stress can be imposed by noting that that the numerical flux contain the term $\jump{\bm{A}_n^T \bm{\tau}}=\jump{\widetilde{\bm{S}} \bm{n}}$ with 
\[
\bm{\widetilde{S}}=\left[\begin{array}{ccccccc}
\tau_{11} & \tau_{12} & \tau_{13} & 0\\
\tau_{21} & \tau_{22} & \tau_{23} & 0\\
\tau_{31} & \tau_{32} & \tau_{33} & 0\\
0 & 0 & 0 & p
\end{array}
\right].
\]
For a face which lies on the top surface of the domain, the free surface boundary or zero traction boundary conditions can be imposed by setting
\[
\jump{\bm{A}_n^T \bm{\tau}}=\jump{\widetilde{\bm{S}} \bm{n}}=-2\widetilde{\bm{S}}^- \bm{n}=-2\bm{A}_n^T \bm{\tau}^-,~~\bm{v}^+=\bm{v}^-\implies \jump{\bm{v}}=0.
\]
For problems which require the truncation of infinite or large domains, basic absorbing boundary conditions can be imposed by setting 
\[
\jump{\bm{A}_n^T \bm{\tau}}=\jump{\widetilde{\bm{S}} \bm{n}}=-\widetilde{\bm{S}}^- \bm{n}=-\bm{A}_n^T \bm{\tau}^-,~~\bm{v}^+=0 \implies\jump{\bm{v}}=-\bm{v}^-.
\]
In addition to the above boundary conditions, more accurate absorbing boundary conditions can be also imposed using perfectly matching layers (PML) \cite{zeng2001} or high order absorbing boundary conditions (HABC) \cite{hagstrom2004,modave2017}.  However, the implementation of such boundary conditions results in an augmented system of PDEs which can become very expensive.  For example, the implementation of PML for the poroelastic wave equations results in a system of thirty PDEs \cite{zeng2001}.

In all cases, the boundary conditions are imposed by computing the numerical fluxes based on the modified jumps. This imposition guarantees energy stability for free surface and absorbing boundary conditions.  

\subsection{Energy stability}
One can show that the DG formulation in (\ref{eq17}) is energy stable in the absence of external forces $(\bm{f}=0)$, and free-surface and absorbing boundary conditions. Integrating by parts the velocity equation in (\ref{eq17}) gives
\begin{equation}
\begin{aligned}
\label{eq18}
\sum_{D^k \in \Omega_h} \left( \bm{Q}_s \myfrac{\partial \bm{\tau}}{\partial t} , \bm{h} \right)_{L^2(D^k)}=&\sum_{{D^k \in \Omega_h}} \left(\left(\sum_{i=1}^{d}\bm{A}_i \myfrac{\partial \bm{v}}{\partial \bm{x}_i} , \bm{h}   \right)_{L^2(D^k)} + \left \langle \myfrac{1}{2} \bm{A_n}\jump{\bm{v}} + \myfrac{\alpha_{\bm{\tau}}}{2} \bm{A_n}\bm{A}_n^T \jump{\bm{\tau}}, \bm{h}  \right \rangle_{L^2(\partial D^k)}\right) \\
\sum_{D^k \in \Omega_h} \left( \bm{Q}_v \myfrac{\partial \bm{v}}{\partial t} , \bm{g} \right)_{L^2(D^k)}=&\sum_{{D^k \in \Omega_h}} \Biggl(-\left(\sum_{i=1}^{d} \bm{\tau},{\bm{A}_i} \myfrac{\partial \bm{g}}{\partial \bm{x}_i}\right)_{L^2(D^k)} \\ & + \left \langle  \bm{A_n}^T \avg{\bm{\tau}} + \myfrac{\alpha_{\bm{v}}}{2} \bm{A}_n^T \bm{A_n} \jump{\bm{v}}, \bm{g}  \right. \rangle_{L^2(\partial D^k)}  
+ \left(\bm{D}\bm{v}, \bm{v}\right)_{L^2(D^k)}\Biggr)
\end{aligned}
\end{equation}
Taking $\left(\bm{h}, \bm{g}\right)=(\bm{\tau}, \bm{v})$ and adding both equations together yields
\begin{align*}
\sum_{D^k \in \Omega_h} \myfrac{1}{2}&\myfrac{\partial}{\partial t}\left( (\bm{Q}_s \bm{\tau}, \bm{\tau})_{L^2(D^k)} + (\bm{Q}_v \bm{v}, \bm{v})_{L^2(D^k)}\right) \\
=&\sum_{D^k \in \Omega_h} \left \langle \myfrac{1}{2} \bm{A_n}\jump{\bm{v}} + \myfrac{\alpha_{\bm{\tau}}}{2} \bm{A_n}\bm{A}_n^T \jump{\bm{\tau}}, \bm{h}  \right \rangle_{L^2(\partial D^k)}+ \left \langle  \bm{A_n}^T \avg{\bm{\tau}} + \myfrac{\alpha_{\bm{v}}}{2} \bm{A}_n^T \bm{A_n} \jump{\bm{v}}, \bm{g}  \right \rangle_{L^2(\partial D^k)}\\ & + \left(\bm{D}\bm{v}, \bm{v}\right)_{L^2(D^k)}\\
=&\sum_{D^k \in \Omega_h}  \sum_{f \in \partial D^k}  \int_f\left(\myfrac{1}{2} \bm{\tau}^T \bm{A}_n \jump{\bm{v}} + \myfrac{\alpha_{\bm{\tau}}}{2} \bm{\tau}^T \bm{A_n} \bm{A}_n^T\jump{\tau}+\bm{v}^T \bm{A}_n^T\avg{\bm{\tau}} + \myfrac{\alpha_{\bm{v}}}{2} \bm{v}^T \bm{A}_n^T \bm{A}_n \jump{\bm{v}}\right) \text{d}\bm{x} \\ \nonumber
& + \sum_{D^k \in \Omega_h} \int_{D^k} \bm{v}^T \bm{D} \bm{v}~\text{d}\bm{x},
\end{align*}
where the term 
\[\sum_{D^k \in \Omega_h} \myfrac{1}{2} \myfrac{\partial}{\partial t}\left( (\bm{Q}_s \bm{\tau}, \bm{\tau})_{L^2(D^k)} + (\bm{Q}_v \bm{v}, \bm{v})_{L^2(D^k)}\right)\]
is the total energy of the system. Let $\Gamma_h$ denote the set of unique faces in $\Omega_h$ and let $\Gamma_{\bm{\tau}}$, $\Gamma_{\text{abc}}$ denote boundaries where free-surface and absorbing boundary conditions are imposed, respectively. We separate surface terms into contributions from interior shared faces and from boundary faces. On interior shared faces, we sum the contributions from the two adjacent elements to yield
\begin{align*}
&\sum_{f \in \Gamma_h \setminus \partial \Omega} \int_f\left(\myfrac{1}{2} \bm{\tau}^T \bm{A}_n \jump{\bm{v}} + \myfrac{\alpha_{\bm{\tau}}}{2} \bm{\tau}^T \bm{A_n} \bm{A}_n^T\jump{\tau}+\bm{v}^T \bm{A}_n^T\avg{\bm{\tau}} + \myfrac{\alpha_{\bm{v}}}{2} \bm{v}^T \bm{A}_n^T \bm{A}_n \jump{\bm{v}}\right) \text{d}\bm{x} + \sum_{D^k \in \Omega_h} \int_{D^k} \bm{v}^T \bm{D} \bm{v}~\text{d}\bm{x}\\
&=-\sum_{f \in \Gamma_h \setminus \partial \Omega} \int_f \left( \myfrac{\alpha_{\bm{\tau}}}{2} \left |\bm{A}_n^T \jump{\bm{\tau}}\right|^2 + \myfrac{\alpha_{\bm{v}}}{2} \left|\bm{A}_n \jump{\bm{v}} \right|^2 \right)\text{d}\bm{x} + \sum_{D^k \in \Omega_h} \int_{D^k} \bm{v}^T \bm{D} \bm{v}~\text{d}\bm{x},
\end{align*}
where $\bm{v}^T \bm{D} \bm{v} < 0 $, since $\bm{D}$ is a negative semi-definite matrix.
For faces which lie on $\Gamma_{\bm{\tau}}$ , $\bm{A}_n^T=-2\bm{A}_{n}^T \bm{\tau}^-,~\bm{A}_n^T \avg{\tau}=0$ and $\jump{\bm{v}}=0$ yielding 
\begin{align*}
&\sum_{f \in \Gamma_{\bm{\tau}}} \int_f\left(\myfrac{1}{2} \bm{\tau}^T \bm{A}_n \jump{\bm{v}} + \myfrac{\alpha_{\bm{\tau}}}{2} \bm{\tau}^T \bm{A_n} \bm{A}_n^T\jump{\tau}+\bm{v}^T \bm{A}_n^T\avg{\bm{\tau}} + \myfrac{\alpha_{\bm{v}}}{2} \bm{v}^T \bm{A}_n^T \bm{A}_n \jump{\bm{v}}\right) \text{d}\bm{x} \\
&=-\sum_{f \in \Gamma_{\bm{\tau}}} \int_f \left ( \alpha_{\bm{\tau}} \left| \bm{A}_n^T \bm{\tau^-} \right|^2\right)~\text{d} \bm{x}.
\end{align*}
Finally, for faces in $\Gamma_{\text{abc}}$, we have $\bm{A}_n^T \avg{\tau}=\myfrac{1}{2} \bm{A}_n^T \tau{^-}$, $\bm{A}_n^T\jump{\tau}=-\bm{A}_n^T \bm{\tau}^- $ and $\jump{\bm{v}}=-\bm{v^-}$, yielding
\begin{align*}
&\sum_{f \in \Gamma_h \setminus \partial \Omega} \int_f\left(\myfrac{1}{2} \bm{\tau}^T \bm{A}_n \jump{\bm{v}} + \myfrac{\alpha_{\bm{\tau}}}{2} \bm{\tau}^T \bm{A_n} \bm{A}_n^T\jump{\tau}+\bm{v}^T \bm{A}_n^T\avg{\bm{\tau}} + \myfrac{\alpha_{\bm{v}}}{2} \bm{v}^T \bm{A}_n^T \bm{A}_n \jump{\bm{v}}\right) \text{d}\bm{x}\\
&=-\sum_{f \in \Gamma_{\text{abc}}} \int_f \left( \myfrac{\alpha_{\bm{\tau}}}{2} \left |\bm{A}_n^T \bm{\tau}^-\right|^2 + \myfrac{\alpha_{\bm{v}}}{2} \left|\bm{A}_n \bm{v^-} \right|^2 \right)\text{d}\bm{x} ,
\end{align*}
Combining contributions from all faces and dissipation in the system yields the following result:
\begin{theorem}
The DG formulation in (\ref{eq17}) is energy stable for $\alpha_{\bm{\tau}}, \alpha_{\bm{v}} \geq 0$ such that
\begin{align}
\label{eq:th1}
\sum_{D^k \in \Omega_h} & \myfrac{1}{2} \myfrac{\partial}{\partial t}\left( (\bm{Q}_s \bm{\tau}, \bm{\tau})_{L^2(D^k)} + (\bm{Q}_v \bm{v}, \bm{v})_{L^2(D^k)}\right)=-\sum_{f \in \Gamma_h \setminus \partial \Omega} \int_f \left( \myfrac{\alpha_{\bm{\tau}}}{2} \left |\bm{A}_n^T \jump{\bm{\tau}}\right|^2 + \myfrac{\alpha_{\bm{v}}}{2} \left|\bm{A}_n \jump{\bm{v}} \right|^2 \right)\text{d}\bm{x} \nonumber \\
&-\sum_{f \in \Gamma_{\bm{\tau}}} \int_f \left ( \alpha_{\bm{\tau}} \left| \bm{A}_n^T \bm{\tau^-} \right|^2\right)~\text{d} \bm{x}-\sum_{f \in \Gamma_{\text{abc}}} \int_f \left( \myfrac{\alpha_{\bm{\tau}}}{2} \left |\bm{A}_n^T \bm{\tau}^-\right|^2 + \myfrac{\alpha_{\bm{v}}}{2} \left|\bm{A}_n \bm{v^-} \right|^2 \right)\text{d}\bm{x} \nonumber \\ & + \sum_{D^k \in \Omega_h} \int_{D^k} \bm{v}^T \bm{D} \bm{v}~\text{d}\bm{x} \le 0. 
\end{align}
\end{theorem}

Since $\bm{Q}_s$ and $\bm{Q}_v$ are positive definite, the left hand side of (\ref{eq:th1}) is an $L^2$-equivalent norm on $(\bm{\tau}, \bm{v})$ and Theorem 1 implies that magnitude of the DG solution is non-increasing in time. This also shows that dissipation is present for positive penalization parameters, i.e. ${\alpha_{\bm{\tau}}}, {\alpha_{\bm{v}}} > 0$.

\subsection{The semi-discrete matrix system for DG}

Let ${\{\phi_i\}}_{i=1}^{N_p}$ be a basis for $P^N\left(\widehat{D}\right)$. In our implementation, we use nodal basis functions located at Warp and Blend interpolation points \cite{hesthaven2007}, which are defined implicitly using an orthogonal polynomial basis on the reference simplex. We define the reference mass matrix $\widehat{\bm{M}}$ and the physical mass matrix $\bm{M}$ for an element $D^k$ as
\[
\left(\widehat{\bm{M}}\right)_{ij}=\int_{\widehat{D}} \phi_j \phi_i~\text{d}\bm{x}, \qquad (\bm{M}_{ij})=\int_{D^k} \phi_j \phi_i~\text{d}\bm{x}=\int_{\widehat{D}} \phi_j \phi_i J~\text{d}\bm{\widehat{x}}.
\]

For affine mappings, $J$ is constant and $\bm{M}=J\bm{\widehat{M}}$. We also define weak differentiation matrices $\bm{S}_{k}$ and face mass matrices $\bm{M}_f$ such that
\[
\left(\bm{S}_{k}\right)_{ij}=\int_{D^k} \myfrac{\partial \phi_j}{\partial \bm{x}_k} \phi_i~\text{d} \bm{x}, ~~~ (\bm{M}_f)_{ij}=\int_f \phi_j \phi_i \text{d}\bm{x}=\int_{\hat{f}} \phi_j \phi_i J^f \text{d}\hat{\bm{x}},
\]
where $J^f$ is the Jacobian of the mapping from the reference face $\widehat{f}$ to $f$. For affinely mapped simplices, $J^f$ is also constant and $M_f=J^f\widehat{\bm{M}}_f$, where the definition of the reference face mass matrix $\widehat{\bm{M}}_f$ is analogous to the definition of the reference mass matrix $\widehat{\bm{M}}$.

Finally, we introduce weighted mass matrices.  Let $w(\bm{x}) \in \mathbb{R}$ and $\bm{W(x)} \in \mathbb{R}^{m \times n}$. Then, scalar and matrix-weighted mass matrices $\bm{M}_w$ and $\bm{M_W}$ are defined as 
\begin{align}
(\bm{M}_w)_{ij}=\int_{D^k} w(\bm{x}) \phi_j(\bm{x}) \phi_i(\bm{x})~ \text{d}\bm{x},\qquad\bm{M_W}=\begin{bmatrix}\bm{M_{W_{1,1}}} & \dots & \bm{M_{W_{1,n}}}\\
\vdots &\ddots & \vdots\\
\bm{M_{W_{m,1}}} & \dots & \bm{M_{W_{m,n}}},
\end{bmatrix}
\label{eq:matweight}
\end{align}
where $\bm{M_{W_{i,j}}}$ is the scalar weighted mass matrix weighted by the $(i,j)^{\text{th}}$ element of $\bm{W}$.  Note that $\bm{M}_w, \bm{M}_{\bm{W}}$ are positive definite if $w(x), \bm{W}$ are pointwise positive definite.

Local contributions to the DG variational form may be evaluated in a quadrature-free manner using matrix-weighted mass matrices as defined above. Let $\bm{\Sigma_i},~\bm{V}_i$ denote vectors containing degrees of freedom for solutions components $\bm{\tau}_i$ and $\bm{v}_i$, such that
\begin{align*}
\bm{v}_i (\bm{x}, t)&=\sum_{j=1}^{N_p} (\bm{V}_i(t))_j\phi_j(\bm{x}), \qquad\qquad1 \le i \le 6 \\
\bm{\tau}_i (\bm{x}, t)&=\sum_{j=1}^{N_p} (\bm{\Sigma_i(t)})_j \phi_j(\bm{x}), \qquad\qquad1 \le i \le 7 
\end{align*}
Then, the local DG formulation can be written as a block system of ordinary differential equations by concatenating $\bm{\Sigma}_i,~\bm{V}_i$ into single vectors $\bm{\Sigma}$ and $\bm{V}$ and using the Kronecker product $\otimes$
\begin{align} 
 \bm{M}_{\bm{Q}_s} \myfrac{\partial \bm{\Sigma}}{\partial t} &= \sum_{i=1}^{d} \left(\bm{A}_i \otimes \bm{S_i}\right)\bm{V} + \sum_{f \in \partial D^k} \left(\bm{I} \otimes \bm{M}_f\right) \bm{F}_\tau \label{eq21}\\
\bm{M}_{\bm{Q}_v} \myfrac{\partial \bm{V}}{\partial t}&=\sum_{i=1}^{d}\left({\bm{A}_i}^T \otimes \bm{S_i}\right)\bm{\Sigma} + \sum_{f \in \partial D^k} \left(\bm{I} \otimes \bm{M}_f\right) \bm{F}_v + \bm{M}_{\bm{D}}\bm{V} \label{eq22},
\end{align}
where $\bm{F}_v$ and $\bm{F}_\tau$ denote the degrees of freedom for the velocity and stress numerical fluxes.

In order to apply a time integrator, we must invert $\bm{M}_{\bm{Q}_s}$ and $\bm{M}_{\bm{Q}_v}$. While the inversion of $\bm{M}_{\bm{Q}_s}$ and $\bm{M}_{\bm{Q}_v}$ can be parallelized from element to element, doing so typically requires either the precomputation and storage of the dense matrix inverses or on-the-fly construction and solution of a large dense matrix system at every time step. The former option requires a large amount of storage, while the latter option is computationally expensive and difficult to parallelize among degrees of freedom. This cost can be avoided when ${\bm{Q}_s}$ and $\bm{Q}_v$ are constant over an element $D^k$. In this case, $\bm{M}_{\bm{Q}_s}$ reduces to
\begin{align*} 
\bm{M}_{\bm{Q}_s}^{-1}=\begin{bmatrix}{\bm{Q}_s}_{(1,1)}\bm{M} & \dots & {\bm{Q}_s}_{(1,N_d)}\bm{M}\\
\vdots &\ddots & \vdots\\
{\bm{Q}_s}_{(N_d,1)}\bm{M} & \dots & {\bm{Q}_s}_{(N_d,N_d)}\bm{M}
\end{bmatrix}^{-1} =\left({\bm{Q}_s} \otimes \bm{M}  \right)^{-1}={\bm{Q}_s}^{-1} \otimes \left(\myfrac{1}{J} \bm{\widehat{M}^{-1}}\right).
\end{align*} 
Similarly ${\bm{M}^{-1}_{\bm{Q}_v}}$ and ${\bm{M}^{-1}_{\bm{D}}}$ can be expressed as $\bm{M}^{-1}_{\bm{Q}_v}={\bm{Q}_v}^{-1} \otimes \left(\myfrac{1}{J} \bm{\widehat{M}^{-1}}\right)$, and $\bm{M}^{-1}_{\bm{D}}={\bm{D}}^{-1} \otimes \left(\myfrac{1}{J} \bm{\widehat{M}^{-1}}\right)$ respectively.  Applying these observations to (\ref{eq21}) and (\ref{eq22}) yields the following sets of local ODEs over each element
\begin{align}
 \myfrac{\partial \bm{\Sigma}}{\partial t} &= \sum_{i=1}^{d} \left({\bm{Q}_s}^{-1}\bm{A} \otimes \bm{D_i}\right)\bm{V} + \sum_{f \in \partial D^k} \left({\bm{Q}_s}^{-1} \otimes \bm{M}_f\right) \bm{F}_\tau,\\
 \myfrac{\partial \bm{V}}{\partial t}&=\sum_{i=1}^{d}\left({\bm{Q}_v} ^{-1}{\bm{A}}^T \otimes \bm{D_i}\right)\bm{\Sigma} + \sum_{f \in \partial D^k} \left(\bm{I} \otimes \bm{M}_f\right) \bm{F}_v + {\bm{M}^{-1}_{\bm{Q}_v}}\bm{M}_{\bm{D}}\bm{V},
\end{align}
where we have introduced the differentiation matrix $\bm{D}_i=\bm{M}{^{-1}}\bm{S}_i$ and lift matrix $\bm{L}_f=\bm{M}^{-1}\bm{M}_f$. For affine elements, both derivative and lift matrices are applied using products of geometric factors and reference derivative and lift matrices.

Unfortunately, if $\bm{Q}_s$ and $\bm{Q}_v$ varies spatially with in the element, then above approach can no longer be used to invert $\bm{Q}_s$ and $\bm{Q}_v$.  Here, we follow the approach of \cite{chan2018}, where $\bm{M}_{\bm{Q}_s},\bm{M}_{\bm{Q}_v}$ are replaced with  weight-adjusted approximations. These approximations are low storage, simple to invert, and yield an energy stable and high order accurate DG method to approximate the matrix-weighted $L^2$ inner product (and corresponding matrix-weighted mass-matrices $\bm{Q}_s$ and $\bm{Q}_v$).  

We also note that, material coefficients $\bm{Q_s},~\bm{Q}_v$ appear only the left hand side of (\ref{eq 7}). The right hand side of (\ref{eq 7}) is equivalent to the discretization of a constant coefficient system.  This provides additional advantages in that the right hand side can be evaluated using efficient techniques for DG discretizations of constant-coefficient problems \cite{chan2017gpu, guo2018}.

\section{Weight-adjusted discontinuous Galerkin (WADG) formulation for poroelastic wave propagation}
\label{sec:WADGForm}
We wish to apply weight-adjusted approximations to avoid the inversion of $\bm{M}_{\bm{Q}{_s}}$ and $\bm{M}_{\bm{Q}{_v}}$. Weight-adjusted inner products are high order approximations of weighted $L^2$ inner products. These weight-adjusted inner products result in weight-adjusted mass matrices whose inverses approximate the inverses of weighted $L^2$ mass matrices.  We briefly review scalar and matrix-valued weight-adjusted mass matrices in the following section.  Matrix weight-adjusted inner products and weight-adjusted approximations with matrix weights are discussed in more detail in \ref{appendix:A0} and \ref{appendix:A02}.

\subsection{Approximation of weighted mass matrix inverses}
The advantage of weighted-adjusted inner products is that the corresponding weight-adjusted mass matrices are straightforward to invert. For scalar weights, weight-adjusted mass matrices are given as follows \cite{chan2018}
\[
\bm{M}_w\approx \bm{M} \bm{M_{1/w}^{-1}} \bm{M}, \qquad\bm{M}_w^{-1} \approx \bm{M}^{-1} \bm{M}_{1/w} \bm{M}^{-1}.
\] 
By evaluating $\bm{M}_{1/w}$ in a matrix-free fashion using a sufficiently accurate quadrature rule, the inverse of the weight-adjusted mass matrix $\bm{M}^{-1} \bm{M}_{1/w} \bm{M}^{-1}$ yields a low storage implementation. Let $\bm{\hat{x}}_i, \bm{\hat{w}}_i$ denote quadrature points and weights on the reference element, and let $\bm{V}_q$ denote the generalized Vandermonde matrix
\begin{align*}
(\bm{V}_q)_{ij}=\phi_j(\bm{\hat{x}}_i)
\end{align*}
whose columns correspond to evaluations of basis functions at quadrature points. Then, for affine elements, $M=J\widehat{\bm{M}}=J\bm{V}_q^T \text{diag}(\hat{\bm{w}}_i)\bm{V}_q$, where $\widehat{\bm{M}}$ is the reference mass matrix and $J$ is the determinant of the Jacobian of the reference-to-physical mapping, which is constant for affine mappings. Additionally, 
\begin{align*}
\bm{M}_{1/w}=J\bm{V}{_q}^T\text{diag}\left(\hat{\bm{w}}_i/w(\hat{\bm{x}_i})\right) \bm{V}_q ,
\end{align*}
where $w(\hat{\bm{x}}_i)$ denotes the evaluation of the weight function $w(\bm{x})$ at quadrature points. Thus, for a vector $\bm{u}$, the inverse of the weight-adjusted mass matrix can be applied as follows
\[
\bm{M}^{-1} \bm{M}_{1/w} \bm{M}^{-1} \bm{u}=\bm{M}^{-1} \bm{V}_q^T\text{diag}(\hat{\bm{w}}_i)\text{diag}\left(\myfrac{1}{w(\bm{x}_i)}\right) \bm{V}_q \bm{M}^{-1}\bm{u}=\bm{P}_q\text{diag}\left(\myfrac{1}{w(\widehat{\bm{x}_i})}\right)\bm{V}_q \myfrac{1}{J}\widehat{\bm{M}}^{-1}\bm{u},
\]
where $\bm{P}_q=\bm{\widehat{M}}^{-1} \bm{V}_q^T\text{diag}(\widehat{\bm{w}}_i)$ is the quadrature-based $L^2$ projection operator on the reference element.

For weight-adjusted inner products with matrix-valued weights, the corresponding weight-adjusted mass matrices approximate weighted $L^2$ mass matrices in a similar fashion
\begin{align*}
\bm{M_W}\approx(\bm{I} \otimes \bm{M}) \bm{M_{W^{-1}}}^{-1} (\bm{I} \otimes \bm{M}),\\
\bm{M_W}^{-1}\approx(\bm{I} \otimes \bm{M}^{-1}) \bm{M_{W^{-1}}} (\bm{I} \otimes \bm{M}^{-1}),
\end{align*}
where $\bm{M_W}$ is the matrix-weighted mass matrix defined in (\ref{eq:matweight}).  We note that $\bm{M_W}$ can be applied in a quadrature based fashion component by component.  

In the context of DG with explicit time-stepping, the factor of $\widehat{\bm{M}}^{-1}$ can be pre-multiplied into the right hand side (the spatial discretization). The application of the weight-adjusted mass matrices then requires only two reference matrices $\bm{V}_q$ and $\bm{P}_q$ and the values of the weight function at quadrature points $w(\bm{\widehat{x}}_i)$. In this study the number of quadrature points is $O(N^3)$ \cite{xiao2010}. The overall storage cost for applying weight-adjusted mass matrices using the above implementation is $O(N^3)$ per element, while the pre-computation and storage of DG operators involving inverses of weighted mass matrices require $O(N^6)$ storage per element.

We derive a weight-adjusted DG method by replacing the $L^2$ inner products in the left hand side of the DG formulation (\ref{eq17}) with weight-adjusted approximations. The right hand side of the WADG formulation is similar to the right hand side of the DG formulation (\ref{eq17}), and preserves a variant of the energy stability in Theorem 1 i.e.,
\begin{align*}
\sum_{D^k \in \Omega_h} \myfrac{1}{2} \myfrac{\partial }{\partial t} \left(\left(  T_{{\bm{Q}_s}^{-1}} \bm{\tau}, \bm{\tau}  \right)_{L^2(D^k)} + \left(  T_{{\bm{Q}_v}^{-1}} \bm{v}, \bm{v}  \right)_{L^2(D^k)}\right) \le 0,
\end{align*}
where $T_{{\bm{Q}_s}^{-1}}$ are $T_{{\bm{Q}_v}^{-1}}$ are weighting operators defined in \ref{appendix:A02}.

We replace the weighted $L^2$ mass matrices in (\ref{eq21}) and (\ref{eq22}) by their weight-adjusted approximations. Inverting these weight-adjusted mass matrices yields the following local system of ODEs for $\bm{\Sigma}$ and $\bm{V}$:
\begin{align*} 
\myfrac{\partial \bm{\Sigma}}{\partial t} &=\left(\bm{I} \otimes \bm{M}^{-1}\right)\bm{M}_{\bm{Q}_{s}^{-1}} \left(\sum_{i=1}^{d} \left(\bm{A}_i \otimes \bm{S_i}\right)\bm{V} + \sum_{f \in \partial D^k} \left(\bm{I} \otimes \bm{M}_f\right) \bm{F}_\tau\right), \\
\myfrac{\partial \bm{V}}{\partial t}&=\left(\bm{I} \otimes \bm{M}^{-1}\right)\bm{M}_{\bm{Q}_{v}^{-1}}\left(\sum_{i=1}^{d}\left({\bm{A}_i}^T \otimes \bm{S_i}\right)\bm{\Sigma} + \sum_{f \in \partial D^k} \left(\bm{I} \otimes \bm{M}_f\right) \bm{F}_v + \bm{M}_{\bm{D}}\bm{V}\right).
\end{align*}

Matrices $\left(\bm{I} \otimes \bm{M}^{-1}\right)\bm{M}_{\bm{Q}_{s}^{-1}}$ and $\left(\bm{I} \otimes \bm{M}^{-1}\right)\bm{M}_{\bm{Q}_{v}^{-1}}$ are applied in a matrix-free fashion using reference element matrices and values of $\bm{Q}_s^{-1}$ and $\bm{Q}_v^{-1}$ at quadrature points in (\ref{eq8}) and (\ref{eq9}). The convergence analysis of the numerical scheme is performed in \ref{appendix:A}.

\begin{table*}
	\caption{Material properties for several poroelastic media used in the examples \cite{carcione1996}}
	\centering
	\medskip
	\begin{tabular}{c c c c c} 
		\hline
		Properties  & Sandstone  & Epoxy-glass & Sandstone & Shale  \\ 
		                 & (Orthotropic)  & (Orthotropic)   & (Isotropic) & (Isotropic)  \\ [0.5ex] 
		\hline
		$K_s~$(GPa) & 80 &  40 & 40 & 7.6  \\ 
		$\rho_s~$ $($kg/m$^3$$)$ & 2500&1815  & 2500 & 2210\\
		$c_{11}~$(GPa) & 71.8 & 39.4 & 36 & 11.9\\ 
		$c_{12}~$(GPa) & 3.2 & 1.2 & 12 & 3.96\\ 
		$c_{13}~$(GPa) & 1.2 & 1.2 & 12 & 3.96\\ 
		$c_{33}~$(GPa) & 53.4 & 13.1 & 36 & 11.9\\ 
		$c_{55}~$(GPa) & 26.1& 3 & 12& 3.96\\ 
		$\phi~$ & 0.2 &0.2  & 0.2&0.16 \\
		$\kappa_1~(10^{-15}~{\rm m}^2)$ & 600 &600 &600&100 \\
		$\kappa_3~(10^{-15}~{\rm m}^2)$ & 100 & 100& 600&100\\
		$T_1~$ & 2 &2 & 2&2\\
		$T_3~$ & 3.6 & 3.6 & 2&2\\
		${K_f}~$(GPa) & 2.5 & 2.5 & 2.5&2.5 \\
		${\rho_f}~$(Kg/m$^3$) & 1040 &1040 & 1040 &1040\\
		${\eta}~$$(10^{-3}$ Kg/m.s) & 1 & 1 &1&1\\ [1ex]
		\hline
	\end{tabular}
\end{table*}
\section{Numerical experiments}
\label{sec:results}
\label{sec:num}
In this section, we present several numerical experiments validating the stability and accuracy of proposed method in two and three dimensions. The convergence of the new DG formulation in piecewise constant isotropic poroelastic media is confirmed. Finally, the method is applied to problems with anisotropy and micro-heterogeneities (sub-element variations).

In all experiments, we follow \cite{chan2017IV} and compute the application of weight-adjusted mass matrices using a quadrature which is exact for polynomials of degree $(2N+1)$.  Time integration is performed using the low-storage 4$^\text{th}$ order five-stage Runge-Kutta scheme of Carpenter and Kennedy \cite{carpenter1994}, and the time step is chosen based on the global estimate
\begin{align}
\label{eqdt}
dt=\min_k\myfrac{C_{CFL}}{\sup_{x \in \Omega} \left \Vert \bm{C(\bm{x})} \right \Vert_2 C_N \left \Vert J^f \right \Vert_{L^{\infty}(\partial D^k)} \left \Vert J^{-1} \right \Vert_{L^{\infty}(D^k)} }
\end{align}
where $C_N=O(N^2)$ is the order-dependent constant in the surface polynomial trace inequality \cite{chan2016I} and $C_{\rm CFL}$ is a tunable global CFL constant. This estimate is derived by bounding the eigenvalues of the spatial DG discretization matrix appearing in the semi-discrete system of ODEs.  This choice of $dt$ is very conservative as it is derived based on an upper bound on the spectral radius.  

\subsection{Spectra and choice of penalty parameter}
We first verify the energy stability of proposed DG formulation. Let $\bm{A}_h$ denote the matrix induced by the global semi-discrete DG formulation, such that time evolution of the solution $\bm{\tau}, \bm{v}$ is governed by 
\[
\myfrac{\partial \bm{Q}}{\partial t}=\bm{A}_h \bm{Q},
\]
where $\bm{Q}$ denotes a vector of degrees of freedom for $(\bm{\tau}, \bm{v})$. We show in Figure 1 eigenvalues of $\bm{A}_h$ for $\alpha_{\bm{\tau}}=\alpha_{\bm{v}}=0$ and $\alpha_{\bm{\tau}}=\alpha_{\bm{v}}=1$ with material parameters of isotropic sandstone (given in Table 1). The discretization parameters are $N=3$ and $h=1/2$. In both cases, the largest real part of any eigenvalues is $O(10^{-14})$, which suggests that the semi-discrete scheme is indeed energy stable.

For practical simulations, the choice of $\alpha_{\bm{\tau}},\alpha_{\bm{v}}$ remains to be specified. Taking $\alpha_{\bm{\tau}},\alpha_{\bm{v}}>0$ results in damping of of under-resolved spurious components of the solutions. However, a naive selection of $\alpha_{\bm{\tau}}, \alpha_{\bm{v}}$ can result in an overly restrictive time-step restriction for stability.  A guiding principle for determining appropriate values of the penalty parameters $\alpha_{\bm{\tau}}, \alpha_{\bm{v}}$ is to ensure that the spectral radius is the same magnitude as the case when $\alpha_{\bm{\tau}} = \alpha_{\bm{v}} = 0$.   For example, the spectral radius of $\bm{A}_h$, $\rho(\bm{A}_h)$ is $18.0034$ for $\alpha_{\bm{\tau}},\alpha_{\bm{v}}=0$ which is $O(N^2/h)$.  The spectral radius $\rho(\bm{A}_h)$ is $19.1217$ for $\alpha_{\bm{\tau}},\alpha_{\bm{v}}=0.5$, while the spectral radius for $\alpha_{\bm{\tau}},\alpha_{\bm{v}}=1$ is $\rho(\bm{A}_h) = 44.44$.  Since the maximum stable timestep is proportional to the spectral radius, taking $\alpha_{\bm{\tau}},\alpha_{\bm{v}}=1$ in this case results in a more restrictive CFL condition.  This phenomena is related to observations in \cite{chan2017II} that large penalty parameters result in extremal eigenvalues of $\bm{A}_h$ with very large negative real parts.  
The optimal choice of scaling also depends on the media heterogeneities on elements adjacent to an interface.  

In all following experiments, we use ${\alpha}_{\bm{\tau}}={\alpha}_{\bm{v}}=1$ unless specified otherwise.

\begin{figure}
\centering
\subfloat[$\alpha_{\bm{\tau}}=\alpha_{\bm{v}}=0$]{
\includegraphics[trim={4.0cm, 8.5cm, 4cm, 9.2cm}, clip, width=0.475\textwidth]{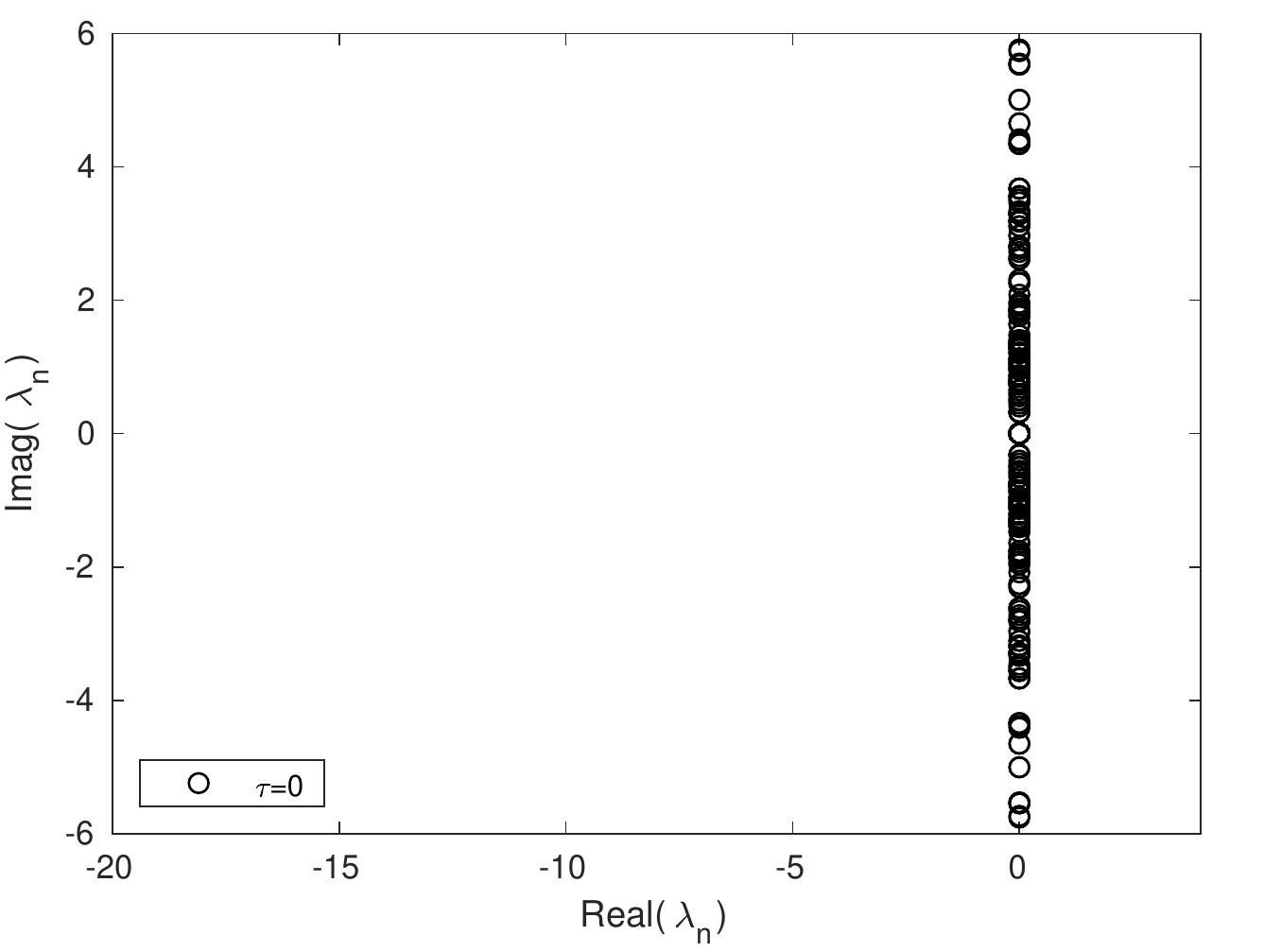}
\label{subfig:central flux}}
\hspace{.5em}
\subfloat[$\alpha_{\bm{\tau}}=\alpha_{\bm{v}}=1$]{
\includegraphics[trim={4.0cm, 8.5cm, 4cm, 9.2cm}, clip, width=0.475\textwidth]{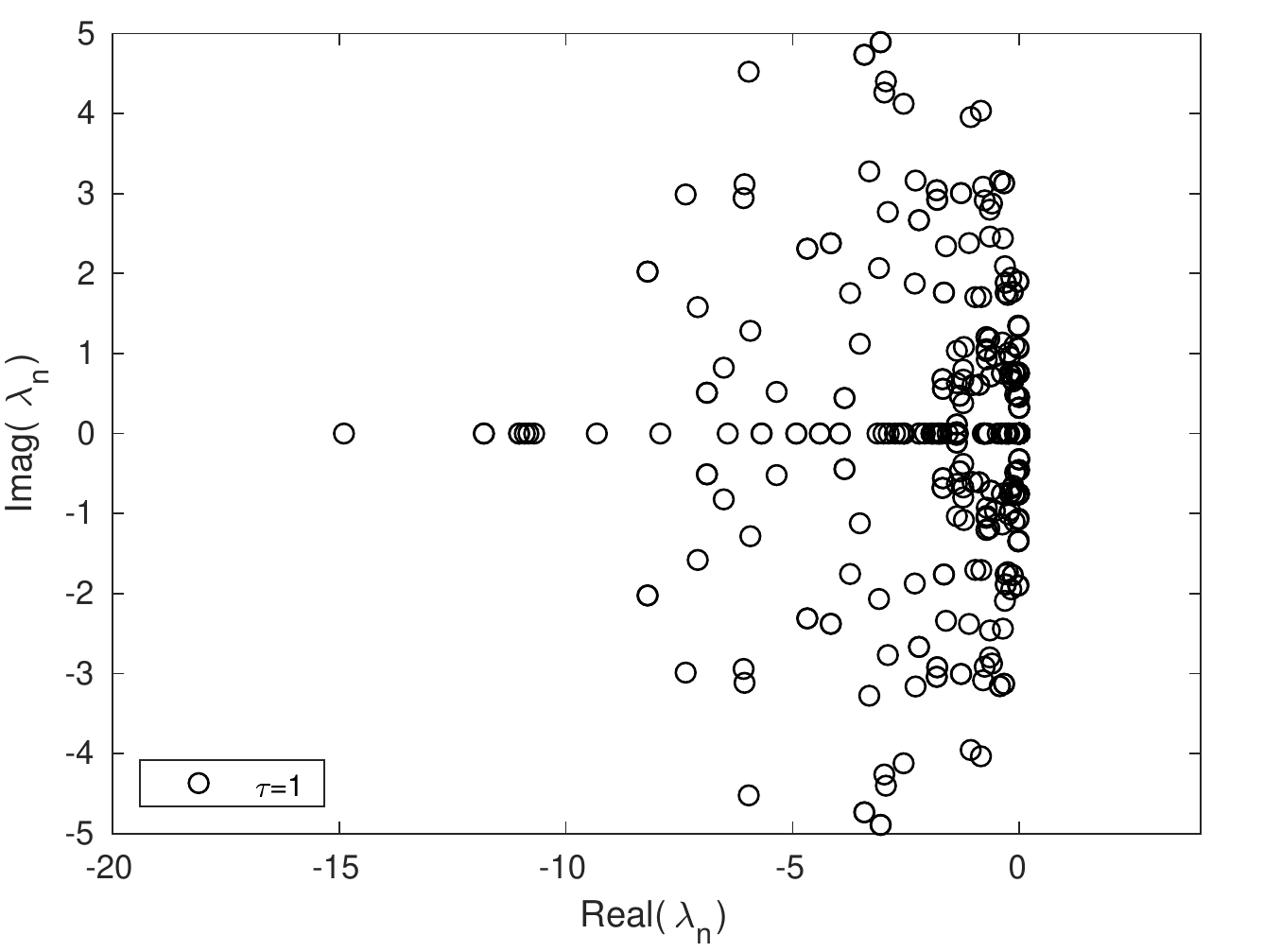}
\label{subfig:penalty flux}}
\caption{Spectra for $N=3$ and $h=1/2$ with a material property of isotropic Sandstone (Table 1). For $\alpha_{\bm{\tau}}$ = $\alpha_{\bm{v}}=0$ and $\alpha_{\bm{\tau}}$ = $\alpha_{\bm{v}}=1$, the largest real part of spectra are 1.6431e-14 and 6.3412e-15, respectively.}
\label{fig:DG spectra}
\end{figure}

\subsection{Analytic solution}
Next, we study the accuracy and convergence of the our DG method for a plane wave propagating in an isotropic porous sandstone with material properties given in Table 1 (Column 4). Unless otherwise stated, we report relative $L^2$ errors for all components of the solution $\bm{U}=(\bm{\tau}, \bm{v})$
\[
\myfrac{\left \Vert \bm{U} -\ \bm{U}_h \right \Vert_{L^2(\Omega)}}{\left \Vert \bm{U} \right \Vert_{L^2(\Omega)}}=\myfrac{\left( \sum_{i=1}^m \left \Vert \bm{U_i} -\bm{U}_{\bm{i},\bm{h}} \right \Vert ^2_{L^2(\Omega)}\right)^{1/2}}{\left(\sum_{i=1}^m \left \Vert \bm{U}_i \right \Vert^2_{L^2(\Omega)} \right)^{1/2}}.
\]
\subsection{Plane wave in a poroelastic medium}
The analytical solution to (\ref{eq 7}) for a plane wave is given as
\begin{align}
\label{eq31}
\bm{Q}_n(\bm{x},t)=\bm{Q}_n^0 \exp[\text{i} \cdot (\omega t- \bm{k}\cdot \bm{x})], \qquad n=1...13,
\end{align}
where $\bm{Q}_n^0$ is the initial amplitude vector of stress and velocity components; $\omega$ are wave frequencies; $\bm{k}=(k_x, k_y, k_z)$ is the wave-number vector.
To achieve realistic poroelastic behavior, we superimpose three plane waves, of the form given by (\ref{eq31}), corresponding to a fast P-wave, an S-wave and a slow P-wave.

Now, we briefly describe how we determine the the wave frequencies $\omega$. Substituting (\ref{eq31}) into (\ref{eq 7}) yields
\begin{align}
\label{eq32}
\omega \bm{Q}_n^0=(\bm{A} k_x + \bm{B} k_y + \bm{C} k_z -\text{i} \bm{E})\bm{Q}_n^0
\end{align}
Solving the three eigenvalues problem for in (\ref{eq31}) for each wave mode $l$ yields in matrix of right eigenvectors $(R^{(l)}_{mn})$ and eigenvalue $(\omega_{l})$.  Following \cite{toro1999,de2008}, the solution of (\ref{eq 7}) can be constructed as
\begin{align}
\label{eq33}
\bm{Q}_n(\bm{x},t)=\sum_{l=1}^3 R^{(l)}_{mn} \gamma_{n}^{(l)} \exp[\text{i} \cdot (\omega^{(l)} t- \bm{k^{(l)}}\cdot \bm{x})],
\end{align}
where $\gamma_{n}^{(l)}$ is a amplitude coefficient with $\gamma_{1}^{(1)}=\gamma_{2}^{(2)}=\gamma_{4}^{(3)}=100$.

For completeness, we perform convergence analyses for inviscid and viscid media separately.  

\subsubsection{Inviscid case ($\eta=0$)}

The error is computed for an inviscid brine-filled ($\eta=0)$ isotropic sandstone.  The inviscid nature of the fluid implies that $\bm{D}=\bm{0}$. We show in Figure 2 the $L^2$ errors computed at T=1 and CFL=1, using uniform triangular meshes constructed by bisecting an uniform mesh of quadrilaterals along the diagonal. Figure 2a  and 2b show the convergence plot using the central flux $(\alpha_{\bm{v}}=\alpha_{\bm{\tau}}=0)$ and penalty flux $(\alpha_{\bm{v}}=\alpha_{\bm{\tau}}=1)$ respectively. For $N=1,...,5$, $O(h^{N+1})$, rates of convergence are observed. We note that for $N=4$ and $N=5$, we observe results for both fluxes which are better than the $4^{\text{th}}$ order accuracy of our time-stepping scheme. This is most likely due to the benign nature of the solution in time and the choice of time step (\ref{eqdt}) which scales as $O(h/N^2)$. For $N=4,5$, the results of Figure 2 suggest that the resulting time step is small enough such that temporal errors of $O(dt^4)$ are small relative to spatial discretization errors of $O(h^{N+1})$.  
\begin{figure}
\centering
\subfloat[$\tau_v=\tau_\sigma=0$ (central flux)]{
\begin{tikzpicture}
\begin{loglogaxis}[
	legend cell align=left,
	width=.475\textwidth,
    xlabel={Mesh size $h$},
       ylabel={$L^2$ errors in $\bm{Q}$}, 
    xmin=.0125, xmax=0.75,
    ymin=1e-12, ymax=1.5,
    legend pos=south east, legend cell align=left, legend style={font=\tiny},
    xmajorgrids=true,
    ymajorgrids=true,
    grid style=dashed,
] 

\addplot[color=blue,mark=*,semithick, mark options={fill=markercolor}]
coordinates{(0.5,0.8188)(0.25,0.3454)(0.125,0.0763)(0.0625,0.0196)};
\logLogSlopeTriangleFlip{0.53}{0.125}{0.858}{1.96}{blue}

\addplot[color=red,mark=square*,semithick, mark options={fill=markercolor}]
coordinates{(0.5,0.2156)(0.25,0.0375)(0.125,0.0048)(0.0625,6.8857e-04)};
\logLogSlopeTriangleFlip{0.53}{0.125}{0.742}{2.80}{red}

\addplot[color=black,mark=triangle*,semithick, mark options={fill=markercolor}]
coordinates{(0.5, 4.6561e-02)(0.25,2.7666e-03)(0.125,2.036e-04)(0.0625,1.0882e-05)};
\logLogSlopeTriangleFlip{0.53}{0.125}{0.597}{4.23}{black}

\addplot[color=magenta,mark=diamond*,semithick, mark options={fill=markercolor}]
coordinates{(0.5,6.4012e-03)(0.25,2.2910e-04)(0.125,7.5020e-06)(0.0625,1.9958e-07)};
\logLogSlopeTriangleFlip{0.53}{0.125}{0.460}{5.23}{magenta}

\addplot[color=violet,mark=pentagon*,semithick, mark options={fill=markercolor}]
coordinates{(0.5,8.4501e-04)(0.25,1.3536e-05)(0.125,2.0354e-07)(0.0625,3.4217e-09)};
\logLogSlopeTriangleFlip{0.53}{0.125}{0.310}{5.94}{violet}

\legend{$N=1$,$N=2$,$N=3$,$N=4$,$N=5$}
\end{loglogaxis}
\end{tikzpicture}
}
\subfloat[$\tau_v=\tau_\sigma=1$ (penalty flux)]{
\begin{tikzpicture}
\begin{loglogaxis}[
	legend cell align=left,
	width=.475\textwidth,
    xlabel={Mesh size $h$},
       ylabel={$L^2$ errors in $\bm{Q}$}, 
    xmin=.0125, xmax=0.75,
    ymin=1e-12, ymax=1.5,
    legend pos=south east, legend cell align=left, legend style={font=\tiny},
    xmajorgrids=true,
    ymajorgrids=true,
    grid style=dashed,
] 

\addplot[color=blue,mark=*,semithick, mark options={fill=markercolor}]
coordinates{(0.5,5.9239e-01)(0.25,1.3794e-01)(0.125,2.9369e-02)(0.0625,6.8824e-03)};
\logLogSlopeTriangleFlip{0.53}{0.125}{0.825}{2.09}{blue}

\addplot[color=red,mark=square*,semithick, mark options={fill=markercolor}]
coordinates{(0.5,1.5494e-01)(0.25,2.3803e-02)(0.125,3.3863e-03)(0.0625,4.0058e-04)};
\logLogSlopeTriangleFlip{0.53}{0.125}{0.726}{3.08}{red}

\addplot[color=black,mark=triangle*,semithick, mark options={fill=markercolor}]
coordinates{(0.5, 2.8130e-02)(0.25,1.8190e-03)(0.125,1.1702e-04)(0.0625,6.8214e-06)};
\logLogSlopeTriangleFlip{0.53}{0.125}{0.580}{4.10}{black}

\addplot[color=magenta,mark=diamond*,semithick, mark options={fill=markercolor}]
coordinates{(0.5,4.2389e-03)(0.25,1.3346e-04)(0.125,4.1682e-06)(0.0625,1.3188e-07)};
\logLogSlopeTriangleFlip{0.53}{0.125}{0.445}{4.98}{magenta}

\addplot[color=violet,mark=pentagon*,semithick, mark options={fill=markercolor}]
coordinates{(0.5,5.3616e-04)(0.25,8.5470e-06)(0.125,1.2937e-07)(0.0625,2.0485e-09)};
\logLogSlopeTriangleFlip{0.53}{0.125}{0.290}{5.98}{violet}

\legend{$N=1$,$N=2$,$N=3$,$N=4$,$N=5$}
\end{loglogaxis}
\end{tikzpicture}
}
\caption{Convergence of $L^2$ error for plane wave in porolastic media with $\eta=0$-inviscid case }
\end{figure}
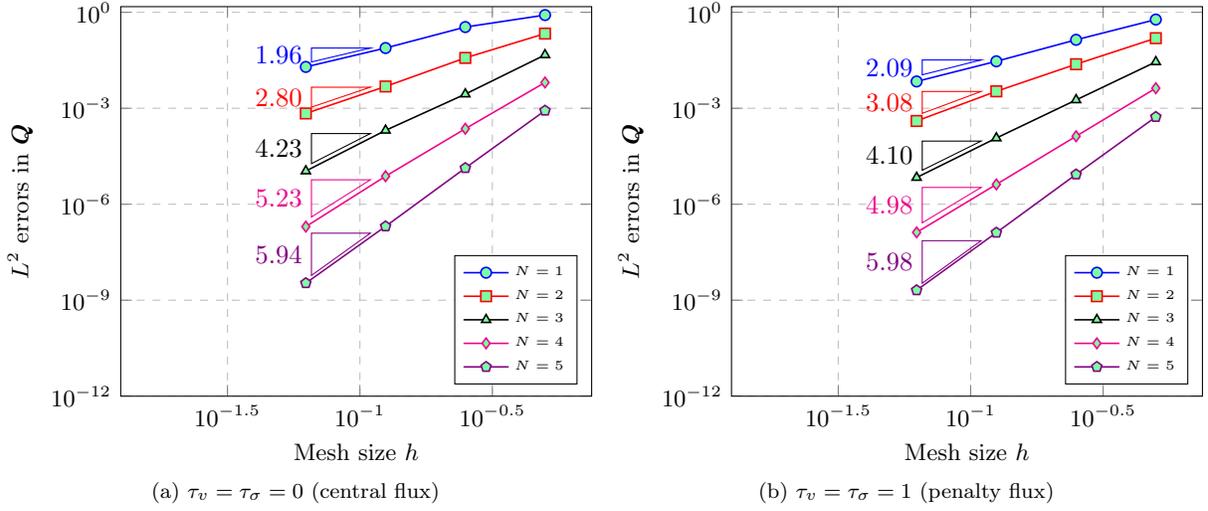

\subsubsection{Viscid case ($\eta \ne 0$)}
The poroelastic formulations used by Carcione \cite{carcione1996} and de la Puente et al.\ \cite{de2008} produce a stiff system of ODEs due to the presence of the dissipative matrix $\bm{D}$.  Consequently, Carcione used a Strang's 2$^\text{nd}$ order operator splitting approach \cite{strang1968} to avoid small $\Delta t$, while de la Puente et al. \cite{de2008} circumvented the effect of stiffness using a local implicit time-stepping approach to achieve the convergence rate of $O(h^{N+1})$.  The poroelastic formulation and numerical discretization used in this study do not appear to produce equally stiff systems of ODEs, and Strang splitting is not required.  However, for a complete analysis, the convergence of the proposed DG formulation in viscid case is studied using both a unified DG scheme (where we incorporate the dissipative matrix $\bm{D}$ directly into the explicit time-stepping scheme) as well as a second order Strang operator splitting \cite{strang1968}.  

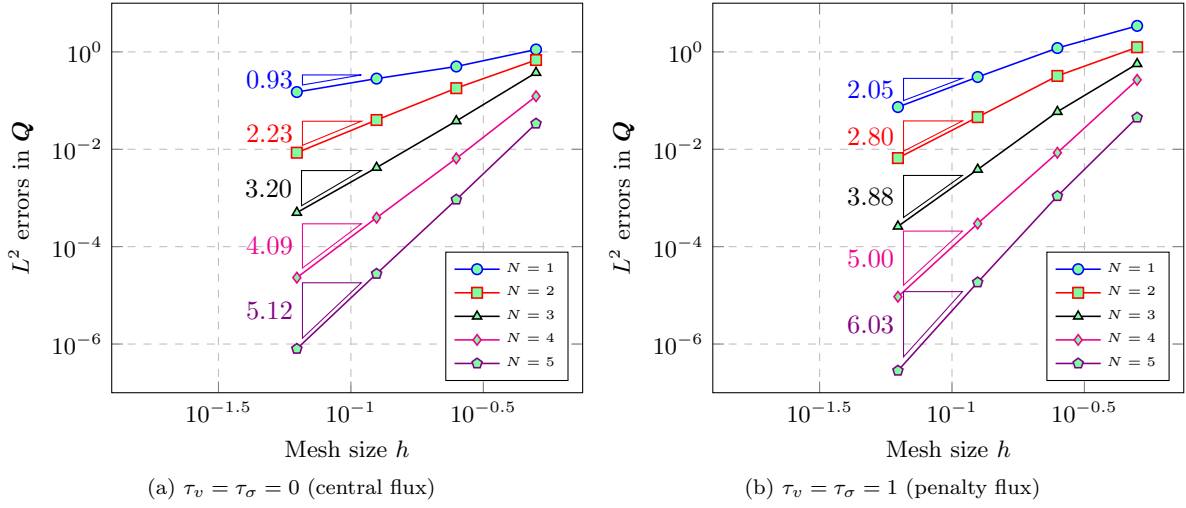
\begin{figure}
\centering
\subfloat[$\tau_v=\tau_\sigma=0$ (central flux)]{
\begin{tikzpicture}
\begin{loglogaxis}[
	legend cell align=left,
	width=.475\textwidth,
    xlabel={Mesh size $h$},
       ylabel={$L^2$ errors in $\bm{Q}$}, 
    xmin=.0125, xmax=0.75,
    ymin=1e-7, ymax=10,
    legend pos=south east, legend cell align=left, legend style={font=\tiny},
    xmajorgrids=true,
    ymajorgrids=true,
    grid style=dashed,
] 

\addplot[color=blue,mark=*,semithick, mark options={fill=markercolor}]
coordinates{(0.5,1.1197)(0.25,0.5012)(0.125,0.2827)(0.0625,0.1498)};
\logLogSlopeTriangleFlip{0.53}{0.125}{0.790}{0.93}{blue}

\addplot[color=red,mark=square*,semithick, mark options={fill=markercolor}]
coordinates{(0.5,0.6809)(0.25,0.1798)(0.125,0.0399)(0.0625, 0.0085)};
\logLogSlopeTriangleFlip{0.53}{0.125}{0.635}{2.23}{red}

\addplot[color=black,mark=triangle*,semithick, mark options={fill=markercolor}]
coordinates{(0.5, 0.3741)(0.25,0.0381)(0.125,0.0042)(0.0625, 4.9974e-04)};
\logLogSlopeTriangleFlip{0.53}{0.127}{0.480}{3.20}{black}

\addplot[color=magenta,mark=diamond*,semithick, mark options={fill=markercolor}]
coordinates{(0.5,0.1228)(0.25, 0.0065)(0.125,3.9282e-04)(0.0625,2.3099e-05)};
\logLogSlopeTriangleFlip{0.53}{0.125}{0.320}{4.09}{magenta}

\addplot[color=violet,mark=pentagon*,semithick, mark options={fill=markercolor}]
coordinates{(0.5,0.0337)(0.25,9.2807e-04)(0.125,2.7704e-05)(0.0625,7.9856e-07)};
\logLogSlopeTriangleFlip{0.53}{0.125}{0.140}{ 5.12}{violet}

\legend{$N=1$,$N=2$,$N=3$,$N=4$,$N=5$}
\end{loglogaxis}
\end{tikzpicture}
}
\subfloat[$\tau_v=\tau_\sigma=1$ (penalty flux)]{
\begin{tikzpicture}
\begin{loglogaxis}[
	legend cell align=left,
	width=.475\textwidth,
    xlabel={Mesh size $h$},
       ylabel={$L^2$ errors in $\bm{Q}$}, 
    xmin=.0125, xmax=0.75,
    ymin=1e-7, ymax=10,
    legend pos=south east, legend cell align=left, legend style={font=\tiny},
    xmajorgrids=true,
    ymajorgrids=true,
    grid style=dashed,
] 

\addplot[color=blue,mark=*,semithick, mark options={fill=markercolor}]
coordinates{(0.5,3.4213e+00)(0.25,1.1994e+00)(0.125,3.0470e-01)(0.0625,7.3608e-02)};
\logLogSlopeTriangleFlip{0.53}{0.125}{0.750}{2.05}{blue}

\addplot[color=red,mark=square*,semithick, mark options={fill=markercolor}]
coordinates{(0.5,1.2444e+00)(0.25,3.2086e-01)(0.125,4.5737e-02)(0.0625, 6.5921e-03)};
\logLogSlopeTriangleFlip{0.53}{0.125}{0.620}{2.80}{red}

\addplot[color=black,mark=triangle*,semithick, mark options={fill=markercolor}]
coordinates{(0.5, 5.6768e-01)(0.25,5.9024e-02)(0.125,3.8388e-03)(0.0625, 2.6091e-04)};
\logLogSlopeTriangleFlip{0.53}{0.125}{0.450}{3.88}{black}

\addplot[color=magenta,mark=diamond*,semithick, mark options={fill=markercolor}]
coordinates{(0.5,2.6650e-01)(0.25, 8.4785e-03)(0.125,2.9762e-04)(0.0625,9.3503e-06)};
\logLogSlopeTriangleFlip{0.53}{0.125}{0.276}{5.00}{magenta}

\addplot[color=violet,mark=pentagon*,semithick, mark options={fill=markercolor}]
coordinates{(0.5,4.4851e-02)(0.25,1.0933e-03)(0.125,1.8481e-05)(0.0625,2.8396e-07)};
\logLogSlopeTriangleFlip{0.53}{0.125}{0.092}{ 6.03}{violet}

\legend{$N=1$,$N=2$,$N=3$,$N=4$,$N=5$}
\end{loglogaxis}
\end{tikzpicture}
}
\caption{Convergence of $L^2$ error for plane wave in poroelastic media  for viscid case 
($\eta \ne 0$) with  unified DG scheme}
\end{figure}

Operator splitting separates the dissipative term from the conservative term at each time step. We solve the dissipative part $( \text{stiff})$ part of the system analytically and solve the conservative part by using the proposed DG formulation.  We rewrite the system with non-zero forcing function as follows
\begin{align}
\label{eq35}
 \myfrac{\partial \bm{Q}}{\partial t}={\bm{A}_h}\bm{Q} + \bm{f},
\end{align}
The formal solution of (\ref{eq35}) is given as

\begin{align}
\bm{Q}(t)=\exp(\bm{A}_ht)\bm{Q_0} + \int_0^t \exp(\tau \bm{A}_h) \bm{f}(t-\tau) d\tau,
\end{align}
where $\exp(\bm{A}_ht)$ is an evolution operator.

Using the operator splitting approach, the propagation matrix can be partitioned as
\begin{align} 
\bm{A}_h=\bm{A}_{c} + \bm{A}_{d},
\end{align}
where the subscript $c$ indicates the matrix representing the conservative part of the system, and the subscript $d$ indicates the diffusive matrix, representing low order terms coupled with the viscosity $(\eta)$. Subsequently, the evolution operator can be expressed as 
\begin{align}
\label{eq38}
\exp(\bm{A}_h t)=\exp[(\bm{A}_{h_c} + \bm{A}_{h_d})t].
\end{align}
By using the product formula, (\ref{eq38}) is expressed as
\begin{align}
\label{eq39}
\exp(\bm{A}_h dt)=\exp\left(\myfrac{1}{2} \bm{A}_{h_d} dt\right)\exp( \bm{A}_{h_c} dt)\exp\left(\myfrac{1}{2} \bm{A}_{h_d} dt\right).
\end{align} 
It can be shown that (\ref{eq38}) is second order accurate in $dt$.  Thus, (\ref{eq39}) allows us to solve the stiff part separately.  The solution of the stiff part of the system can be derived analytically and is shown in Appendix A.

\begin{figure}[H]
\centering
\subfloat[$\tau_v=\tau_\sigma=0$ (central flux)]{
\begin{tikzpicture}
\begin{loglogaxis}[
	legend cell align=left,
	width=.475\textwidth,
    xlabel={Mesh size $h$},
       ylabel={$L^2$ errors in $\bm{Q}$}, 
    xmin=.0125, xmax=0.75,
    ymin=1e-7, ymax=10,
    legend pos=south east, legend cell align=left, legend style={font=\tiny},
    xmajorgrids=true,
    ymajorgrids=true,
    grid style=dashed,
] 

\addplot[color=blue,mark=*,semithick, mark options={fill=markercolor}]
coordinates{(0.5,3.3359e+00)(0.25,1.5147e+00)(0.125,7.4787e-01)(0.0625,3.3447e-01)};
\logLogSlopeTriangleFlip{0.53}{0.125}{0.830}{1.16}{blue}

\addplot[color=red,mark=square*,semithick, mark options={fill=markercolor}]
coordinates{(0.5,2.0314e+00)(0.25,5.3700e-01)(0.125,1.1957e-01)(0.0625, 2.5297e-02)};
\logLogSlopeTriangleFlip{0.53}{0.125}{0.695}{2.24}{red}

\addplot[color=black,mark=triangle*,semithick, mark options={fill=markercolor}]
coordinates{(0.5, 1.1175e+00)(0.25,1.1385e-01)(0.125,1.2548e-02)(0.0625, 1.4956e-03)};
\logLogSlopeTriangleFlip{0.53}{0.125}{0.538}{3.07}{black}

\addplot[color=magenta,mark=diamond*,semithick, mark options={fill=markercolor}]
coordinates{(0.5,3.6657e-01)(0.25, 1.9445e-02)(0.125,1.1760e-03)(0.0625,6.9114e-05)};
\logLogSlopeTriangleFlip{0.53}{0.125}{0.376}{4.10}{magenta}

\addplot[color=violet,mark=pentagon*,semithick, mark options={fill=markercolor}]
coordinates{(0.5,1.0083e-01)(0.25,2.7767e-03)(0.125,8.2784e-05)(0.0625,2.4046e-06)};
\logLogSlopeTriangleFlip{0.53}{0.125}{0.200}{ 5.11}{violet}

\legend{$N=1$,$N=2$,$N=3$,$N=4$,$N=5$}
\end{loglogaxis}
\end{tikzpicture}
}
\subfloat[$\tau_v=\tau_\sigma=1$ (penalty flux)]{
\begin{tikzpicture}
\begin{loglogaxis}[
	legend cell align=left,
	width=.475\textwidth,
    xlabel={Mesh size $h$},
       ylabel={$L^2$ errors in $\bm{Q}$}, 
    xmin=.0125, xmax=0.75,
    ymin=1e-7, ymax=10,
    legend pos=south east, legend cell align=left, legend style={font=\tiny},
    xmajorgrids=true,
    ymajorgrids=true,
    grid style=dashed,
] 

\addplot[color=blue,mark=*,semithick, mark options={fill=markercolor}]
coordinates{(0.5,3.4215e+00)(0.25,1.1994e+00)(0.125,3.0469e-01)(0.0625,7.3607e-02)};
\logLogSlopeTriangleFlip{0.53}{0.125}{0.750}{2.05}{blue}

\addplot[color=red,mark=square*,semithick, mark options={fill=markercolor}]
coordinates{(0.5,1.2444e+00)(0.25,3.2085e-01)(0.125,4.5737e-02)(0.0625, 6.5920e-03)};
\logLogSlopeTriangleFlip{0.53}{0.125}{0.620}{2.80}{red}

\addplot[color=black,mark=triangle*,semithick, mark options={fill=markercolor}]
coordinates{(0.5, 5.6767e-01)(0.25,5.9024e-02)(0.125,3.8388e-03 )(0.0625, 2.6091e-04)};
\logLogSlopeTriangleFlip{0.53}{0.125}{0.450}{3.88}{black}

\addplot[color=magenta,mark=diamond*,semithick, mark options={fill=markercolor}]
coordinates{(0.5,2.6650e-01)(0.25, 8.4784e-03)(0.125,2.9763e-04)(0.0625,9.3670e-06)};
\logLogSlopeTriangleFlip{0.53}{0.125}{0.276}{5.00}{magenta}

\addplot[color=violet,mark=pentagon*,semithick, mark options={fill=markercolor}]
coordinates{(0.5,4.4851e-02)(0.25,1.0933e-03)(0.125,1.8512e-05)(0.0625,3.9166e-07)};
\logLogSlopeTriangleFlip{0.53}{0.125}{0.092}{ 5.60}{violet}

\legend{$N=1$,$N=2$,$N=3$,$N=4$,$N=5$}
\end{loglogaxis}
\end{tikzpicture}
}
\caption{Convergence of $L^2$ error for plane wave in poroelastic media  for viscid case 
($\eta \ne 0$) with  paired DG and Strang splitting approach }
\end{figure}
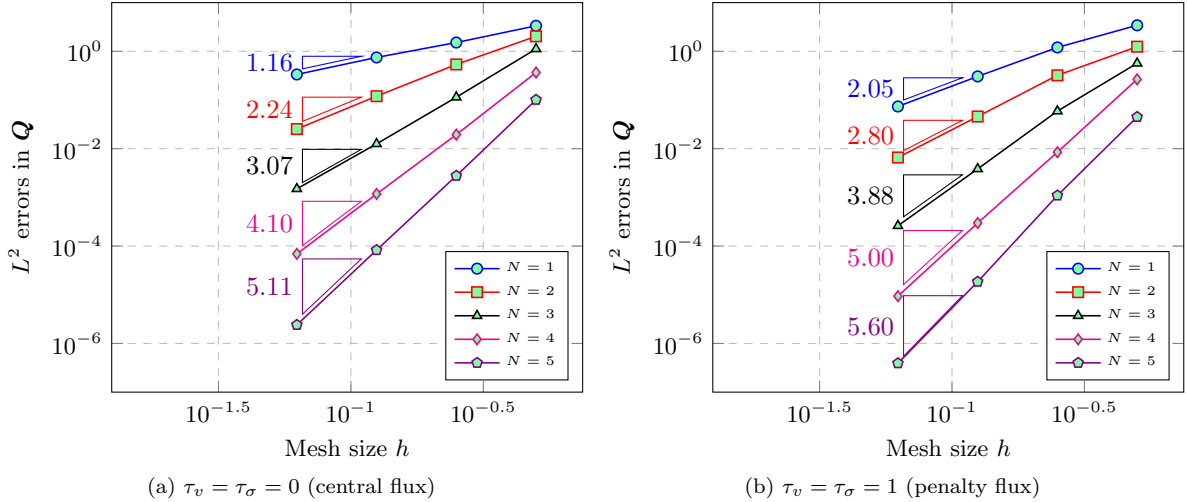
We show in Figure 3 the $L^2$ errors for viscid case, computed at T=1 and CFL=1 using uniform triangular meshes. Figure 3a  and 3b show convergence plots using the central flux $(\alpha_{\bm{v}}=\alpha_{\bm{\tau}}=0)$ and penalty flux $(\alpha_{\bm{v}}=\alpha_{\bm{\tau}}=1)$, respectively.  For $N=1,...,5$, $O(h^{N+1})$ rate of convergence are observed when using the penalty flux $\alpha_{\bm{v}}=\alpha_{\bm{\tau}}=1$.  When using a central flux $(\alpha_{\bm{v}}=\alpha_{\bm{\tau}}=0)$, we observe a so-called ``even-odd" pattern  \cite{hesthaven2007, mercerat2015}, with convergence rate $O(h^N)$ for odd $N$ and between $O(h^{N+1/2})$ and $O(h^{N+1})$ for even $N$.  

Similar to the inviscid case, for $N=4$ and $N=5$, we observe results for both fluxes which are better than the $4^{\text{th}}$ order accuracy of our time-stepping, due to the fact that temporal error are small relative to spatial discretization error. Furthermore, the spectral radius  $\rho(\bm{A}_h)$ for $N=4$ and $h=1/2$ is 25.91 which is of $O(N^2/h)$.   Since ${N^2}/{h} > \vert \vert \bm{D} = 1.667 \vert \vert$, the dissipative term does not increase the stiffness of the high order DG scheme.

We show in Figure 4 the $L^2$ errors for the viscid case using Strang splitting, computed at T=1 and CFL=1. Figure 4a and 4b show convergence plots using the central flux $(\alpha_{\bm{v}}=\alpha_{\bm{\tau}}=0)$ and penalty flux $(\alpha_{\bm{v}}=\alpha_{\bm{\tau}}=1)$, respectively.  Observations on convergence are similar to those in the unified DG scheme. Here, it is worth noting that for $N=5$ the convergence rate is $5.60$, which does not match the optimal theoretical rate of $O(h^{N+1})$.  This is may be due to the second order accuracy of the splitting scheme resulting in larger temporal errors relative to the spatial discretization errors $O(h^{N+1})$.

\subsection{Application examples}

We next demonstrate the accuracy and flexibility of the proposed DG for several application-based problems in linear poroelasticity with micro heterogeneities and anisotropy.  All computations are done using penalty parameters $\alpha_{\bm{v}}=\alpha_{\bm{\tau}}=1$ unless specified otherwise. In the subsequent sections, the field $\bm{b}$ represents the center of mass particle velocity vector \cite{sahay1994}, which is expressed as
\begin{align}
\label{com}
\bm{b}=\bm{v} + \left(\myfrac{\rho_f}{\rho}\right)\bm{q} .
\end{align}
In subsequent simulations, the forcing is applied to both the  $z-$ component of stress $\tau_{zz}$ and the fluid pressure $p$ by a Ricker wavelet point source 
\begin{align}
\label{Ricker}
f(\bm{x}, t)=(1-2(\pi f_0(t-t_0))^2)\exp[-(\pi f_0(t-t_0))^2] \delta(\bm{x}-\bm{x}_0),
\end{align}
where $\bm{x}_0$ is the position of the point source and $f_0$ is the central frequency.

In the following simulations, three types of poroelastic waves are observed: the fast P wave, S wave, and slow P wave. The fast P wave has solid and fluid motion in phase, while the slow P wave (Biot's mode) has the solid and fluid motion out of phase with one another. At low frequencies, the slow P wave is diffusive in character as viscous forces dominate over inertial forces. However, at high frequencies, the dominance of inertial forces over viscous forces results in the propagation of the slow P wave.

\subsubsection{Orthotropic sandstone}
To illustrate the effect of anisotropy on poroelastic wave propagation, we perform a computational experiment in orthotropic sandstone with material properties given in Table 1.  The size of the computational domain is $18.25~\text{m} \times 18.25~\text{m} $. The domain is discretized with uniform triangular element with a minimum edge length of $5~\text{cm}$. Figures 5(a)-(d) represent the $x-$ and $z-$ components of the center of mass particle velocity of the orthotropic sandstone, where (a) and (b) correspond to the inviscid case $(\eta=0)$, and (c) and (d) to the viscid case $(\eta \ne 0)$. The central frequency of the forcing function is $f_0=3730~\text{Hz}$, and polynomials of degree $N=4$ are used. The propagation time is $1.56~\text{ms}$. Three events can be observed: the fast P mode ($\text{P}_ \text{f}$, outer wavefront), the shear wave (S, middle wavefront), and the slow P mode ($\text{P}_ \text{s} $, inner wavefront). In the viscid case, the slow mode diffuses faster and the medium behaves almost as a single phase medium.

\begin{figure}
\centering
\subfloat[Orthotropic Sandstone, $b_x$ with $\eta = 0$]{
\begin{tikzpicture}
\begin{axis}[enlargelimits=false, axis on top, axis equal image, xlabel=x (m), ylabel=y (m), width=0.6\textwidth]

\addplot graphics [xmin=0,xmax=20,ymin=0,ymax=20] {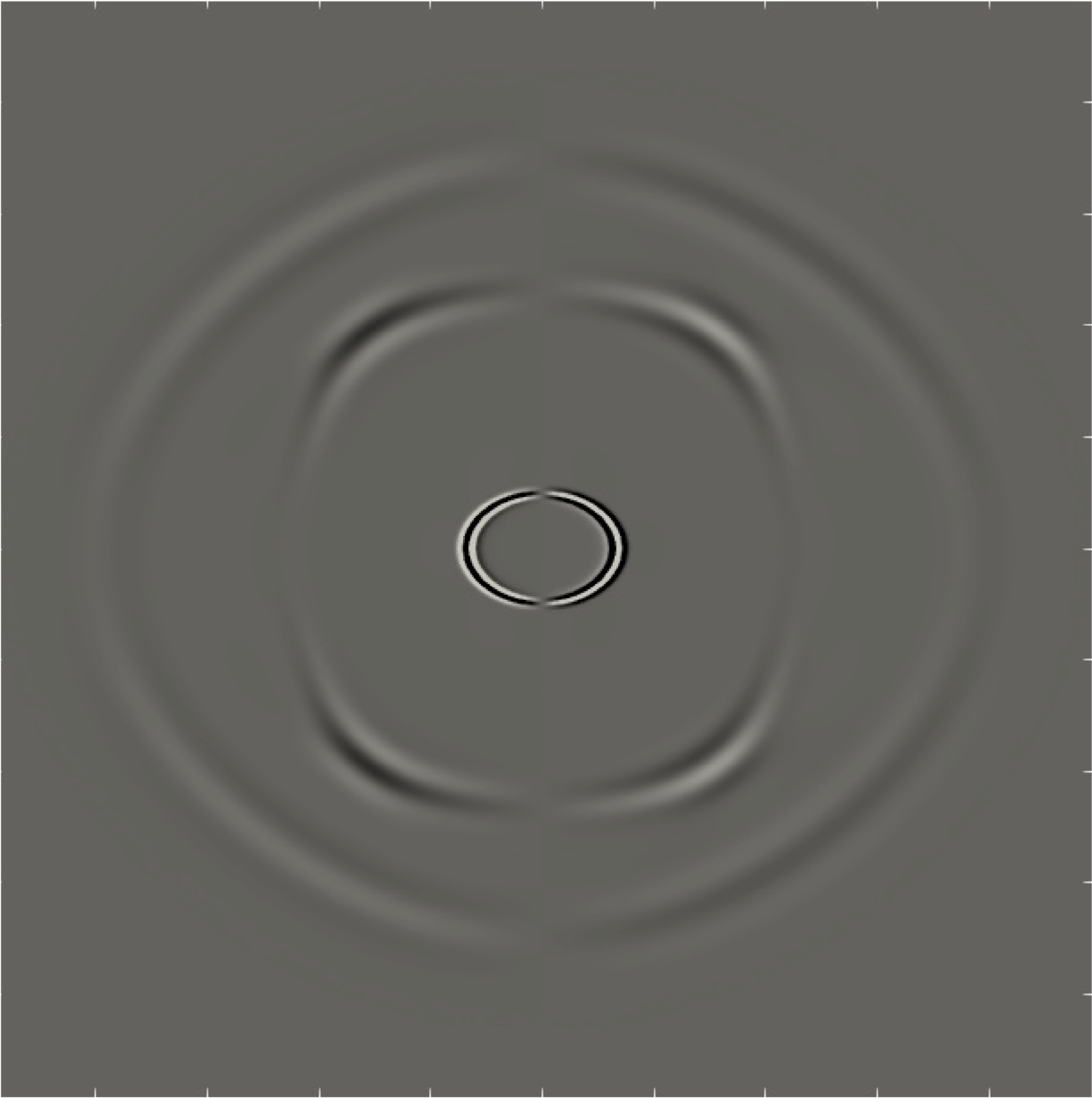};
\end{axis}
\end{tikzpicture}
}
\subfloat[Orthotropic Sandstone, $b_z$ with $\eta = 0$]{
\begin{tikzpicture}
\begin{axis}[enlargelimits=false, axis on top, axis equal image, xlabel=x (m), ylabel=y (m), width=0.6\textwidth]

\addplot graphics [xmin=0,xmax=20,ymin=0,ymax=20] {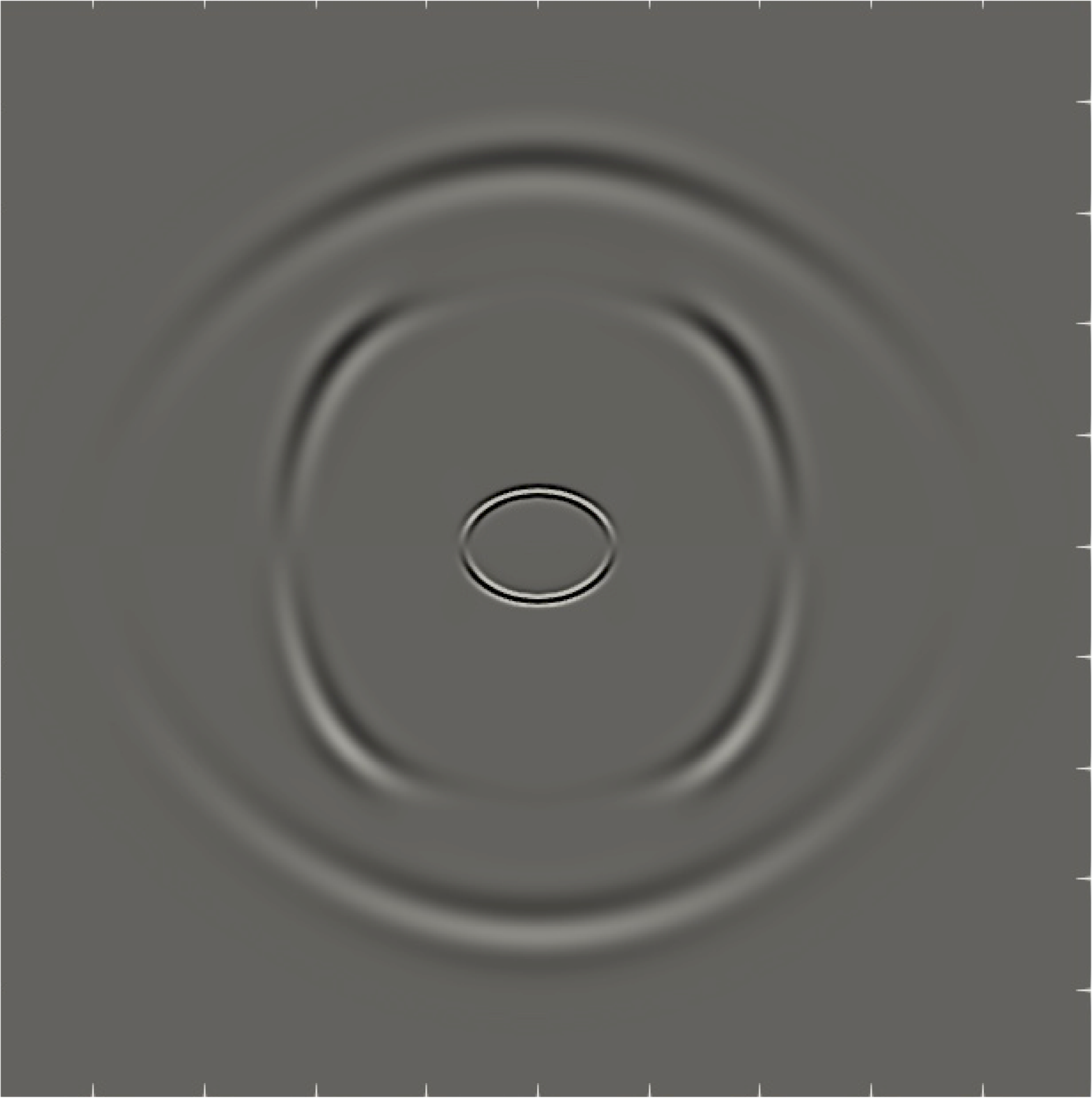};
\end{axis}
\end{tikzpicture}
}\\
\subfloat[Orthotropic Sandstone, $b_x$ with $\eta \ne 0$]{
\begin{tikzpicture}
\begin{axis}[enlargelimits=false, axis on top, axis equal image, xlabel=x (m), ylabel=y (m), width=0.6\textwidth]

\addplot graphics [xmin=0,xmax=20,ymin=0,ymax=20] {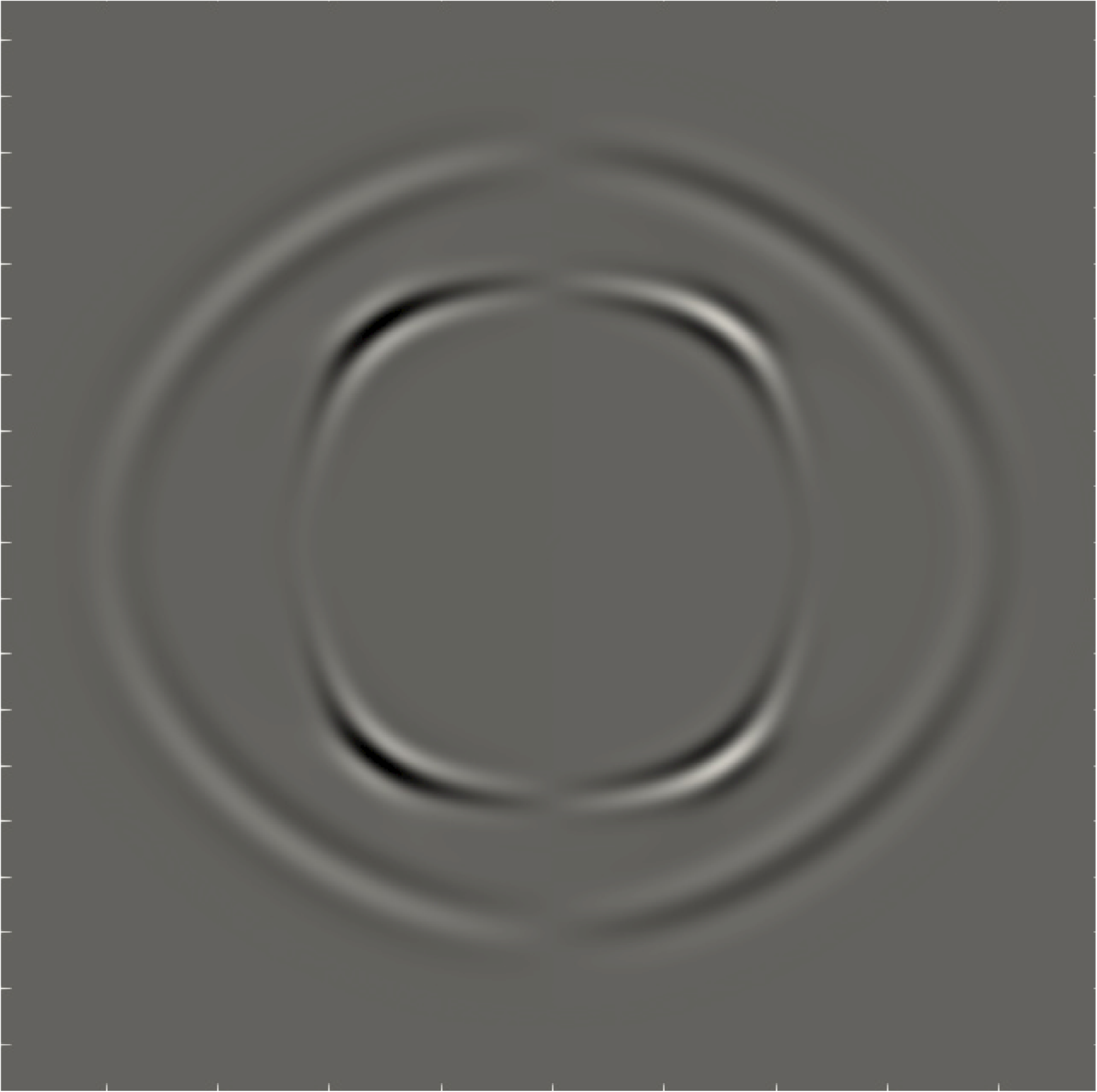};
\end{axis}
\end{tikzpicture}
}
\subfloat[Orthotropic Sandstone, $b_z$ with $\eta \ne 0$]{
\begin{tikzpicture}
\begin{axis}[enlargelimits=false, axis on top, axis equal image, xlabel=x (m), ylabel=y (m), width=0.6\textwidth]

\addplot graphics [xmin=0,xmax=20,ymin=0,ymax=20] {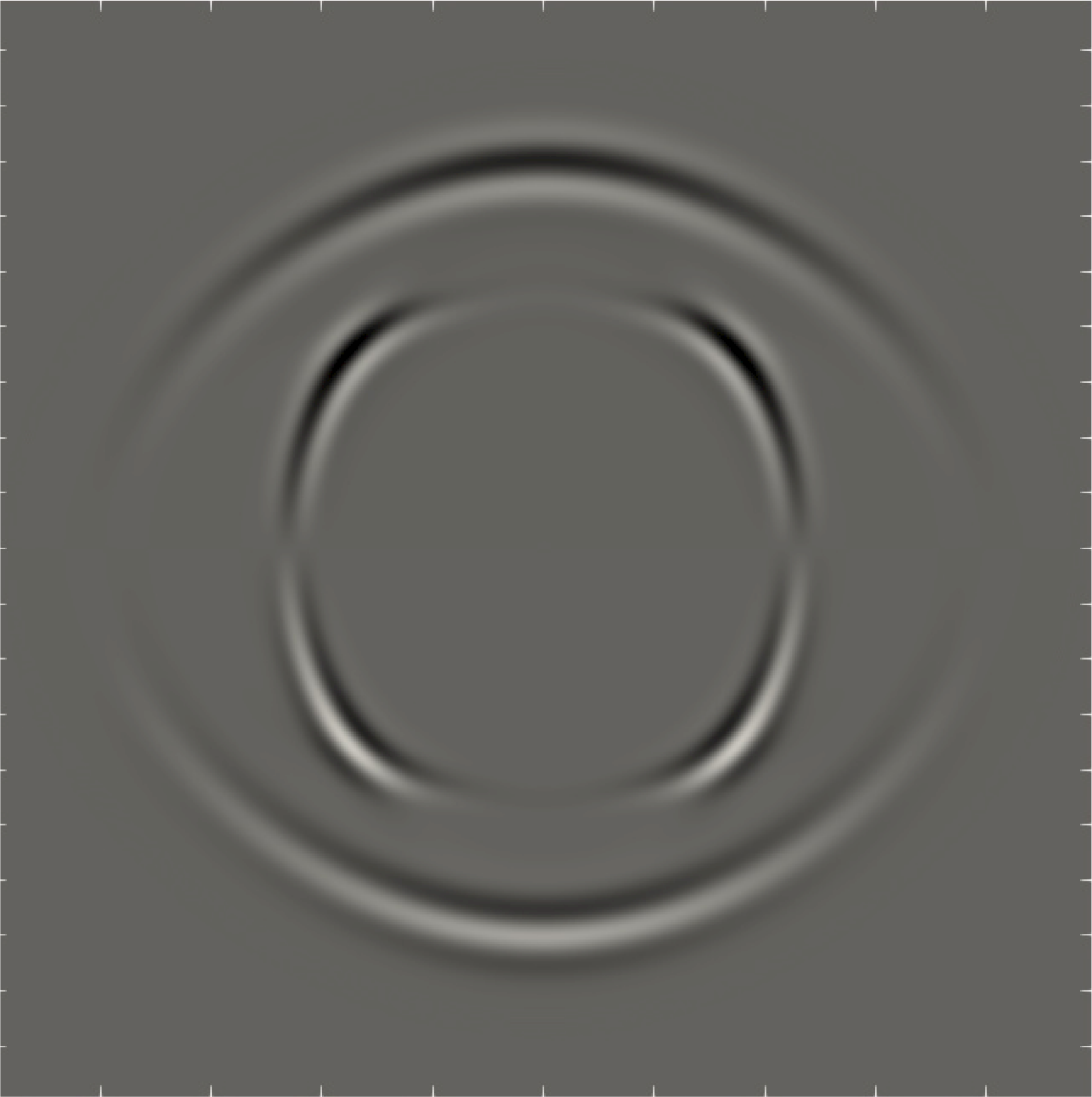};
\end{axis}
\end{tikzpicture}
}
\caption{Snapshots of the centre of mass particle velocity in orthotropic sandstone, computed at $t=1.56~\text{ms}$, where (a) and (b) corresponds to $\eta=0$, and (c) and (d) corresponds to $\eta \ne 0$. The central frequency of the forcing function is  3730 Hz. The solution is computed using polynomials of degree $N=4$. }
\end{figure}

\subsubsection{Epoxy-glass material}
Snapshots of the $x-$ and $z-$ components of the center of mass particle velocity in the epoxy-glass porous medium are shown in Figure 6. Figures 6(a) and (b) correspond to the inviscid case $(\eta=0)$, and (c) and (d) to the viscid case $(\eta \ne 0)$. The central frequency is $ f_c=3135~\text{Hz}$ and the propagation time is $1.8~\text{ms}$. It is worth noting the cuspidal triangles of $\text{S}$ and $\text{P}_ \text{s} $, which are phenomena typical in anisotropic materials.  At $45^o$, the polarization of the $\text{P}_ \text{s}$ mode wave is almost horizontal, which confirms the results shown in Figure 3(b) of \cite{carcione1996}.

\begin{figure}
\centering
\subfloat[Epoxy-glass, $b_x$ with $\eta = 0$]{
\begin{tikzpicture}
\begin{axis}[enlargelimits=false, axis on top, axis equal image, xlabel=x (m), ylabel=y (m), width=0.6\textwidth]

\addplot graphics [xmin=0,xmax=20,ymin=0,ymax=20] {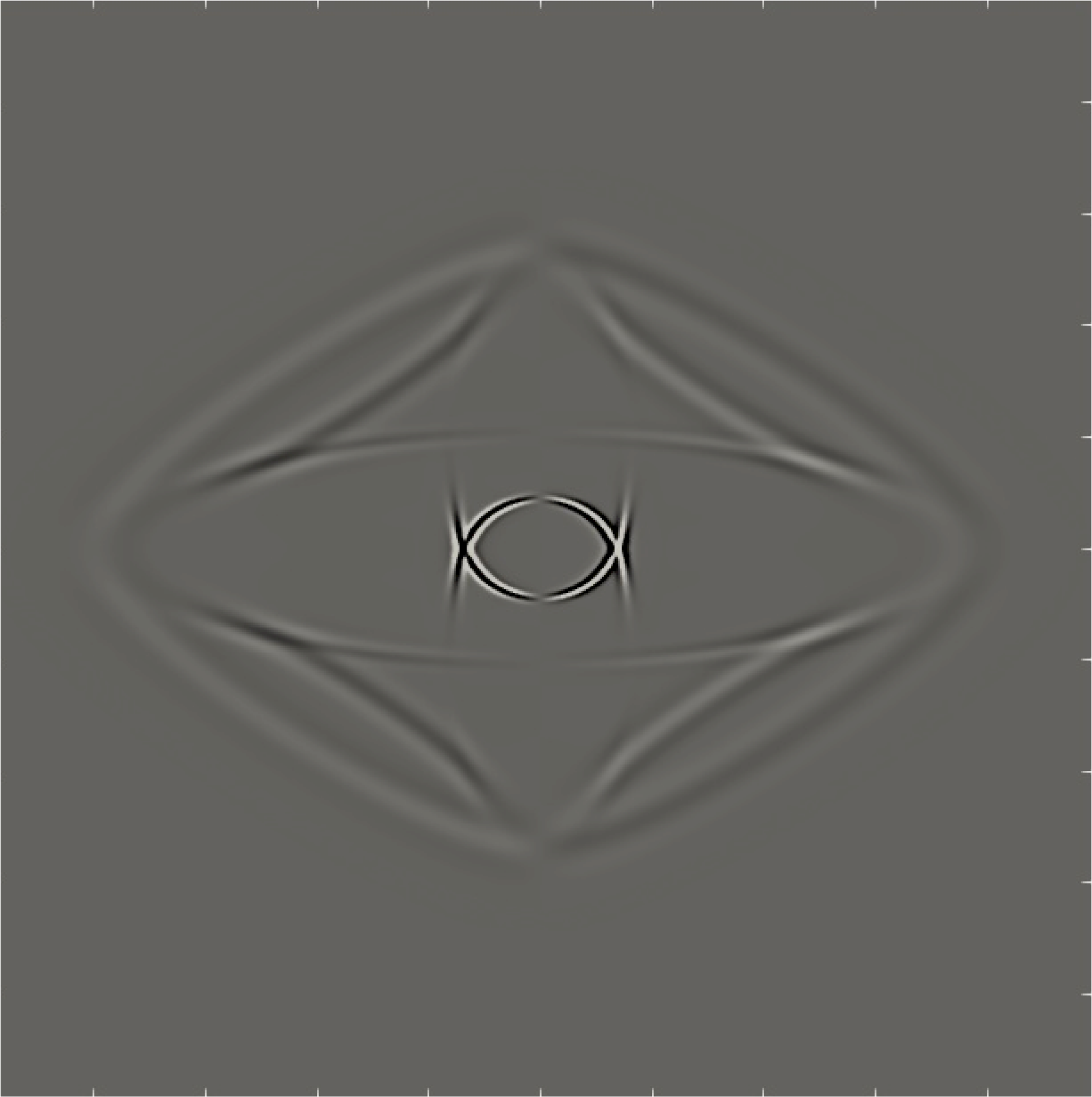};
\end{axis}
\end{tikzpicture}
}
\subfloat[Epoxy-glass, $b_z$ with $\eta = 0$]{
\begin{tikzpicture}
\begin{axis}[enlargelimits=false, axis on top, axis equal image, xlabel=x (m), ylabel=y (m), width=0.6\textwidth]

\addplot graphics [xmin=0,xmax=20,ymin=0,ymax=20] {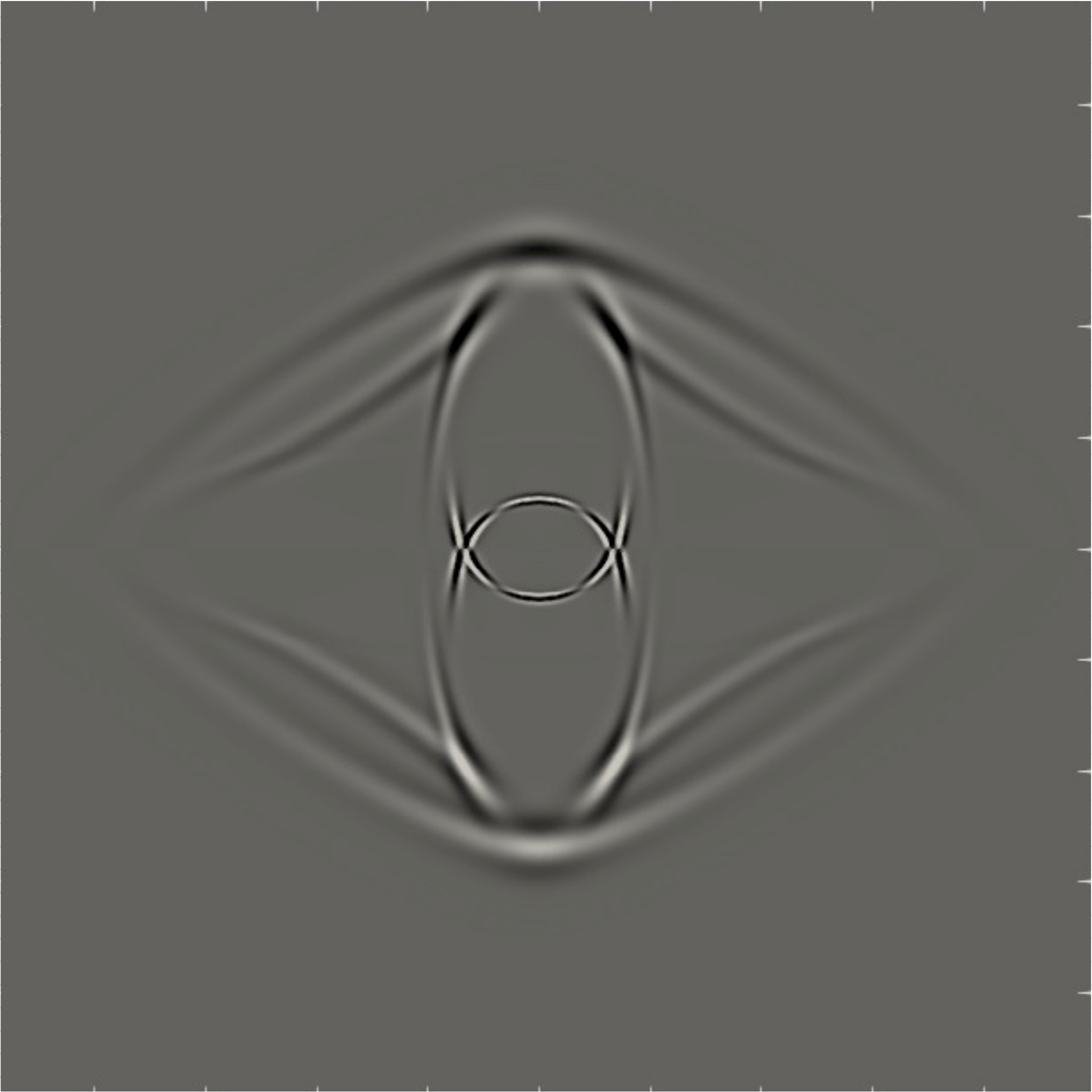};
\end{axis}
\end{tikzpicture}
}\\
\subfloat[Epoxy-glass, $b_x$ with $\eta \ne 0$]{
\begin{tikzpicture}
\begin{axis}[enlargelimits=false, axis on top, axis equal image, xlabel=x (m), ylabel=y (m), width=0.6\textwidth]

\addplot graphics [xmin=0,xmax=20,ymin=0,ymax=20] {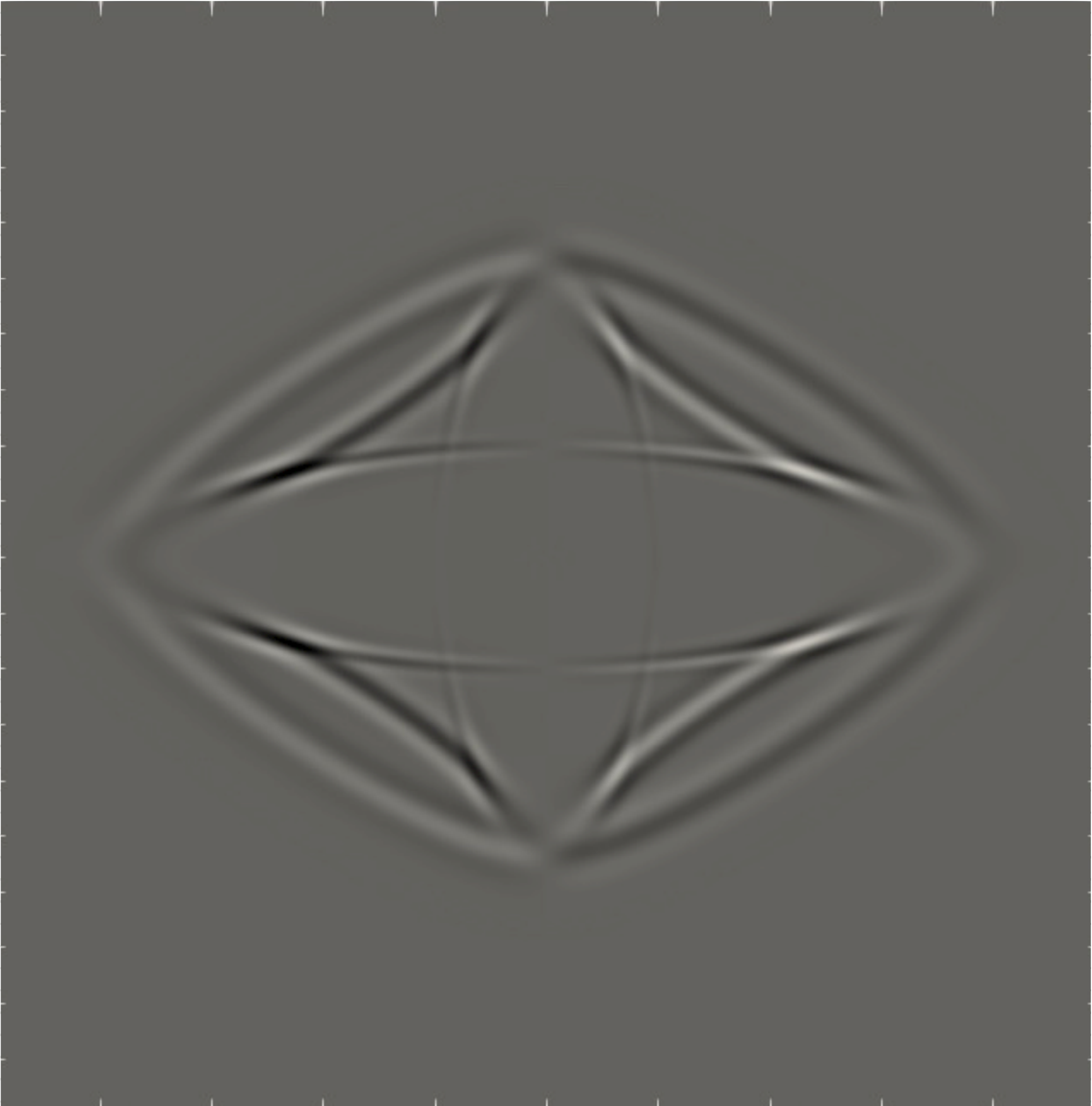};
\end{axis}
\end{tikzpicture}
}
\subfloat[Epoxy-glass, $b_z$ with $\eta \ne 0$]{
\begin{tikzpicture}
\begin{axis}[enlargelimits=false, axis on top, axis equal image, xlabel=x (m), ylabel=y (m), width=0.6\textwidth]

\addplot graphics [xmin=0,xmax=20,ymin=0,ymax=20] {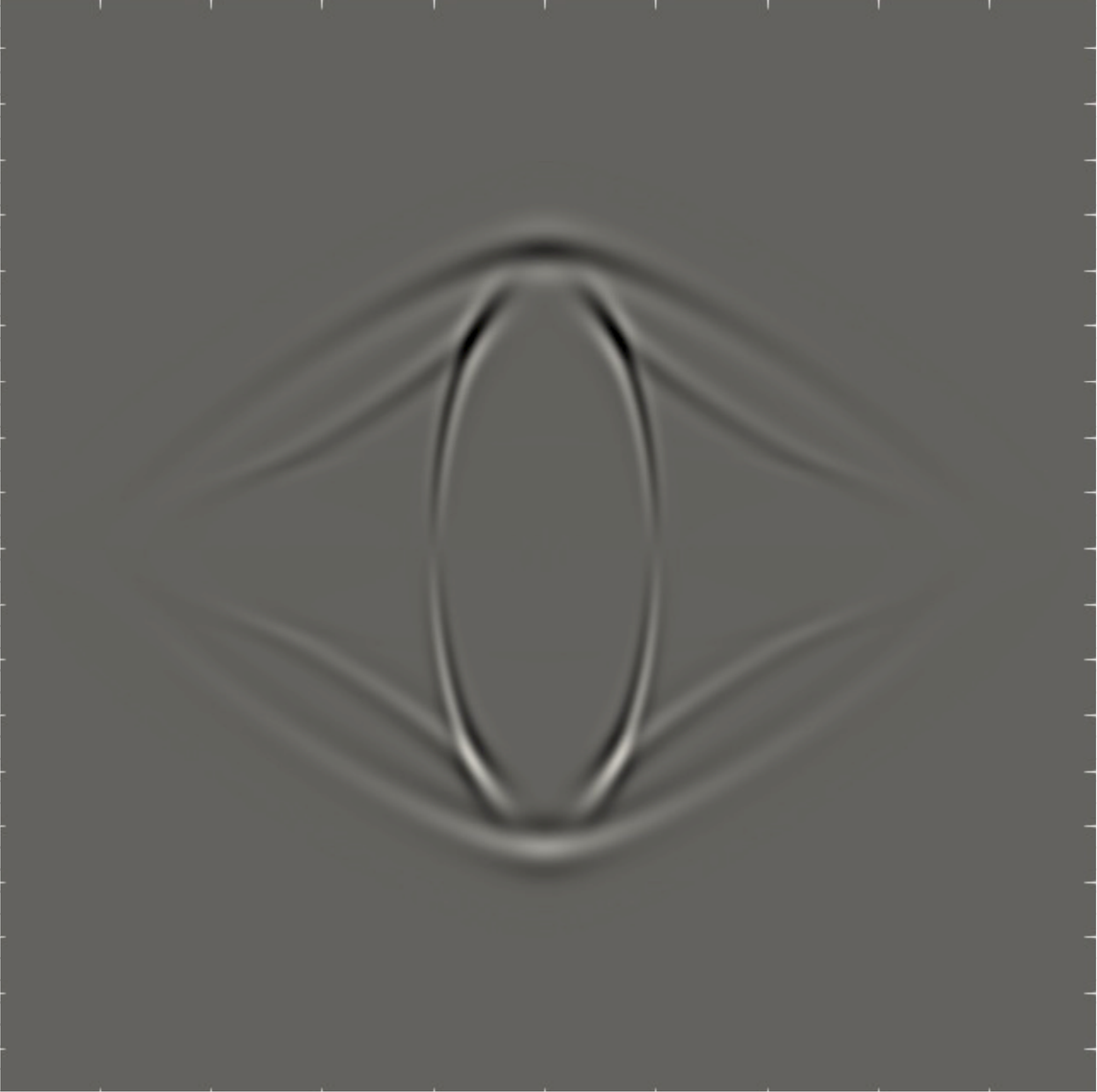};
\end{axis}
\end{tikzpicture}
}
\caption{Snapshots of the centre of mass particle velocity in epoxy-glass, computed at $t=1.8~\text{ms}$, where (a) and (b) corresponds to $\eta=0$, and (c) and (d) corresponds to $\eta \ne 0$. The central frequency of the forcing function is  3135 Hz. The solution is computed using polynomials of degree $N=4$. }
\end{figure}

\subsubsection{Isotropic-isotropic and isotropic-anisotropic layered model}
In this example, we illustrate the effect of an interface between two layers of porous media. 
We constructed two models, each with two layers.  In the first model, the top and bottom layer are isotropic and made of shale and sandstone.  In second model, the sandstone is replaced by orthotropic sand stone. The detailed material properties are given in Table 1. These models are considered to be filled with brine $(\eta=0)$. The size of the computational domain is $1400~\text{m} \times 1500~\text{m}$ in the $x$ and $z$ directions, respectively. The minimum edge size of the triangular elements used to mesh the domain is $8~\text{m}$. The point source is located at $(750~\text{m}, 900~\text{m})$ with a Ricker wavelet of frequency 45 Hz. The propagation time for isotropic and anisotropic model are $0.25~\text{s}$ and $0.22~\text{s}$, respectively. The simulation is performed using polynomials of degree $N=4$. Snapshots of the $z$ component of the center of mass particle velocity are shown in Figures 7a and 7b for inviscid $(\eta=0)$ and viscid fluid $(\eta \ne 0)$, respectively. Figure 7 clearly shows the direct, reflected, and transmitted wavefront, corresponding to all three modes. The slow P wave is more prominent in the shale. The effect of anisotropy on all three modes is clearly seen as wavefronts moves with different phase velocities.

\begin{figure}
\centering
\subfloat[Isotropic layered model, $b_z$ ]{
\begin{tikzpicture}
\begin{axis}[enlargelimits=false, axis on top, axis equal image, xlabel=x (km), ylabel=y (km), width=0.6\textwidth]

\addplot graphics [xmin=0,xmax=1.4,ymin=0,ymax=1.5] {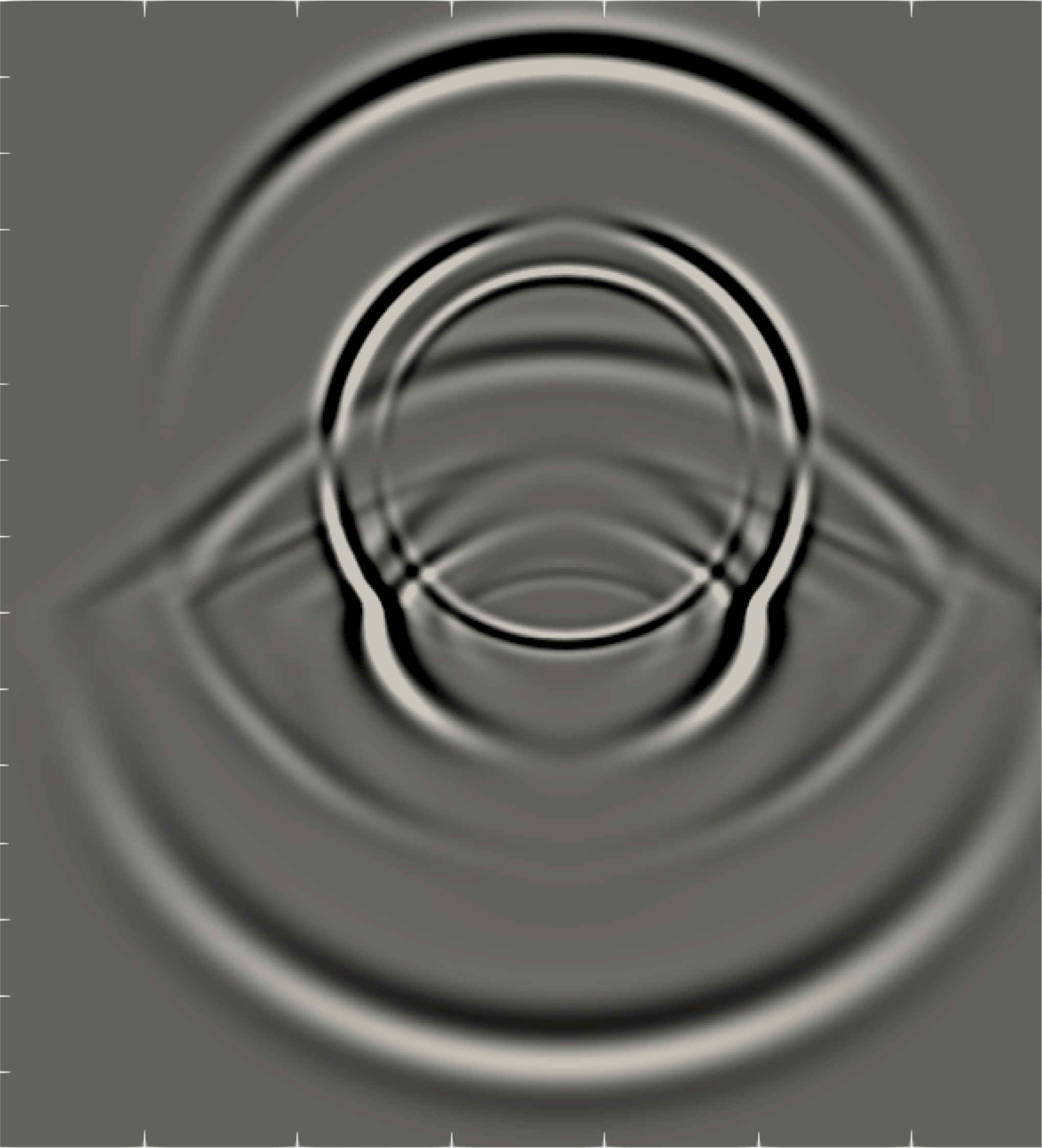};
\end{axis}
\end{tikzpicture}
}
\subfloat[Anisotropic layered model, $b_z$]{
\begin{tikzpicture}
\begin{axis}[enlargelimits=false, axis on top, axis equal image, xlabel=x (km), ylabel=y (km), width=0.6\textwidth]

\addplot graphics [xmin=0,xmax=1.4,ymin=0,ymax=1.5] {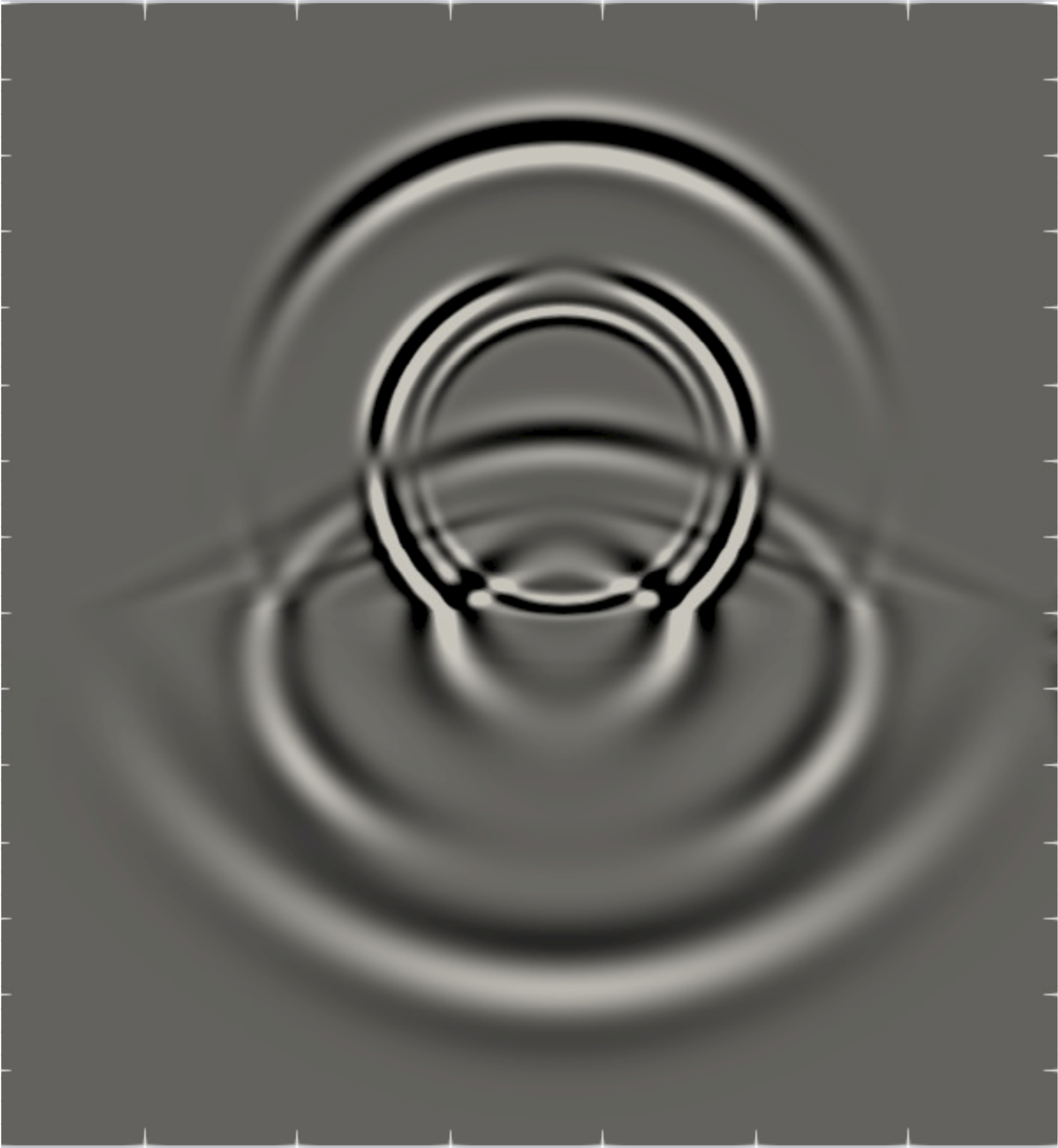};
\end{axis}
\end{tikzpicture}
}

\caption{Snapshots of the $z-$ component of centre of mass particle velocity in (a) isotropic layered model computed at $t=0.25~\text{s}$ (b) anisotropic layered model computed at $t=0.22~\text{s}$. The central frequency of the forcing function is 45 Hz. The solution is computed using polynomials of degree $N=4$.}
\end{figure}

Here, we illustrate the effect of viscosity of fluid on a slow P wave using an isotropic two layer model with the same discretization parameters as those used in Figure 7a. We compute numerical solutions for two cases. In first case, both layers are assumed to be brine filled with $(\eta=0)$ and are represented by inviscid-inviscid model. However, in the second case we assume that the bottom layer is filled with a viscous fluid and represented as inviscid-viscid model. Figure 8a shows the snapshot of the $z$ component of the center of mass particle velocity for the inviscid-inviscid model, which shows a direct, reflected, and transmitted slow P wave. On the other hand, Figure 8b shows a snapshot of the $z$ component of the center of mass particle velocity for the inviscid-viscid model, which shows only a direct and reflected slow P wave. The transmitted slow P wave is missing in the inviscid-viscid model, which confirms Biot's observation on slow P waves for the low frequency regime \cite{biot1956I}.

\begin{figure}
\centering
\subfloat[Inviscid-inviscid layered model, $b_z$  ]{
\begin{tikzpicture}
\begin{axis}[enlargelimits=false, axis on top, axis equal image, xlabel=x (km), ylabel=y (km), width=0.6\textwidth]

\addplot graphics [xmin=0,xmax=1.4,ymin=0,ymax=1.5] {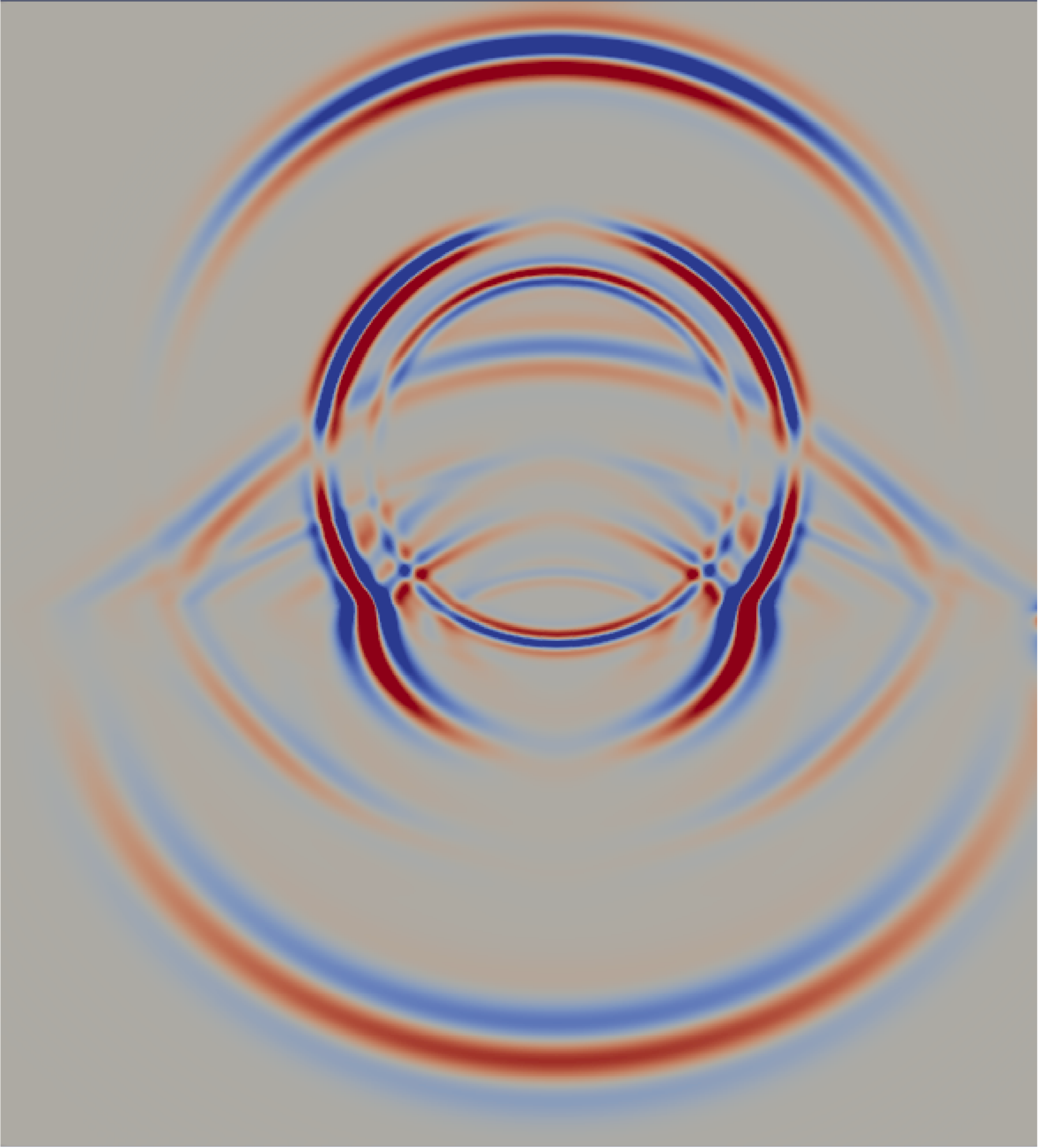};
\end{axis}
\end{tikzpicture}
}
\subfloat[Inviscid-viscid layered model, $b_z$]{
\begin{tikzpicture}
\begin{axis}[enlargelimits=false, axis on top, axis equal image, xlabel=x (km), ylabel=y (km), width=0.6\textwidth]
\addplot graphics [xmin=0,xmax=1.4,ymin=0,ymax=1.5] {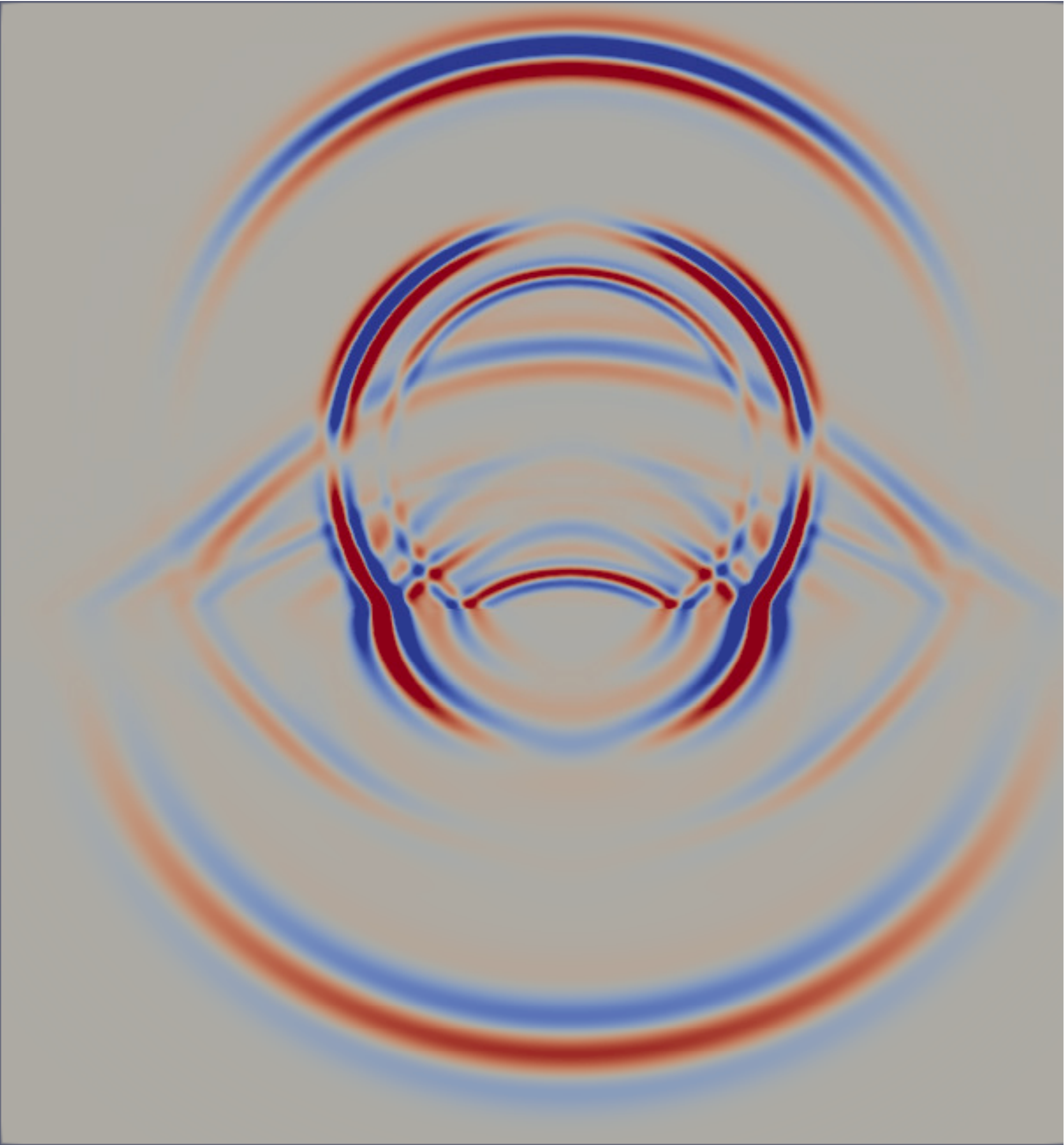};
\end{axis}
\end{tikzpicture}
}

\caption{Snapshots of the $z-$ component of centre of mass particle velocity in (a) Inviscid-inviscid layered model (b) Inviscid-viscid layered model, computed at $t=0.25~\text{s}$. The central frequency of the forcing function is 45 Hz for a viscid case $(\eta=0)$. The solution is computed using polynomials of degree $N=4$.}
\end{figure}

\subsubsection{Piecewise constant vs sub-cell heterogeneities}
To address the effect of micro-heterogenities on poroelastic wave propagation, we illustrate the efficacy of WADG method over a traditional DG method using piecewise constant coefficients on each element.  We modulate the material properties contained in $\bm{Q}^{-1}_s$ and $\bm{Q}^{-1}_v$ with
\[
\rho=1+.5\sin(2\pi \bm{x})\sin(2\pi \bm{y}),
\]
where $\bm{x}$ and $\bm{y}$ are coordinates of of cubature nodes.  We construct piecewise constant approximations to $\bm{Q}^{-1}_s$ and $\bm{Q}^{-1}_v$ by taking the local average over each element.

Figure 9a and 9b show piecewise smooth and high order heterogeneous approximations of bulk density, respectively. To demonstrate the effect of subelement variations on wave propagation, we first compute the solutions approximated by polynomials of degree $N=2$ with $K_{\rm 1D}=64$, where $K_{\rm  1D}$ denotes the number of elements in one direction. These solutions are shown in Figures 10a (piecewise constant) and 10b (WADG).  The solutions are similar, as the lower resolution due to the piecewise constant approximation of media is offset by the larger number of elements.  Next, we compute the solutions for $N=4$ and $K_{1D}=32$. The effects of subelement variations are now clearly marked by spurious reflections in Figures 10c and 10d. These effects are more prominent for higher orders as shown in Figures 10e and f, which display the piecewise constant and WADG solutions for $N=8$ and $K_{1D}=16$.

\begin{figure}
\centering
\subfloat[Piecewise smooth bulk-density $\bm{\rho}$ ]{
\includegraphics[width=0.5\textwidth]{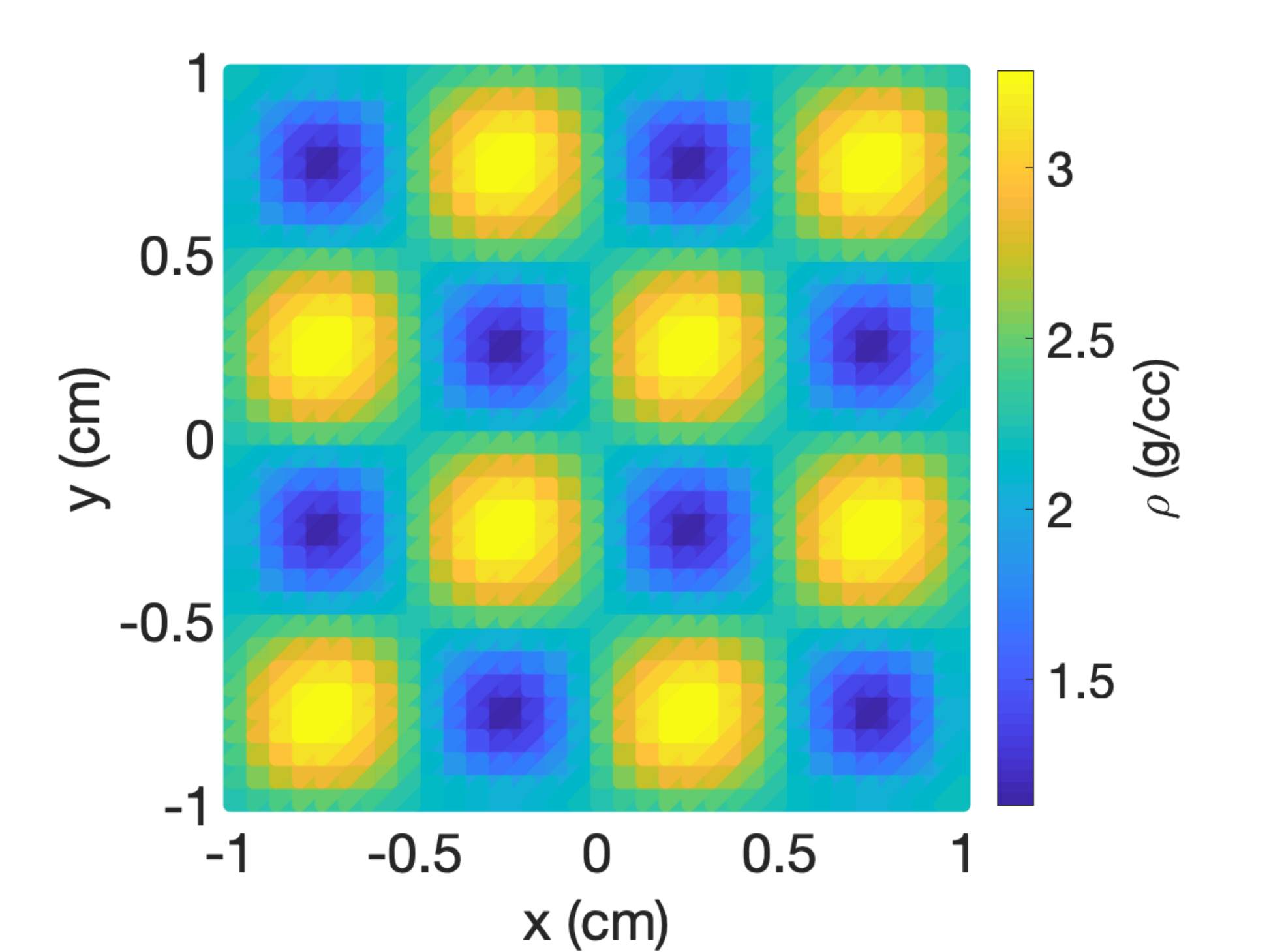}
}
\subfloat[heterogeneous bulk-density $\bm{\rho}$ ]{
\includegraphics[width=0.50\textwidth]{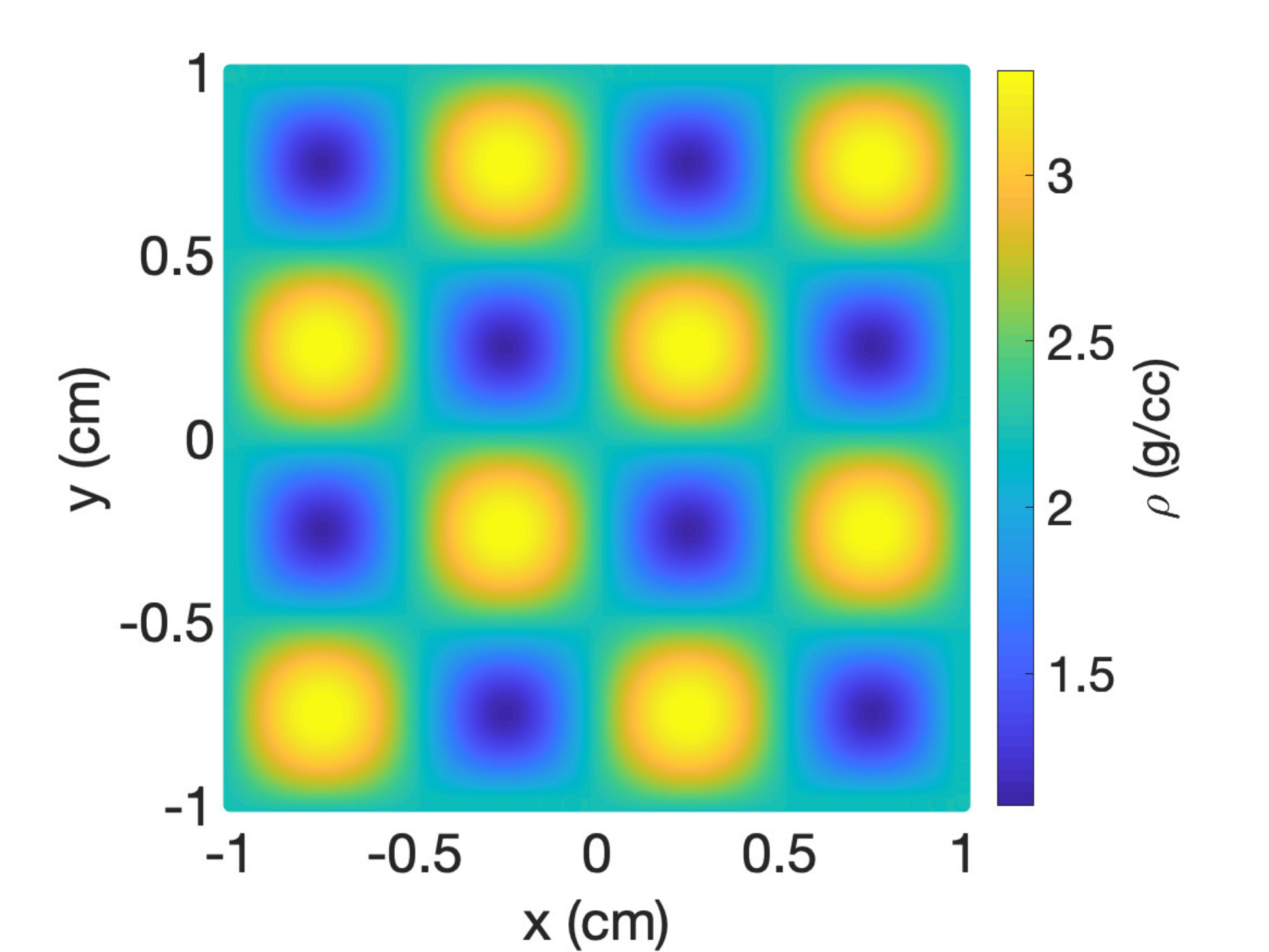}
}
\caption{Heterogeneous model of bulk density $(\rho)$ (a) obtained from piecewise approximation with high order sub-cell variations, used for the piecewise constant solution (b) high order sub-cell heterogeneous distribution, used for WADG solution.}
\end{figure}

\begin{figure}
\centering
\subfloat[$N=2,~K_{1D}=64,$ and piecewise constant ]{
\includegraphics[width=0.5\textwidth]{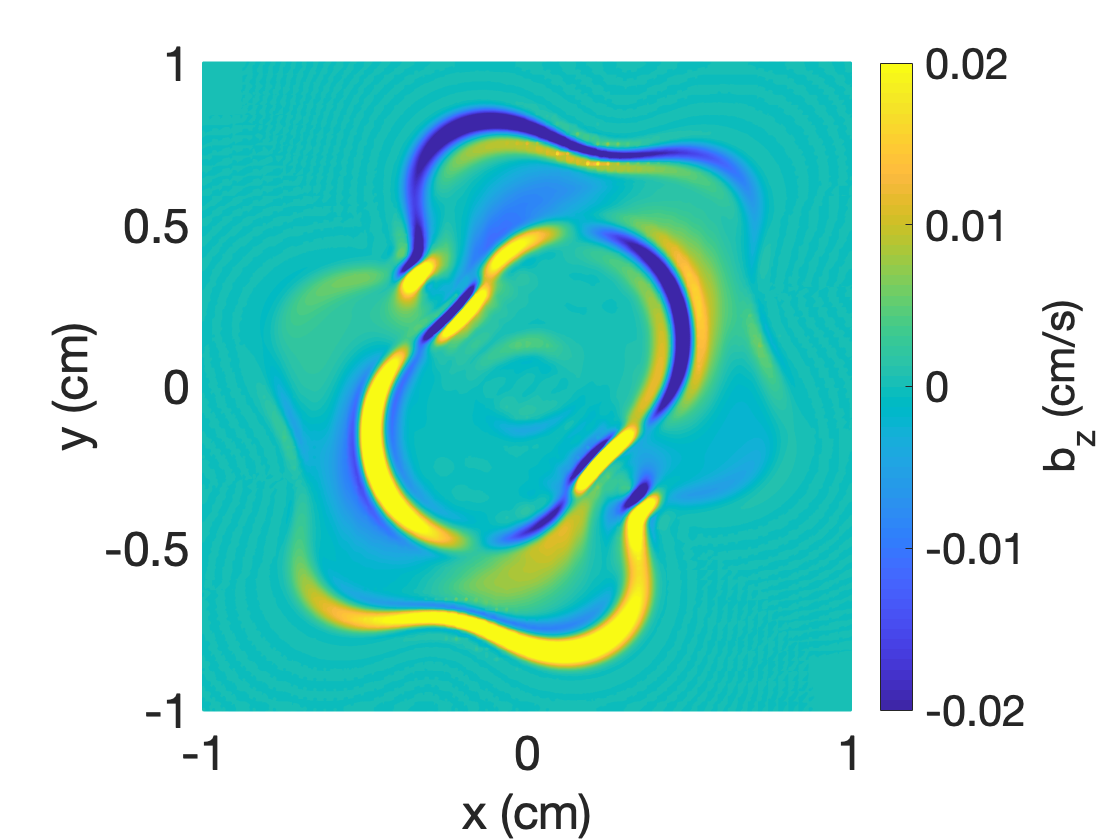}
}
\subfloat[$N=2,~K_{1D}=64,$ and WADG ]{
\includegraphics[width=0.50\textwidth]{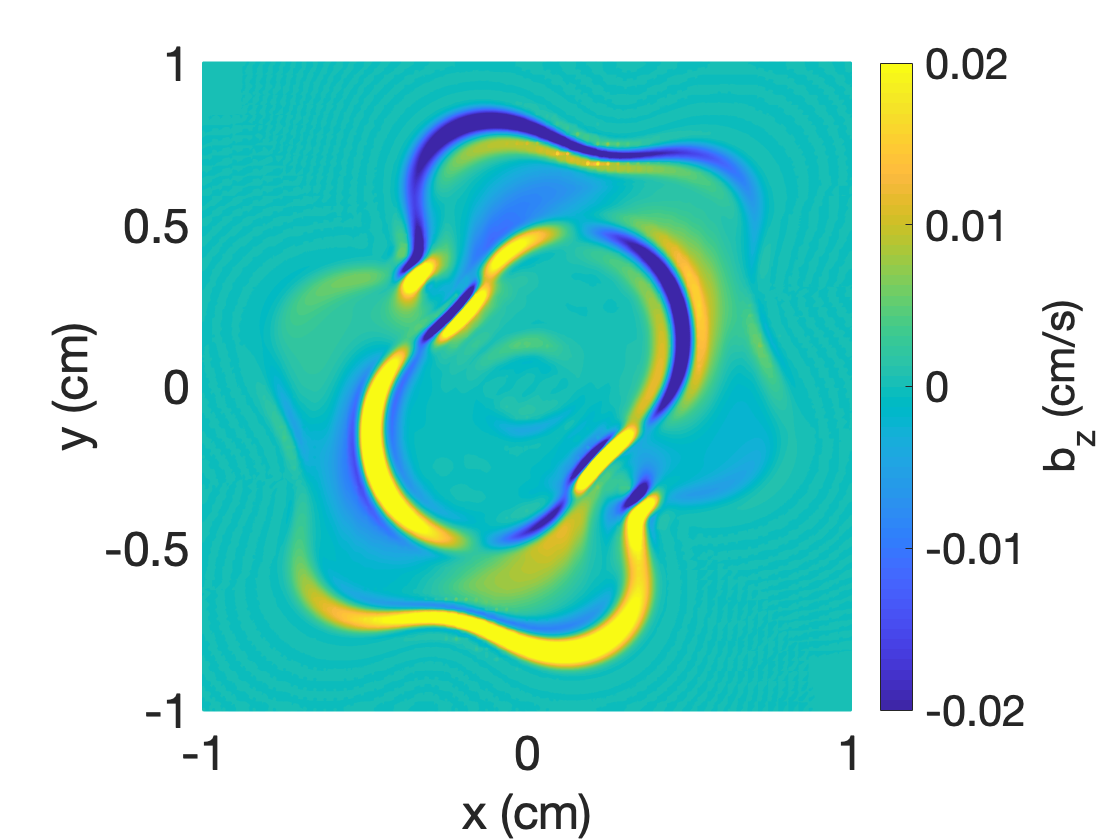}
}\\
\subfloat[$N=4,~K_{1D}=32,$ and piecewise constant ]{
\includegraphics[width=0.5\textwidth]{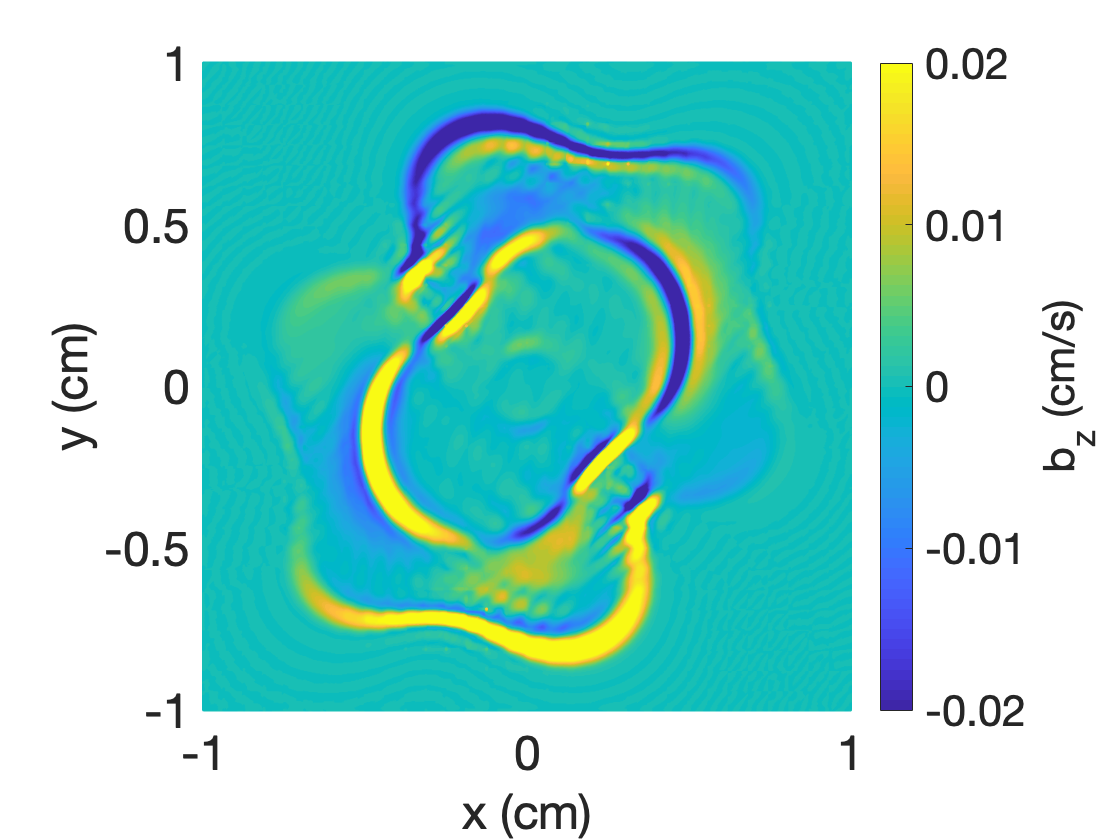}
}
\subfloat[$N=4,~K_{1D}=32,$ and WADG ]{
\includegraphics[width=0.50\textwidth]{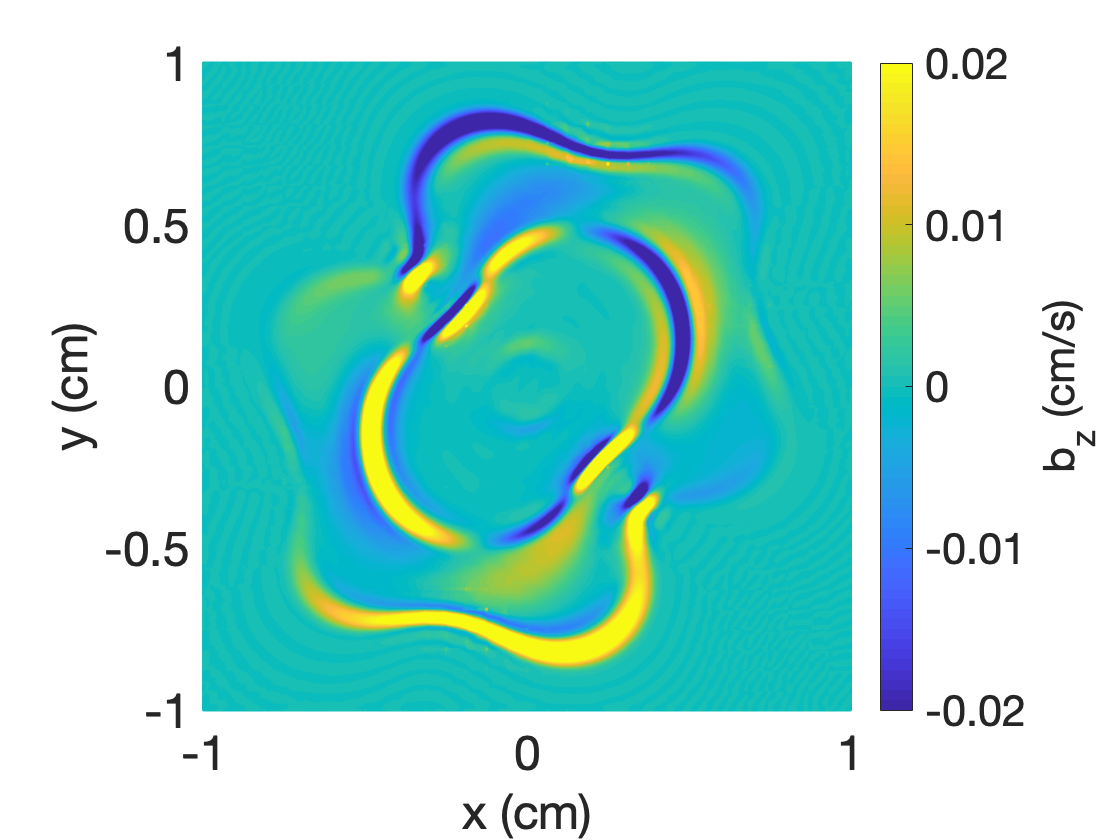}
}\\
\subfloat[$N=8,~K_{1D}=16,$ and piecewise constant ]{
\includegraphics[width=0.5\textwidth]{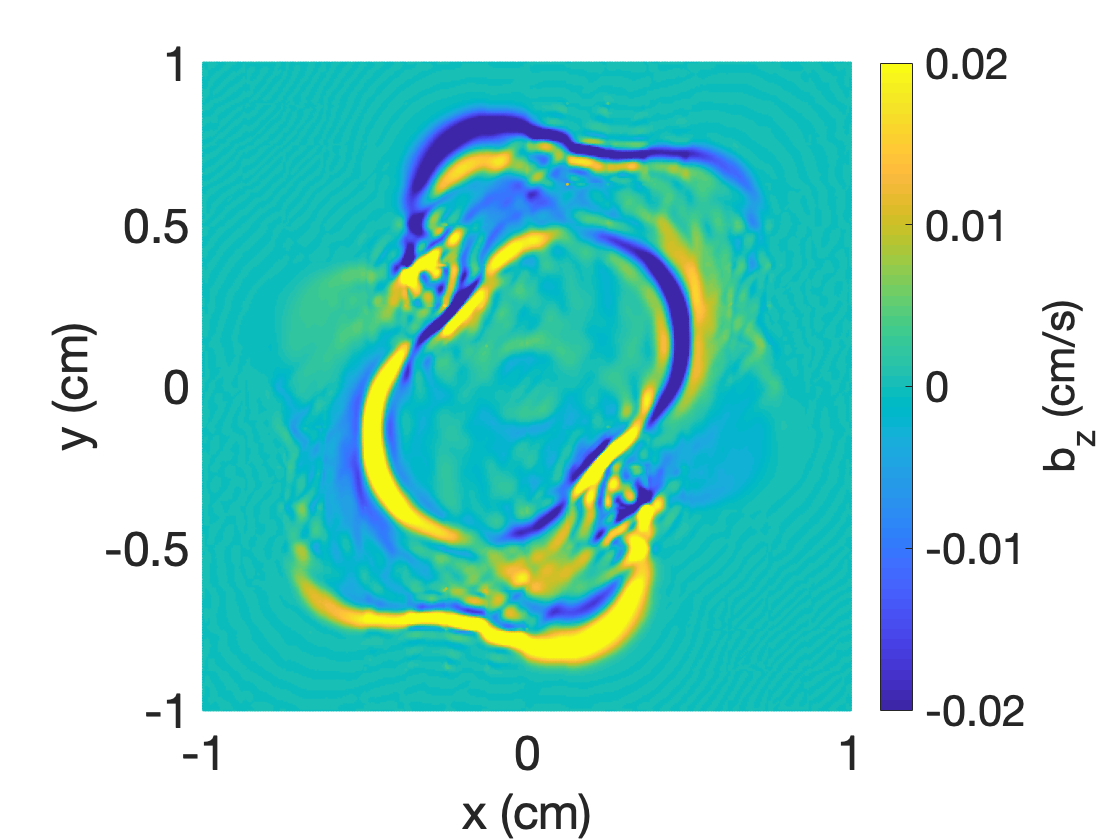}
}
\subfloat[$N=8,~K_{1D}=16,$ and WADG ]{
\includegraphics[width=0.50\textwidth]{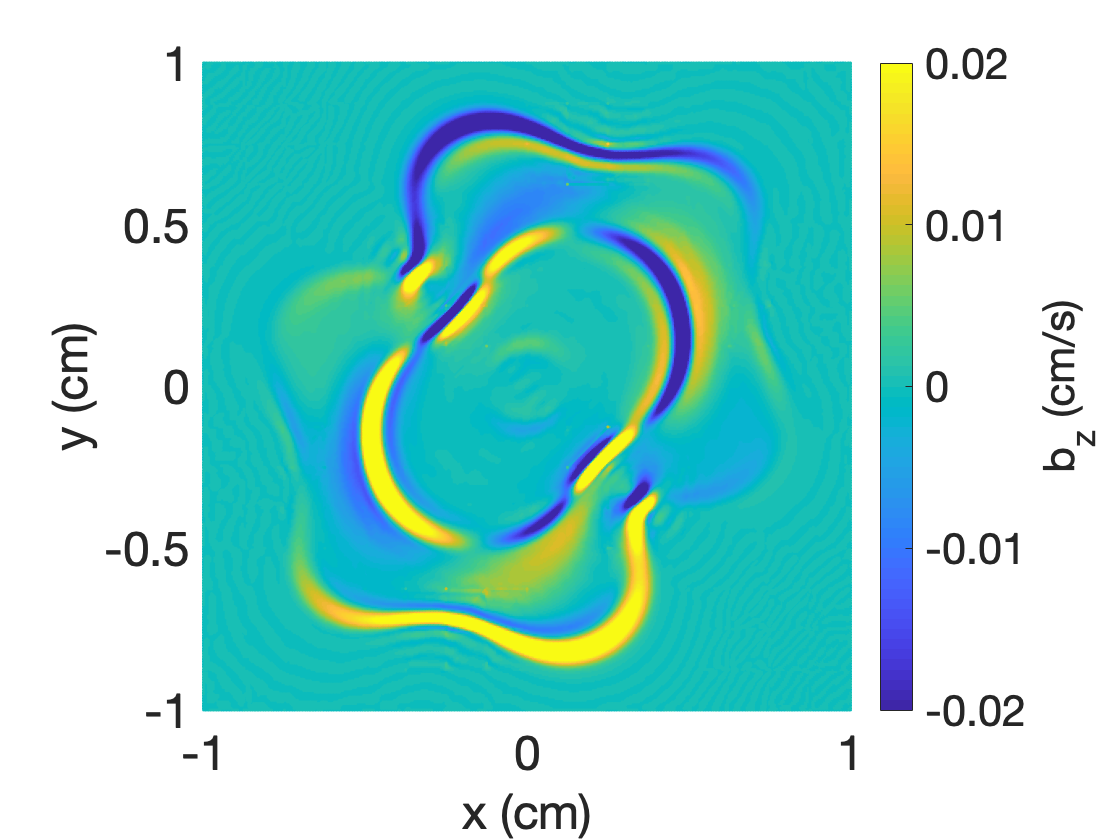}
}
\caption{A comparison between the solutions obtained from piecewise constant media and sub-cell heterogeneities treated using WADG.}
\end{figure}

\begin{table*}
	\caption{Material properties used for comparison of solution with SPECFEM implementations \cite{morency2008}}
	\centering
	\medskip
	\begin{tabular}{c c c c } 
		\hline
		Properties  & Medium I  & Medium II & Medium III \\ [0.5ex] 
		\hline
		$\rho_s~$ $($kg/m$^3$$)$ & 2650 & 2200 & 2650 \\
		$\rho_f~$ $($kg/m$^3$$)$ & 880 & 950 & 750\\
		$\phi~$ & 0.1 & 0.4 & 0.2 \\
		$T_1~$ & 2 & 2 & 2 \\
		$T_3~$ & 2 & 2 & 2\\
		$K_s~$(GPa) & 12.2 & 6.9 & 6.9  \\
		$K_f~$(GPa) & 1.985 & 2.0 & 2.0  \\
		$\kappa_{fr}~$(GPa)$^{*}$ & 9.6 & 6.7 & 6.7\\ 
		$\mu_{fr}^{**}~$(GPa)$^{**}$ & 5.1 & 3.0 & 3.0\\
		${\eta}~$$(10^{-3}$ Kg/m.s)  & 0  & 0 & 0 \\ [1ex]
		\hline
		$^*$: Frame bulk modulus,~$^{**}$: Frame shear modulus.
	\end{tabular}
\end{table*}

\subsubsection{Validation of numerical solution}
Here, we show comparisons between numerical solutions of the poroelastic wave equation obtained from the proposed DG method and from the spectral element method \cite{morency2008}. Morency and Tromp \cite{morency2008} used the spectral element method to solve the poroelastic wave equation in second order displacement form for 2D isotropic porous materials. First, we perform a validation study in a homogeneous porous medium. The domain dimensions are $1000~\text{m}\times1000~\text{m}$. The domain is discretized with triangular and quadrilateral elements with minimum edge length $5~\text{m}$ for the DG and spectral element method, respectively. 

Mornecy and Tromp implemented a pure dilatational explosive source as an external forcing function. To match the assumption of the dilatational explosive source, we have set the forcing functions corresponding to $\tau_{11}$ and $\tau_{22}$ to be equal, while the forcing function corresponding to $\tau_{12}$ is set to zero since a pure dilatational source does not radiate shear waves. It was also assumed in \cite{morency2008} that the fluid and the solid contain equal energy, and thus the forcing function corresponding to the pressure variable $p$ is the same as for $\tau_{11}$ and $\tau_{12}$. 

The time domain response of the forcing function is a Ricker wavelet with a central frequency of $f_c=30~\text{Hz}$. The material properties of the medium are given in column 2 of Table 2. In [3], the material properties were represented using Lame constants $\kappa_{fr}$, $\mu_{fr}$.  These can be converted to compliance coefficients as follows: 
\[
c_{55}=\mu_{fr},~~c_{11}=\kappa_{fr} + \myfrac{4}{3}\mu_{fr},~~c_{12}=c_{13}=\kappa_{fr}-\myfrac{2}{3}\mu_{fr}.
\]

Figure 11a shows a snapshot of the solid particle velocity $v_3$, 
computed at $t=0.14~\text{s}$. $\text{P}_\text{f}$ and $\text{P}_\text{s}$ represent the modes corresponding to fast and slow P waves, respectively. The point source is located at $\bm{x}_0=(500~\text{m},~300~\text{m})$, marked by a white star in Figure 11a. Figure 11b shows a comparison between the normalized solid vertical particle velocity $(v_3)$ obtained from both the DG method and spectral-element method. The time history of the solution is recorded at $(600~\text{m},~400~\text{m})$, denoted by the red diamond in Figure 11a. The agreement between the solutions is good. The root mean square error (RMSE) between the solutions is $8.68\text{e}-03$. A comparison between the solution obtained from the DG method and a pseudo-spectral method \cite{carcione1996} is also shown in Figure 11c. The comparison shows a good agreement between the solutions except at the onset time $0.06~\text{s}$.  The RMS misfit between the solution is $3.18\text{e}-02$.

Next, we validate numerical solutions obtained from the DG method with a spectral-element method for a heterogenous medium. A two layer medium is constructed with material properties given in column 2 (top layer) and 3 (bottom layer) of Table 2. The size of the computational domain is $(4800~\text{m} \times 4800~\text{m})$ and discretized with the triangular and quadrilateral elements for DG and spectral-element method, respectively. The minimum size of the edge of the element is taken as $200~\text{m}$. Figure 12a represents the snapshot of solid particle velocity at $t=0.73~\text{s}$. The pure dilatational point source is located at $\bm{x}_0=(1600~\text{m}, 2900~\text{m})$, denoted by white star in Figure 12a, and implemented in same way as discussed for homogeneous medium. Figure 12b shows a comparison between solutions obtained from the DG and spectral-element methods.  The time history of the solutions is recorded at $(2000~\text{m}, 1867~\text{m})$, marked by a red diamond in Figure 12a.  The comparison shows a good overall agreement, though some discrepancies are present around time $T=1$.  The RMS misfit between the solutions is $6.65\text{e}-02$.

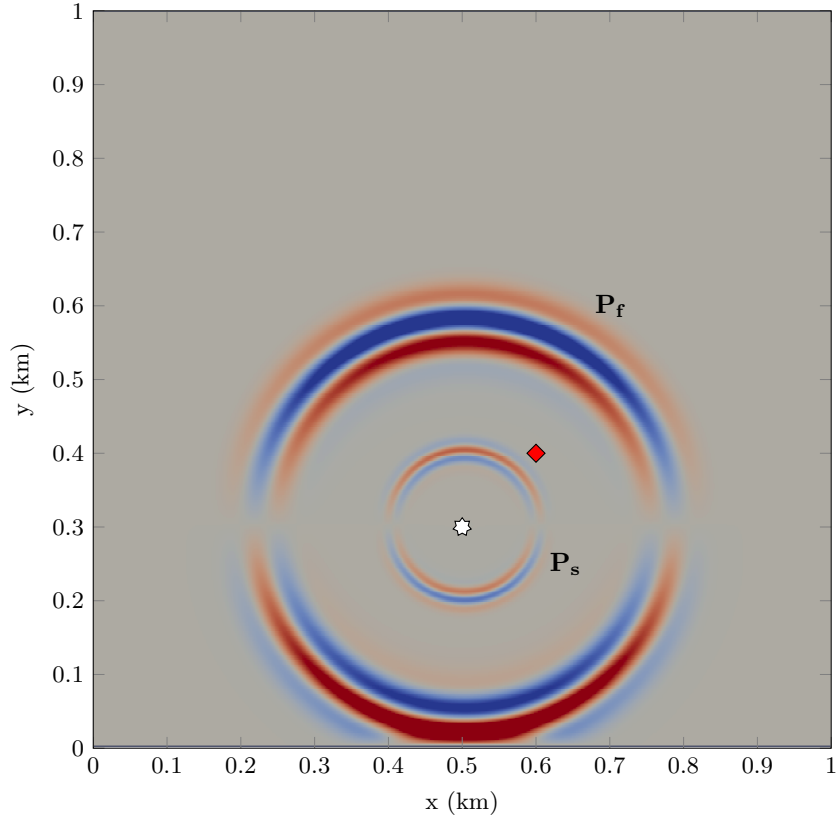
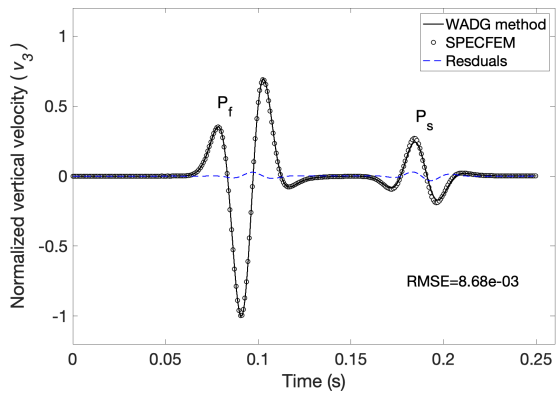
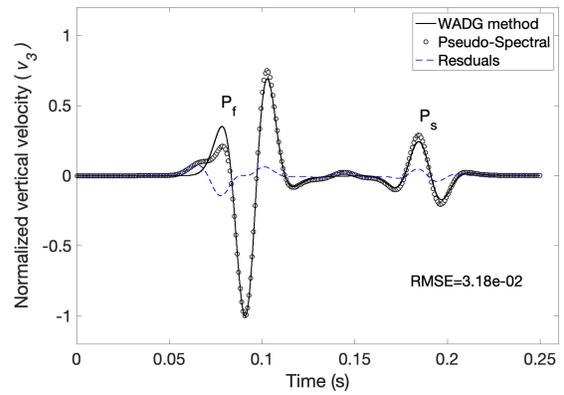
\begin{figure}
\centering
\subfloat[Snapshot of vertical solid particle velocity, $v_3$ at $t=0.14~\text{s}$  ]{
\begin{tikzpicture}
\begin{axis}[enlargelimits=false, axis on top, axis equal image, xlabel=x (km), ylabel=y (km), width=0.8\textwidth]

\addplot graphics [xmin=0,xmax=1,ymin=0,ymax=1] {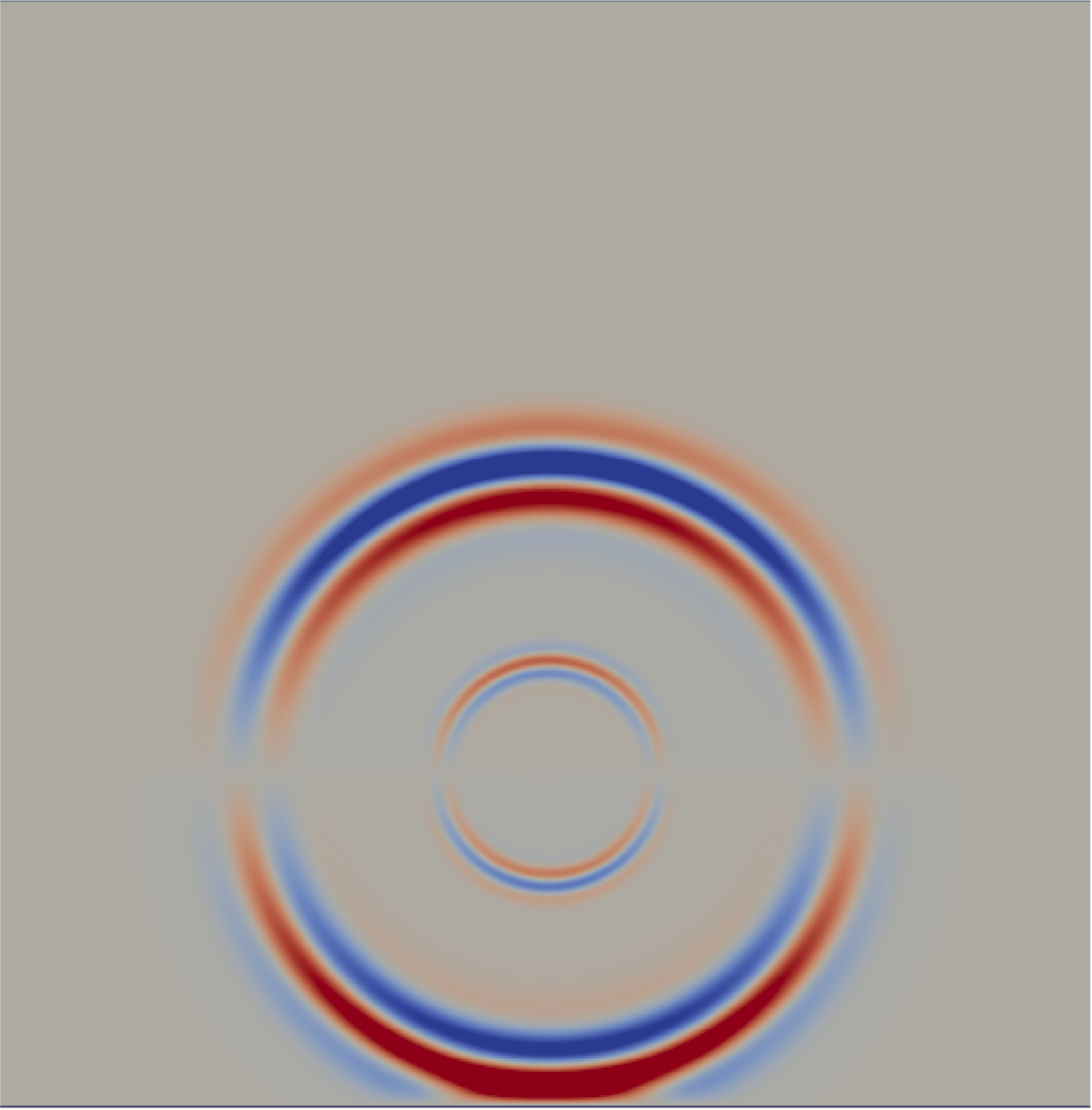};
\node at (axis cs:0.5,0.3)[star,star points=7,draw=black, fill=white, inner sep=0pt,minimum size=7pt]  {};
\node at (axis cs:0.6,0.4)[diamond,draw=black, fill=red, inner sep=0pt,minimum size=7pt] {};
\node at (axis cs:0.7,0.6)[] {$\textbf{P}_\textbf{f}$};
\node at (axis cs:0.64,0.25)[] {$\textbf{P}_\textbf{s}$};
\end{axis}
\end{tikzpicture}
}\\
\subfloat[SPECFEM vs WADG]{
\includegraphics[width=0.5\textwidth]{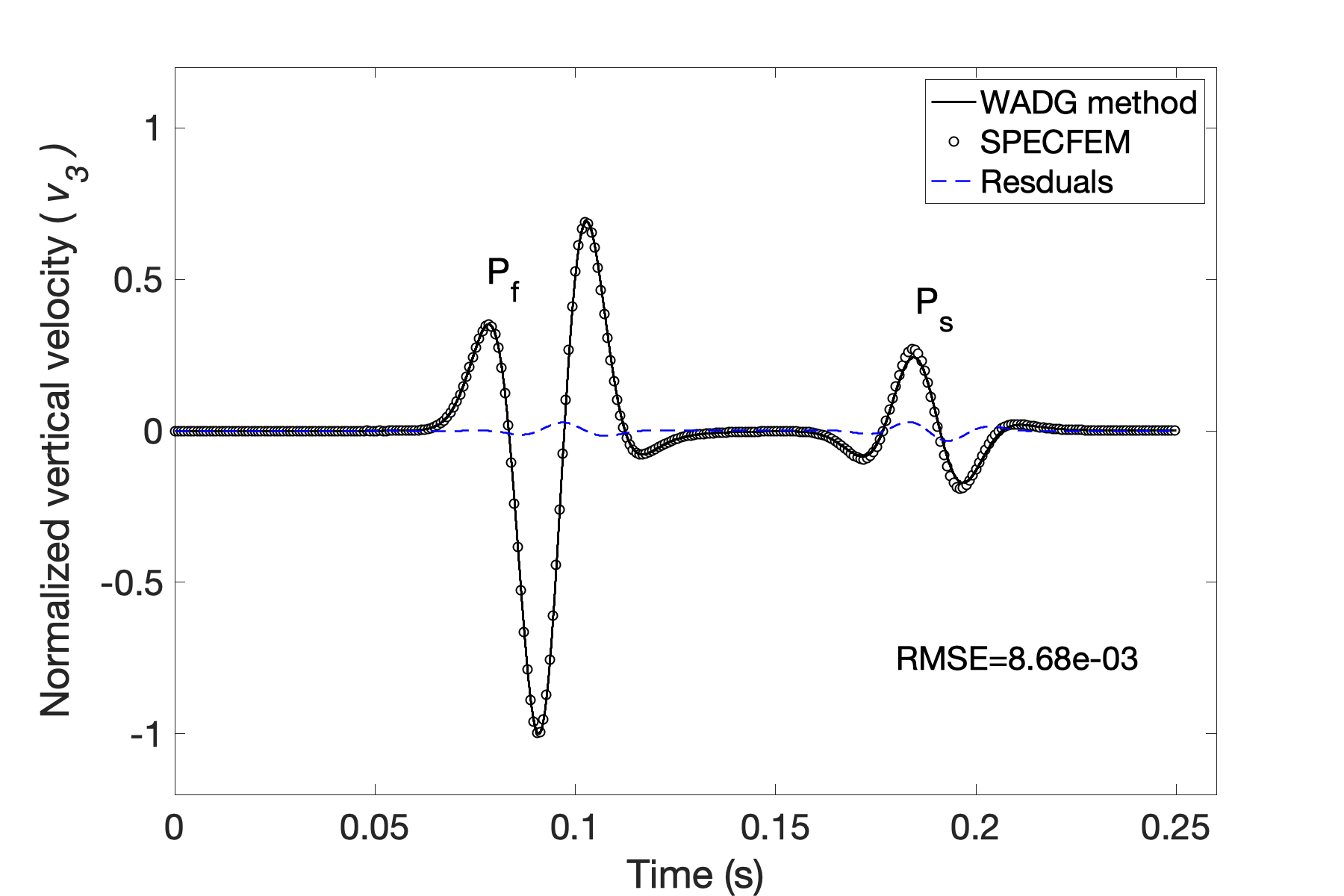}
}
\subfloat[Pseudo-spectral method vs WADG]{
\includegraphics[width=0.5\textwidth]{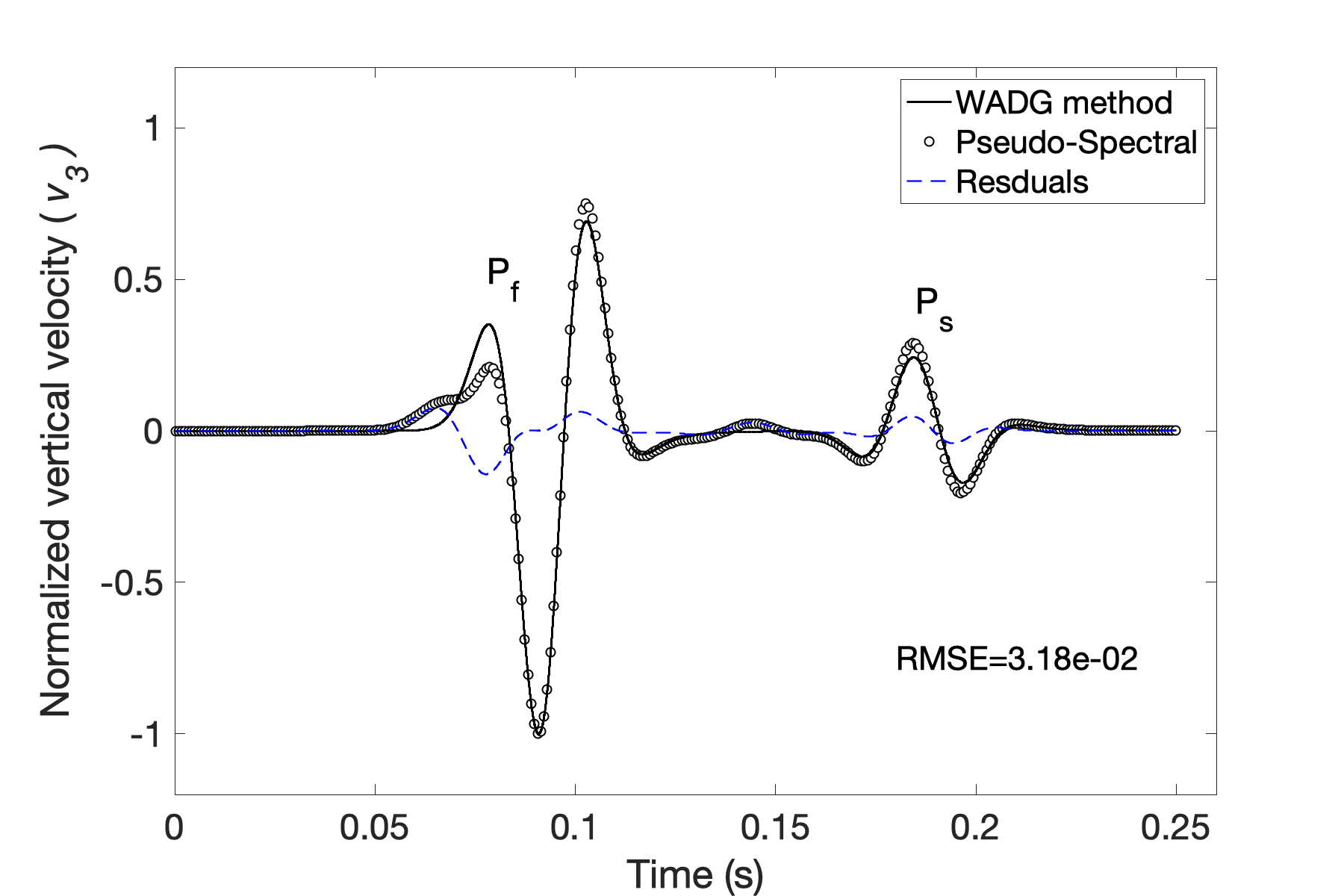}
}

\caption{A comparison between the solutions obtained from the proposed DG method and spectral-element method in a homogeneous porous medium with (a) showing the snapshot of the solutions obtained from the DG method computed with polynomials of degree $N=4$. Figures (b) and (c)  show comparisons between the solutions obtained from the DG/spectral element methods and the DG/pseudo-spectral methods, respectively. $\text{P}_\text{f}$ and $\text{P}_\text{s}$ represent fast and slow P wave, respectively.}
\end{figure}

\begin{figure}
\centering
\subfloat[Snapshot of vertical solid particle velocity, $v_3$ at $t=0.73~\text{s}$  ]{
\begin{tikzpicture}
\begin{axis}[enlargelimits=false, axis on top, axis equal image, xlabel=x (km), ylabel=y (km), width=0.8\textwidth]

\addplot graphics [xmin=0,xmax=4.8,ymin=0,ymax=4.8] {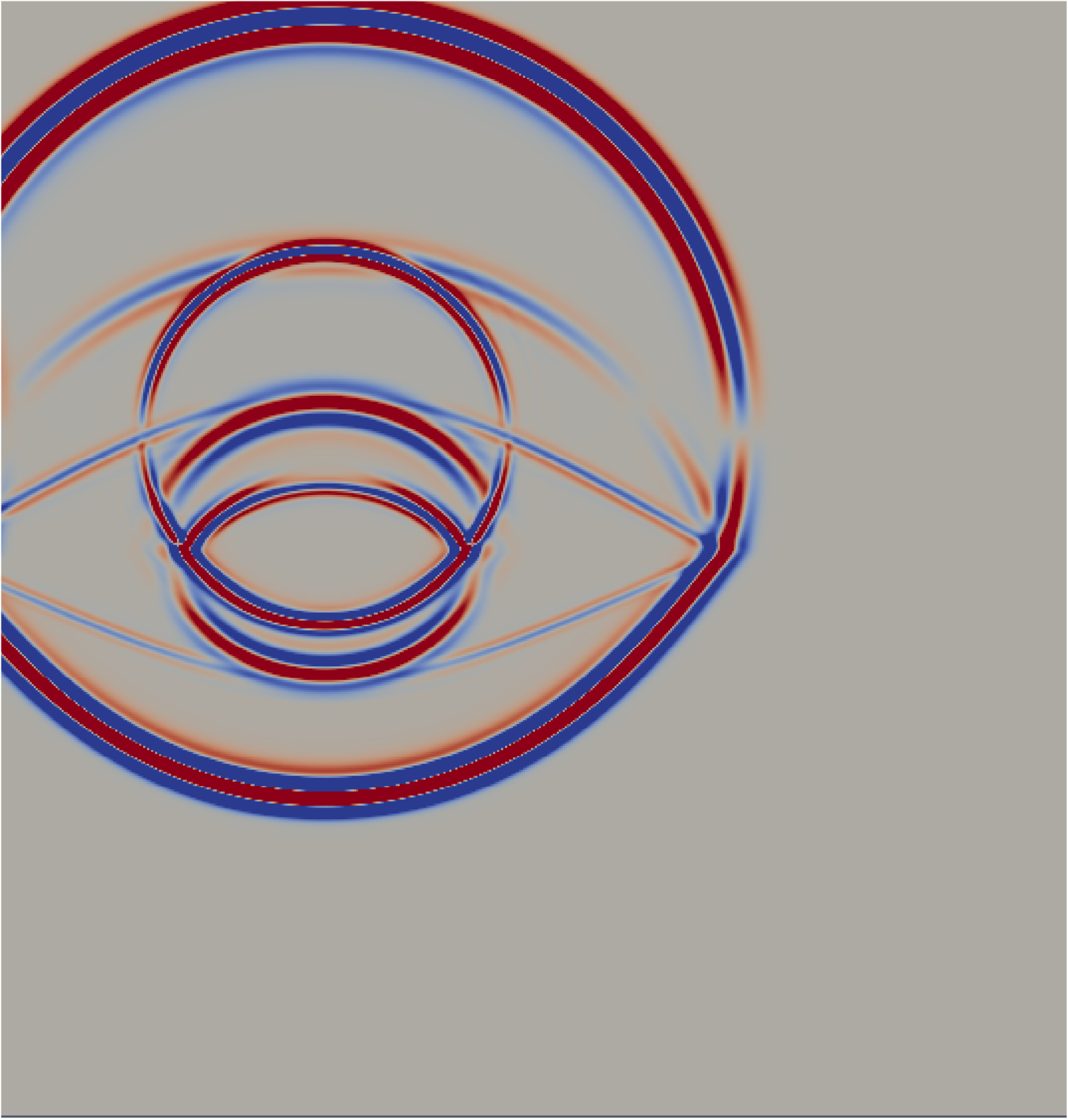};
\node at (axis cs:1.6,2.9)[star,star points=7,draw=black, fill=white, inner sep=0pt,minimum size=7pt]  {};
\node at (axis cs:2.0,1.867)[diamond,draw=black, fill=red, inner sep=0pt,minimum size=7pt] {};
\node at (axis cs:3.0,1.8)[draw=black, fill=white] {$(\text{i})$};
\node at (axis cs:2.5,2.3)[draw=black, fill=white] {$\text{(ii})$};
\node at (axis cs:1.5,1.9)[draw=black, fill=white] {$\text{(iii})$};
\node at (axis cs:2.0,2.4)[draw=black, fill=white] {$\text{(iv})$};
\end{axis}
\end{tikzpicture}
}\\
\subfloat[SPECFEM vs WADG]{
\includegraphics[width=0.75\textwidth]{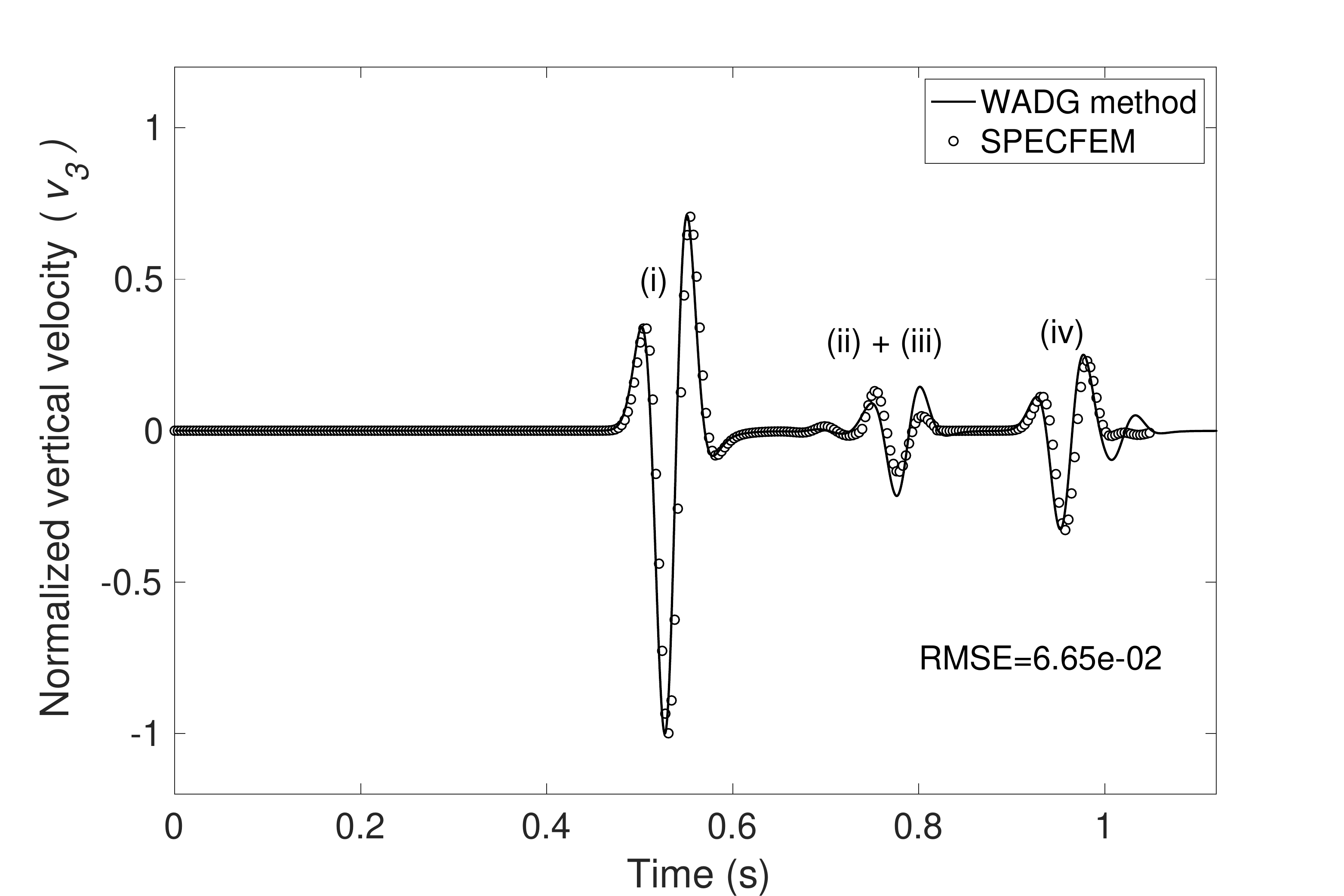}
}

\caption{A comparison between the solutions obtained from the proposed DG method and spectral-element method in a heterogeneous porous medium with (a) showing the snapshot of the solutions obtained from the DG method computed with polynomials of degree $N=4$. Figure (b) shows a comparison between the solutions obtained from the DG/spectral element methods. The labels (i)-(iv) correspond to different transmitted wave modes as follows: (i) denotes the transmitted fast P wave, (ii) denotes the conversion of fast P to S and slow P waves, (iii) denotes converted slow P to S waves, and (iv) denotes converted slow P waves to fast P waves.}
\end{figure}

\subsubsection{3D models and computational results}
Finally, we present numerical solutions of the poroelastic wave equation on 3D domains.  First, we compute the solution in an homogeneous medium and then extend the computation to a realistic synthetic reservoir model constructed from a varying layer with an undulated topography. 

\subsection{Epoxy-glass model}

First, we compute the solution for a 3D cube made of epoxy-glass with material properties given in column 3 of Table 1. The size of the computational domain is $2~\text{km} \times 2~\text{km} \times 2~\text{km}$ in $x$, $y$ and $z$ directions, respectively. The domain is discretized with tetrahedral elements with a minimum edge length of $45~\text{m}$.  Figure 13(a)-(b) represents the $x-$ and $z-$ components of the solid particle velocity in the epoxy-glass medium at $0.25~\text{s}$. The central frequency of the forcing function is $ f_0=35~\text{Hz}$, and polynomials of degree $N=5$ are used.  

Three events can be observed: the fast P mode ($\text{P}_ \text{f}$, outer wavefront), the shear wave (S, middle wavefront), and the slow P mode (inner wavefront).  In an orthotropic material, the effect of anisotropy on wave propagation is observed in two orthogonal planes, whereas isotropic effects are observed only in one plane. Figure 13a clearly shows the effect of anisotropy on wave propagation with an isotropic effect in the plane parallel to the $y$ axis.

\begin{figure}
\centering
\subfloat[Snapshot of horizontal solid particle velocity, $v_1$ at $t=0.25~\text{s}$ ]{
\includegraphics[trim={5.5cm, 1cm, 7cm, 1cm}, clip, width=0.80\textwidth]{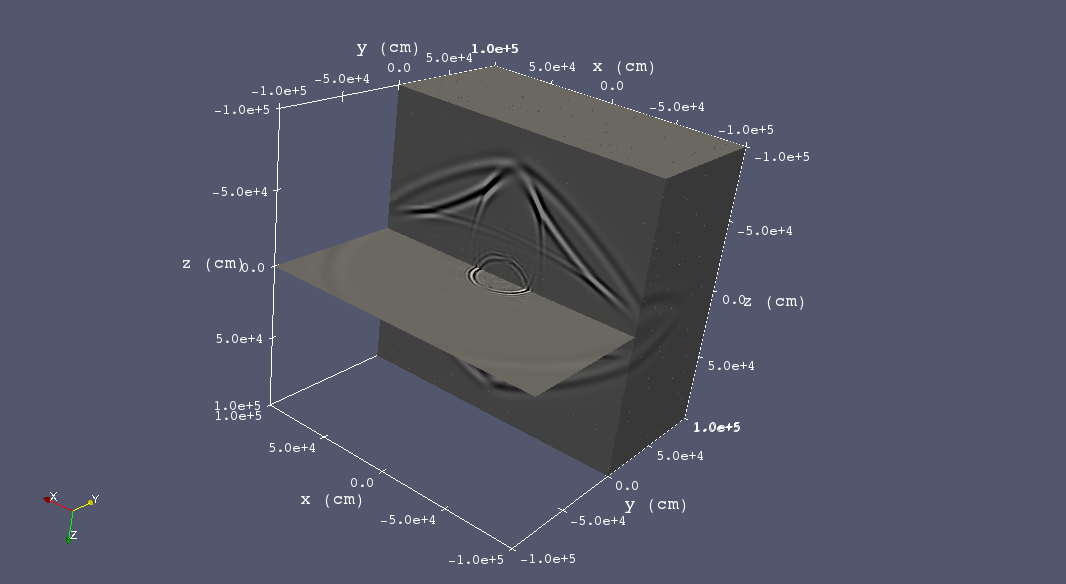}
}\\
\subfloat[Snapshot of vertical solid particle velocity, $v_3$ at $t=0.25~\text{s}$  ]{
\includegraphics[trim={5.5cm, 1cm, 7cm, 1cm}, clip, width=0.80\textwidth]{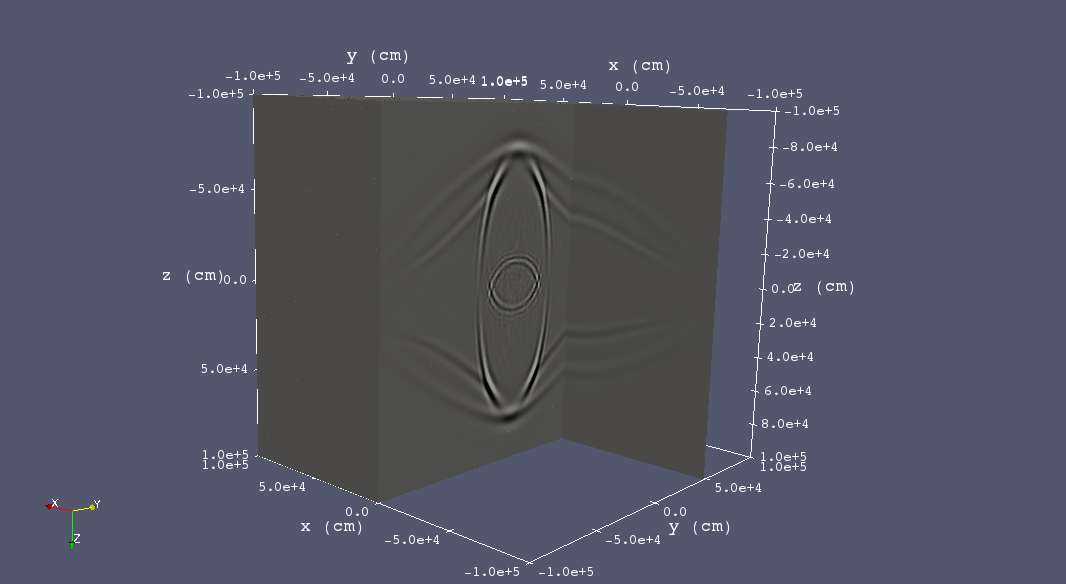}
}
\caption{Snapshots of the $x$ and $z$ components of solid particle velocity in a 3D epoxy-glass material. (a) and (b) corresponds to the snapshots of $v_1$ and $v_3$ at $t=0.25~\text{s}$. The central frequency of the forcing function is 25 Hz. The solution is computed using polynomials of degree $N=5$.}
\end{figure}

\subsection{Reservoir model}
We construct a 3D reservoir model characterized by rock layers, discontinuity, and a surface with undulated topography. The discretized model along with the density distribution is shown in Figures 14a and 14b, respectively. The dimension of the model is $22.8~\text{km}\times17.4~\text{km}\times 7.0~\text{km}$ in $x$, $y$ and $z$ directions, respectively. The domain is discretized with tetrahedral elements with a minimum edge length of $125~\text{m}$. The top surface of the model is perturbed so that the effects of the topography, assumed as a free surface, could be incorporated into numerical simulations. Figure 15(a)-(b) represent the $x-$ and $z-$ components of the solid particle velocity at $3.5~\text{s}$. The central frequency of the forcing function is $f_0=10~\text{Hz}$, and polynomials of degree $N=3$ are used. The various modes of transmissions, reflections and scattering can be clearly seen in Figures 15a and 15b. Figure 16(a)-(b) show the $x-$ and $z-$ components of the fluid particle velocity at $3.5~\text{s}$.  

\begin{figure}
\centering
\subfloat[Discretized reservoir model ]{
\includegraphics[trim={5.5cm, 1cm, 7cm, 1cm}, clip, width=0.80\textwidth]{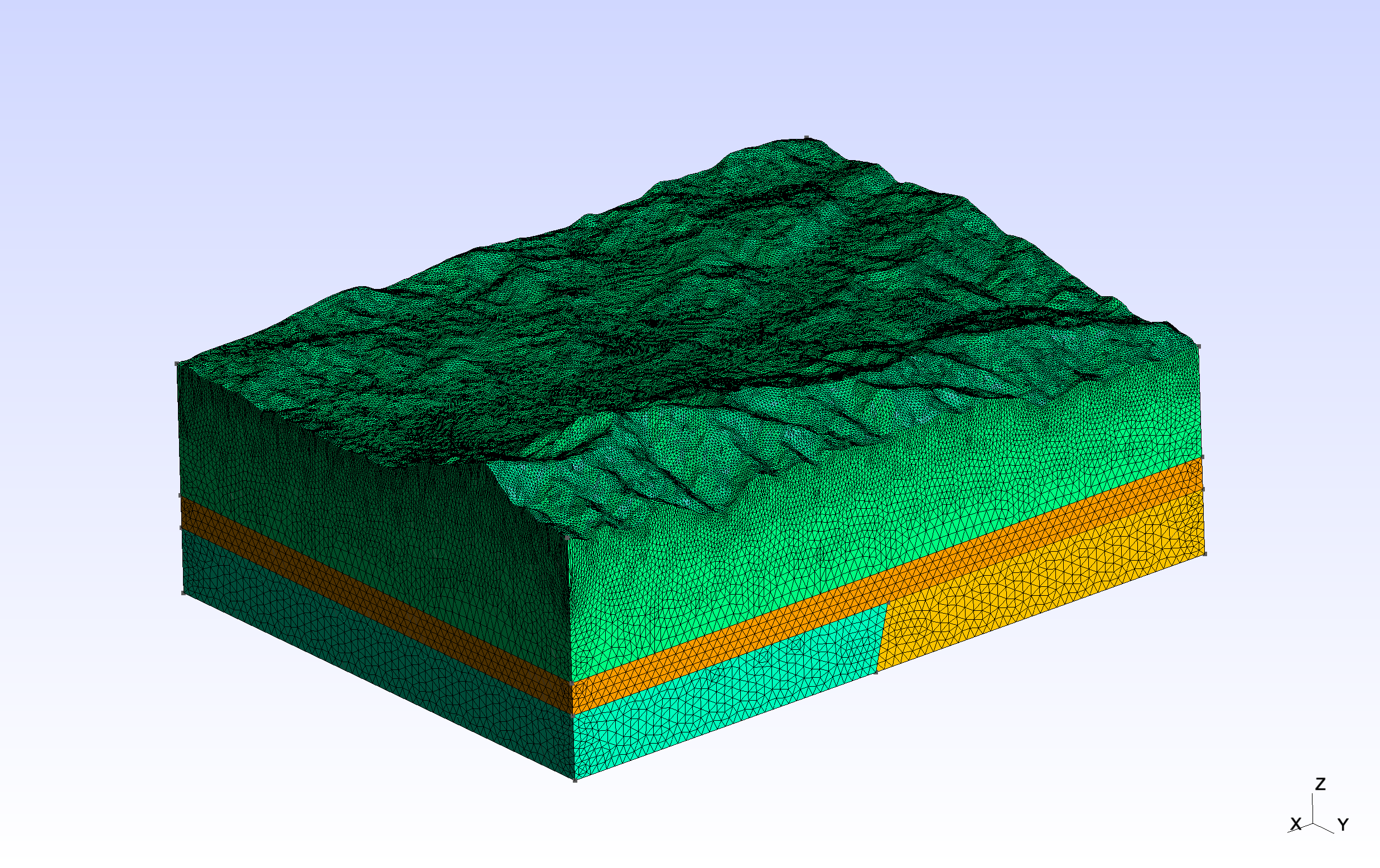}
}\\
\subfloat[Density $(\rho)$ model ]{
\includegraphics[trim={5.5cm, 0cm, 1cm, 1cm}, clip, width=0.80\textwidth]{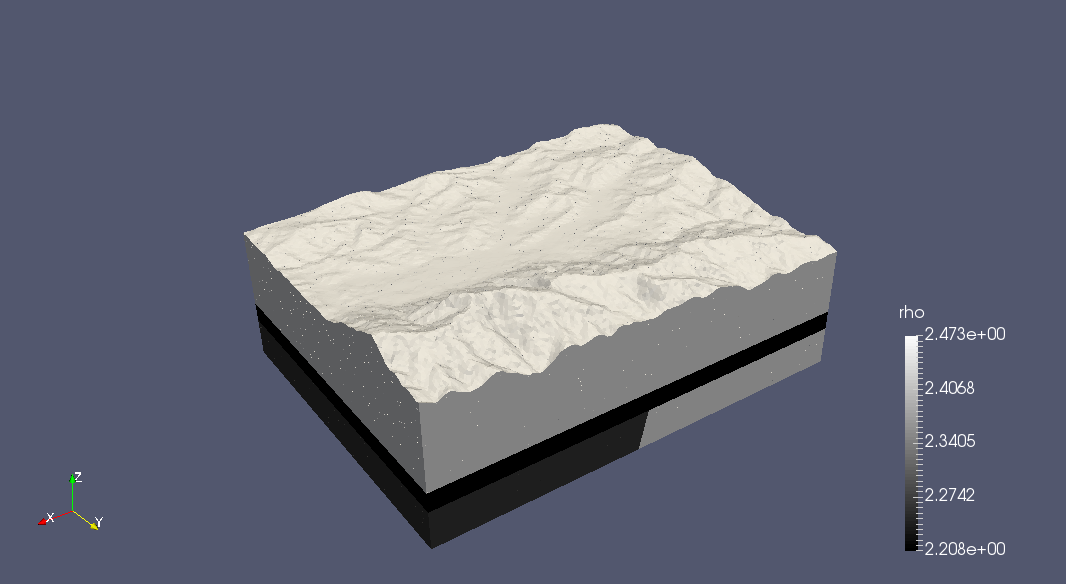}
}
\caption{3D Reservoir model (a) model discretized with tetrahedral elements and (b) distribution of the density $(\rho)$ in the model.}

\end{figure}

\begin{figure}
\centering
\subfloat[Snapshot of horizontal solid particle velocity, $v_1$ at $t=3.5~\text{s}$ ]{
\includegraphics[trim={1cm, 1cm, 7cm, 1cm}, clip, width=0.80\textwidth]{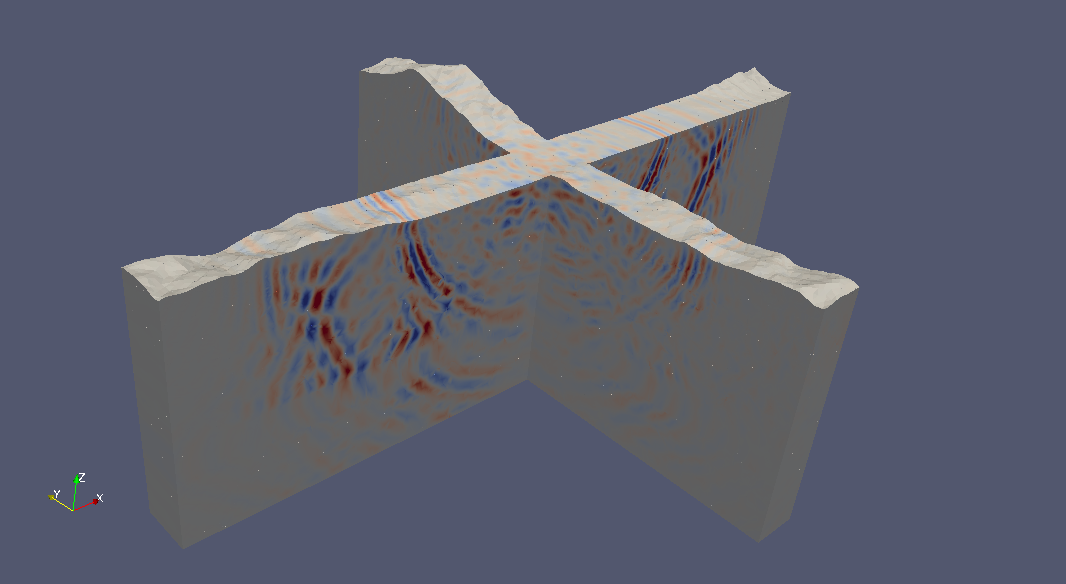}
}\\
\subfloat[Snapshot of vertical solid particle velocity, $v_3$ at $t=3.5~\text{s}$ ]{
\includegraphics[trim={1cm, 1cm, 7cm, 1cm}, clip, width=0.80\textwidth]{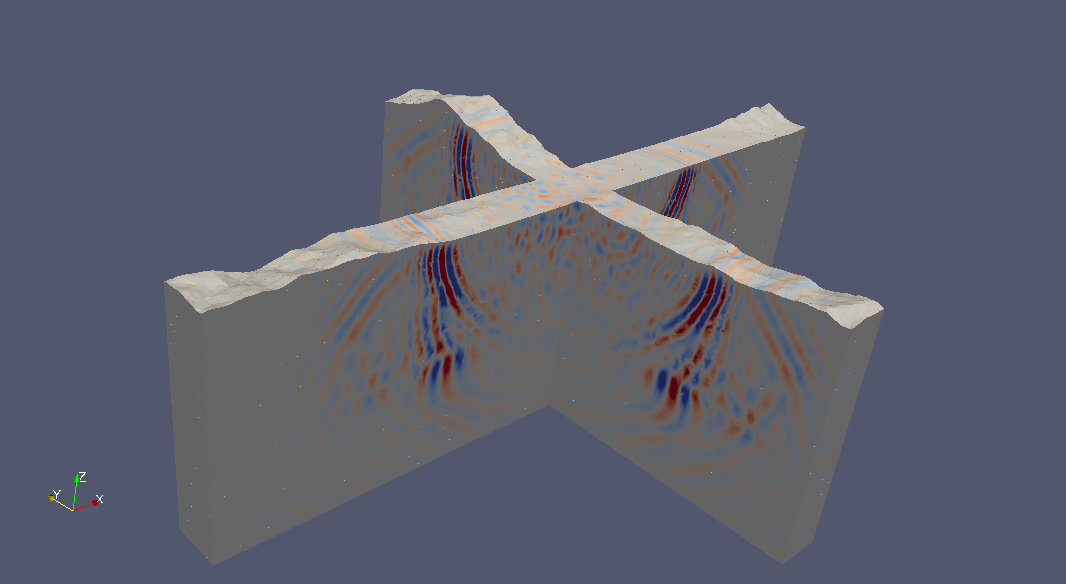}
}
\caption{Snapshots of the $x$ and $z$ components of solid particle velocity in a reservoir model. (a) and (b) corresponds to the snapshots of $v_1$ and $v_3$ at $t=3.5~\text{s}$. The central frequency of the forcing function is 10 Hz. The solution is computed using polynomials of degree $N=3$.}

\end{figure}

\begin{figure}
\centering
\subfloat[Snapshot of horizontal fluid particle velocity, $q_1$ at $t=3.5~\text{s}$ ]{
\includegraphics[trim={1cm, 1cm, 7cm, 1cm}, clip, width=0.80\textwidth]{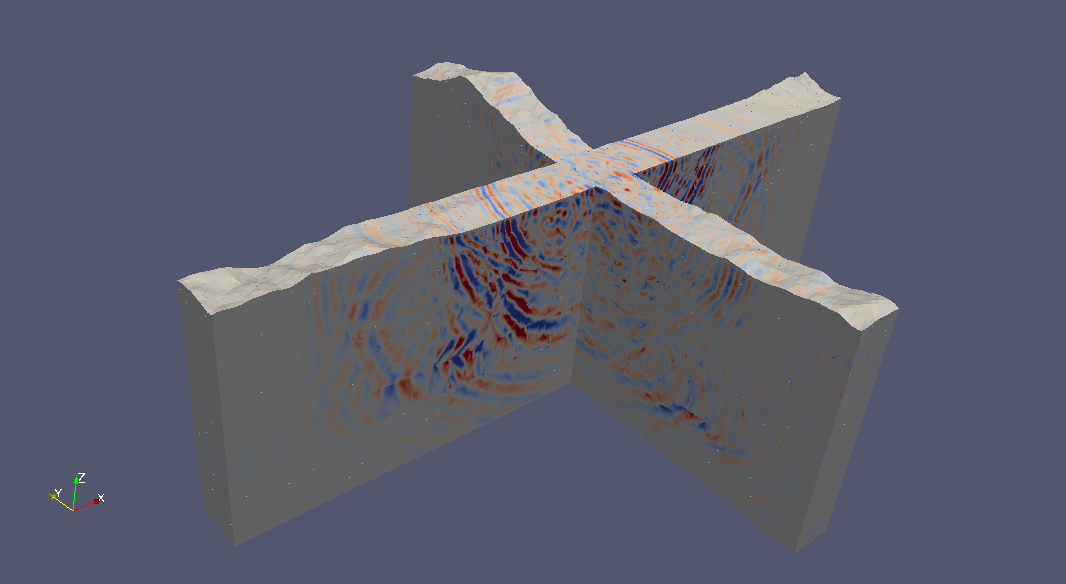}
}\\
\subfloat[Snapshot of vertical fluid particle velocity, $q_3$ at $t=3.5~\text{s}$ ]{
\includegraphics[trim={1cm, 1cm, 7cm, 1cm}, clip, width=0.80\textwidth]{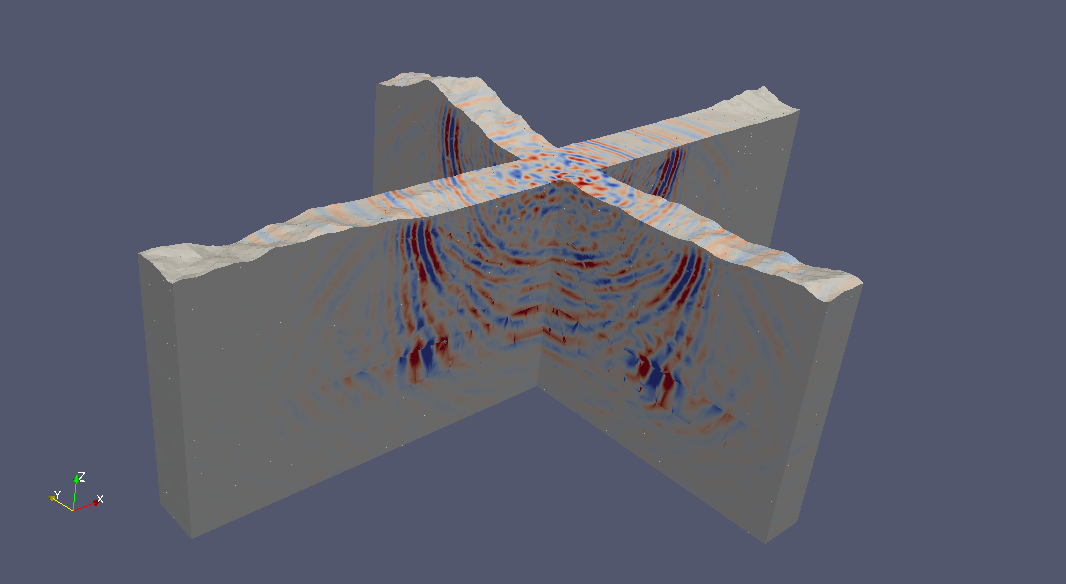}
}
\caption{Snapshots of the $x$ and $z$ components of fluid particle velocity in a reservoir model. (a) and (b) corresponds to the snapshots of $q_1$ and $q_3$ at $t=3.5~\text{s}$. The central frequency of the forcing function is 10 Hz. The solution is computed using polynomials of degree $N=3$.}

\end{figure}

\section{Conclusions}
This work presents a weight-adjusted discontinuous Galerkin (WADG) method for the linear poroelastic wave equations with arbitrary heterogeneous media. The method is energy stable and high order accurate for arbitrary stiffness tensors. The penalty numerical fluxes for this formulation are simple to derive and implement, and their lack of dependence on the stiffness tensors allows for a unified treatment of isotropic and anisotropic media with micro-heterogeneities. 

We confirm the high-order accuracy of the numerical method using an analytic plane wave solution in a poroelastic media. Results obtained using this method also show good agreement with existing results in the literature for both problems involving both homogeneous and heterogeneous media. Finally, we provide computational results for a large 3D model consisting the surface topography. 

We note that the implementation of this method reduces to the application of the weight-adjusted mass matrix inverse and the evaluation of constant-coefficient terms in the DG formulation. The cost of the latter step can be reduced (especially at high orders of approximation) by using fast methods based on Bernstein-Bezier polynomials for the application of derivative and lift matrices for constant-coefficient terms \cite{chan2017gpu}.  Future work will also involve an extension of the proposed formulation for the broad-band Biot's model and coupled acousto-poroelastic media.

\section{Acknowledgments}

The authors gratefully thank sponsors of the Geo-Mathematical Imaging Group at Rice University for providing the resources to carry out this work.  Jesse Chan gratefully acknowledges support from the NSF under grants DMS-1719818 and DMS-1712639.  MVdH gratefully acknowledges support from the Simons Foundation under the MATH + X program and the NSF under grant DMS-1559587. The authors gratefully thank Ruichao Ye, Vitaly Katsnelson, Peter Caday, Christopher Wong and Sundeep Sharma for helpful and informative discussions. We also acknowledge Dr.\ Octavio Castillo from Barcelona Supercomputing Center (BSC), Barcelona, Spain for his help in 3D meshing of reservoir model.

\bibliographystyle{unsrt}
\bibliography{dg}

\newpage
\appendix
\section{Matrix weight-adjusted inner products}
\label{appendix:A0}
Let $D^{\alpha} \bm{v}$ denote the component-wise differentiation of $\bm{v}$ with respect to a $d-$dimensional multi-index $\alpha$. Then, vector $L^p$ Soblev norms for $\bm{v(x)} \in \mathbb{R}$ can be defined as
\[
|\bm{v}|_{W^{k,p}}^p = \sum_{i=1}^m |\bm{v}_i|_{W^{k,p}}^p, \qquad||\bm{v}||_{W^{k,p}}^p = \sum_{i=1}^m ||\bm{v}_i||_{W^{k,p}}^p\qquad1 \le p \le \infty,
\]
\[
|\bm{v}|_{W^{k,\infty}}=\max_i |\bm{v}_i|_{W^{k, \infty}}, \qquad ||\bm{v}||_{W^{k,\infty}}=\max_i ||\bm{v}_i||_{W^{k, \infty}},
\]
where $W^{k,p}$ and $W^{k, \infty}$ are the Sobolev spaces \cite{chan2017III}.
Let $\bm{W}(\bm{x})$ be a matrix-valued weight function which is pointwise symmetric positive definite 
\[
0 < w_{min} \le ||\bm{W}(\bm{x})||_2 \le w_{max} < \infty, ~~~ 0 < \tilde{w}_{min} \le ||\bm{W}^{-1} (\bm{x})||_2 \le \tilde{w}_{max} < \infty,\qquad~\forall \bm{x} \in \Omega.
\] 
The $k^{th}$ order Sobolev norm for $\bm{W}(\bm{x})$ is defined as 
\[
||\bm{W}(\bm{x})||_{k,p,\infty}^p = \sum_{|\alpha| \le k} \sup_{x} ||D^{\alpha}\bm{W}(\bm{x})||_p^p.
\]

\section{Weight-adjusted approximations with matrix weights}
\label{appendix:A02}
Let $\Pi_N\bm{u}$ defined as $L^2$ projection applied to each component of the vector-valued function $\bm{u}$. We then define the operator $T_{\bm{W}}$ as
\begin{align*}
T_{\bm{W}} \bm{v}=\Pi_N(\bm{Wv}).
\end{align*}
The inverse operator $T_W^{-1}$ is defined implicitly via
\[
\left(\bm{W}T_{\bm{W}}^{-1}\bm{v}, \bm{\delta v}\right)_{L^2(D^k)}=(\bm{v}, \bm{\delta v})_{L^2(D^k)}, \qquad~~\forall \bm{\delta v} \in (P^N (D^k))^m.
\]

\begin{lemma}
\label{lem:1}
Let $\Pi_{N}$ denote the component-wise $L^2$ projection and $\bm{W} \in {(L^{\infty})}^{m \times m}$. Then $T_{\bm{W}}$ satisfies the following properties:
\begin{enumerate}
\item $T_{\bm{W}}^{-1} T_{\bm{W}}=\Pi_N$
\item $\Pi_{N} T^{-1}_{\bm{W}}=\bm{T}_{\bm{W}}^{-1}=T_{\bm{W}}^{-1}$
\item $|| T_{\bm{W}}^{-1}||_{L^2(D^k)} \le \tilde{w}_{\max}.$
\item $(T_{W}^{-1}\bm{v}, \bm{w})_{L^2(D^k)}$ forms an inner product on $(P^{N})^m \times (P^N)^m$, which is equivalent to the $L^2$ inner product with equivalence constants $C_1=\tilde{w}_{\min},~C_2=\tilde{w}_{\max}$.
\end{enumerate}
\end{lemma}
The proof of Lemma 1 is given in detail in \cite{chan2017I}.

\section{Convergence analysis of WADG scheme for poroelastic wave equations}
\label{appendix:A}
Using the results in Section 4.2, we can extend the semi-discrete convergence analysis in \cite{houston2002,chan2017IV,warburton2013} in to linear poroelastic wave propagation on meshes of affine elements. Let $\bm{U}$ and $\bm{U}_h$ represent the exact and discrete WADG solutions, respectively. We assume that $\bm{U}$ and $\myfrac{\partial \bm{U}}{\partial t}$ are sufficiently regular such that
\[
\left \Vert \bm{U}\right \Vert_{W^{N+1,2}{(\Omega_h)}} < \infty, \qquad\left \Vert\myfrac{\partial \bm{U}}{\partial t}\right\Vert_{W^{N+1,2}{(\Omega_h)}} < \infty,
\]
with ${\left\Vert \bm{U} \right \Vert}^2_{W^{N+1,2}{(\Omega_h)}}=\sum_k{\left\Vert \bm{U} \right \Vert}^2_{W^{N+1,2}{(D^k)}}$.

The WADG formulation can be written as
\[
\left(T^{-1}_{\bm{A}_0^{-1}} \myfrac{\partial \bm{U}}{\partial t}, \bm{V}\right)_{L^2(\Omega)}
+ a(\bm{U}, \bm{V}) + b(\bm{U}, \bm{V}) =(\bm{f}, \bm{V}),
\]

\[
a(\bm{U}, \bm{V})=\sum_{D^k \in \Omega_h}\left(-\left( \sum_{i=1}^d \bm{\tau}, \bm{A}_i \myfrac{\partial \bm{g}}{\partial \bm{x}_i}   \right)_{L^2(D^k)} + \left( \sum_{i=1}^d \bm{A}_i \myfrac{\partial \bm{v}}{\partial \bm{x}_i}, \bm{g}\right)_{L^2(D^k)} + (\bm{Dv}, \bm{g})_{L^2(D^k)} \right),
\]

\[
b(\bm{U}, \bm{V})=\sum_{D^k \in \Omega_h}\left( \left\langle \bm{A}_n^T \avg{\bm{\tau}} + \myfrac{\alpha_{\bm{v}}}{2}\bm{A}_n^T \bm{A}_n \jump{\bm{v}}, \bm{g} \right \rangle_{L^2(D^k)} + \left\langle \myfrac{1}{2} \bm{A}_n \jump{\bm{v}} + \myfrac{\alpha_{\bm{\tau}}}{2} \bm{A}_n \bm{A}_n^T \jump{\bm{\tau}}, \bm{h} \right \rangle_{L^2(D^k)}\right),
\]
where $\bm{A}_0(\bm{x})=\text{diag} (\bm{Q}_s, \bm{Q}_v)$, $\bm{U}$ and $\bm{V}$ are group variables, defined as $\bm{U}=(\bm{\tau}, \bm{v})$ and $\bm{V}=(\bm{h}, \bm{g}) \in (V_h (\Omega_h))^d \times (V_h(\Omega_h))^{N_d}$.

The DG formulation in (\ref{eq17}) is consistent and thus
\begin{align}
\left(\bm{A}_0^{-1} \myfrac{\partial \bm{U}}{\partial t}, \bm{V}\right)_{L^2(\Omega)}
+ a(\bm{U}, \bm{V}) + b(\bm{U}, \bm{V}) =(\bm{f}, \bm{V}), \label{eqdg} \\
\left(T^{-1}_{\bm{A}_0^{-1}} \myfrac{\partial \bm{U}_h}{\partial t}, \bm{V}\right)_{L^2(\Omega)}
+ a(\bm{U}_h, \bm{V}) + b(\bm{U}_h, \bm{V}) =(\bm{f}, \bm{V}), \label{eqwadg}
\end{align}
$\forall~\bm{V} \in (V_h (\Omega_h))^d \times (V_h(\Omega_h))^{N_d}$.
We decompose the error $\bm{U}-\bm{U}_h$ into a projection error $\bm{\epsilon}$ and a discretization error $\bm{\eta}$.
\[
\bm{U}-\bm{U}_h = (\Pi_N \bm{U}-\bm{U}_h) - (\Pi_N \bm{U}-\bm{U})=\bm{\eta}-\bm{\epsilon}.
\] 

We assume that $\bm{U}_h(\bm{x}, 0)$ is the $L^2$ projection of the exact initial condition, such that $\bm{\eta}|_{t=0}=0$. We also introduce a consistency error $\delta=\bm{A}_0\bm{U}-T^{-1}_{\bm{A}_0^{-1}}\bm{U}$ resulting from the approximation of $\bm{A}_0 \bm{U}$ by a weight-adjusted inner product 
\[
\bm{A}_0 \myfrac{\partial \bm{U}}{\partial t} - T^{-1}_{\bm{A}_0^{-1}}\myfrac{\partial \bm{U}_h}{\partial t}= \myfrac{\partial }{\partial t}\left(\bm{A}_0 \bm{U}-T^{-1}_{\bm{A}_0^{-1}} \bm{U} \right) + \myfrac{\partial }{\partial t} T^{-1}_{\bm{A}_0^{-1}} (\Pi_N \bm{U}-\bm{U}_h)=\myfrac{\partial \delta}{\partial t} + \myfrac{\partial }{\partial t}\left(T^{-1}_{\bm{A}_0^{-1}}\right),
\]
where we have used that $T^{-1}_{\bm{A}_0^{-1}} \Pi_N$. Subtracting the DG and WADG formulations in (\ref{eqdg}) and (\ref{eqwadg}) and setting $\bm{V}=\bm{\eta}$ yields
\begin{align}
\label{eq27}
\myfrac{1}{2}\myfrac{\partial}{\partial t}\left(T^{-1}_{\bm{A}_0^{-1}}\bm{\eta}, \bm{\eta}\right)_{L^2(\Omega)} + (\bm{\eta}, \bm{\eta})=\left(-\myfrac{\partial \bm{\delta}}{\partial t}, \bm{\eta}\right)_{L^2(\Omega)} + a(\bm{\epsilon}, \bm{\eta}) + b(\bm{\epsilon}, \bm{\eta}),
\end{align}
where we have used that $a(\bm{\eta}, \bm{\eta})=0$ from skew-symmetry.
We bound $a(\bm{\epsilon}, \bm{\eta}) + b(\bm{\epsilon}, \bm{\eta})$ in (\ref{eq27}) by integrating by parts the stress equation and using the component-wise $L^2$ orthogonality of $\bm{\epsilon}$ to derivatives of $\bm{\eta}$. This reduces to
\begin{align*}
\sum_{D^k \in \Omega_h} &\left(\left \langle  \bm{A}_n^T \avg{\bm{\epsilon}_{\tau}} + \myfrac{\alpha_{\bm{v}}} {2} \bm{A}_n^T \bm{A}_n \jump{\bm{\epsilon}_v}, \bm{\eta}_v  \right \rangle_{L^2(\partial D^k)} + \left \langle  \bm{A}_n \avg{\bm{\epsilon}_{v}} + \myfrac{\alpha_{\bm{\tau}}} {2} \bm{A}_n\bm{A}_n^T  \jump{\bm{\epsilon}_{\bm{\tau}}}, \bm{\eta}_{\tau}  \right \rangle_{L^2(\partial D^k)} - \left \langle \bm{D}\bm{\epsilon}_{v},\bm{\epsilon}_{v}\right \rangle_{L^2(D^k)}  \right)\\
&
=\myfrac{1}{2}\sum_{D^k \in \Omega_h}\Biggl( \left \langle  \avg{\bm{\epsilon}_{\tau}} - \myfrac{\alpha_{\bm{v}}} {2} \bm{A}_n \jump{\bm{\epsilon}_v}, \bm{A}_n \jump{\bm{\eta}_v}  \right \rangle_{L^2(\partial D^k)} + \left \langle  \avg{\bm{\epsilon}_{v}} - \myfrac{\alpha_{\bm{\tau}}} {2} \bm{A}_n^T \jump{\bm{\epsilon}_{\tau}}, \bm{A}_n^T \jump{\bm{\eta}_{\tau}}  \right \rangle_{L^2(\partial D^k)}\\& - \left \langle \bm{D}\bm{\epsilon}_{v},\bm{\epsilon}_{v}\right \rangle_{L^2(D^k)}\Biggr)\\
&
\le \myfrac{1}{2}\sum_{D^k \in \Omega_h} \left \Vert \avg{\bm{\epsilon}_{\tau}} - \myfrac{\alpha_{\bm{v}}} {2} \bm{A}_n \jump{\bm{\epsilon}_v} \right \Vert_{L^2(\partial D^k)} \left \Vert  \bm{A}_n \jump{\bm{\eta}_v} \right \Vert_{L^2(\partial D^k)} \\&+ \left \Vert \avg{\bm{\epsilon}_{v}} - \myfrac{\alpha_{\bm{\tau}}} {2} \bm{A}_n^T \jump{\bm{\epsilon}_{\tau}}  \right \Vert_{L^2(\partial D^k)} \left \Vert \bm{A}_n^T \jump{\bm{\eta}_{\tau}} \right \Vert_{L^2(\partial D^k)}\\
&
\le C_{\tau} \sum_{D^k \in \Omega_h} \left \Vert \bm{\epsilon} \right \Vert_{L^2\left(\partial D^k\right)} \left( \myfrac{\alpha_{\bm{v}}}{2} \left \Vert \bm{A}_n \jump{\bm{\eta}_v}\right \Vert^2 _{L^2(\partial D^k)} + \myfrac{\alpha_{\bm{\tau}}}{2} \left \Vert \bm{A}_n^T \jump{\bm{\eta}_{\tau}}\right \Vert^2 _{L^2(\partial D^k)}\right)^{1/2}
\end{align*}
where $C_{\tau}$ is maximum of $(\alpha_{\bm{\tau}}, \alpha_{\bm{v}})$ and we have also used the property of negative sime-definite matrix $\bm{D}$. Using the Young's inequality with $\alpha=C_{\tau}/2$ yields the following bound
\[
|a(\bm{\epsilon}, \bm{\eta}) + b(\bm{\epsilon}, \bm{\eta})| \le b(\bm{\eta}, \bm{\eta}) + \myfrac{C_{\tau}^2}{4}\sum_{D^k \in \Omega_h} \left\Vert \bm{\epsilon} \right \Vert_{L^2(\partial D^k)}^2.
\]
Applying this to (\ref{eq27}) and using Cauchy-Schwarz on $\left( \myfrac{\partial \delta}{\partial t}, \bm{\eta}\right)_{L^2(\partial D^k)}$ yields
\begin{align}
\label{eq28}
\myfrac{1}{2} \myfrac{\partial }{\partial t} \left(T^{-1}_{\bm{A}_0^{-1}}\bm{\eta}, \bm{\eta}\right)_{L^2(\Omega)} + b(\bm{\eta}, \bm{\eta}) \le \left \Vert \myfrac{\partial \delta}{\partial t} \right  \Vert_{L^2(\partial D^k)} \left\Vert \bm{\eta} \right \Vert _{L^2(\partial D^k)} + b(\bm{\eta}, \bm{\eta}) + \sum_{D^k \in \Omega_h} \myfrac{C_{\tau}^2}{4} \left \Vert \bm{\epsilon} \right \Vert_{L^2(\partial D^k)}
\end{align}
Using a  $hp$ trace inequality \cite{warburton2003}
\begin{align}
\label{eq28}
\sum_{D^k \in \Omega_h} \left \Vert \bm{\epsilon}^2\right \Vert_{L^2(\partial D^k)}  \le \sum_{D^k \in \Omega_h} Ch^{-1} \left \Vert \bm{\epsilon}^2 \right \Vert_{L^2(\partial D^k)}=Ch^{-1}\left \Vert \bm{\epsilon}_{L^2{(\Omega)}} \right \Vert \le Ch^{2N+1} \left \Vert \bm{U} \right \Vert_{W^{N+1,2}}.
\end{align}
We now use that that $\bm{Q}_s$ and $\bm{Q}_v$ are positive definite such that
\[
0 < s_{\text{min}} \le \bm{u}^T \bm{Q}_s(\bm{x})\bm{u} \le s_{\text{max}}  < \infty
\] 
\[
0 < \tilde{s}_{\text{min}} \le \bm{u}^T \bm{Q}_s(\bm{x})\bm{u} \le \tilde{s}_{\text{max}}  < \infty
\]
\[
0 < v_{\text{min}} \le \bm{u}^T \bm{Q}_v(\bm{x})\bm{u} \le v_{\text{max}}  < \infty
\] 
\[
0 < \tilde{v}_{\text{min}} \le \bm{u}^T \bm{Q}_v(\bm{x})\bm{u} \le \tilde{v}_{\text{max}}  < \infty
\] 
We define $A_{\max}=\max\left({s}_{\text{max}}\tilde{s}_{\text{max}}, {v}_{\text{max}}\tilde{v}_{\text{max}} \right), ~A_{\min}=\min\left(\tilde{s}_{\text{min}}, \tilde{v}_{\text{min}} \right)$ and substitute the result of (\ref{eq28}) in (\ref{eq27}) performing time integration from $[0,T]$ and using  modified Gronwall's inequality (Lemma 1.10 in \cite{dolejvsi2015}), (\ref{eq28}) reduces to
\begin{align*}
\left \Vert \bm{\eta} \right \Vert _{L^2(\Omega)} \le \myfrac{CTh^{N+1/2}}{A_{\min}} \sup_{t \in [0,T]} \left(h^{1/2}A_{\max} \left \Vert \bm{A}_0 \right \Vert_{W^{N+1, \infty} (\Omega_h)}\left \Vert \myfrac{\partial \bm{U}}{\partial t} \right \Vert			_{L^2{(\Omega)}} + \left \Vert \bm{U} \right \Vert_{W^{N+1,2}(\Omega_h)} \right).
\end{align*}
Using the triangle inequality yields
\[
\left \Vert \bm{U} - \bm{U}_h \right \Vert \le (C_1 + C_2 T)h^{N+1/2} \sup_{t \in [0,T]} \left(h^{1/2}A_{\max} \left \Vert \bm{A}_0 \right \Vert_{W^{N+1, \infty} (\Omega_h)}\left \Vert \myfrac{\partial \bm{U}}{\partial t} \right \Vert			_{L^2{(\Omega)}} + \left \Vert \bm{U} \right \Vert_{W^{N+1,2}(\Omega_h)} \right).
\]
where $C_1$ and $C_2$ depends on $A_{\min}$ and $A_{\max}$. These estimates show that $L^2$ errors decrease proportionally to $h^{N+1/2}$. In practice we often observe the rate of convergence to be $O(h^{N+1})$.

\section{Analytical solution of diffusive part of the poroelastic system}
From (\ref{eq9}), the diffusive part of the system is expressed as
\label{appendix:B}
\begin{align}
\bm{Q}_v \myfrac{\partial \bm{v}}{\partial t}=\bm{D}\bm{v}.
\end{align}
The component wise expression for intermediate velocity vector is
\begin{align}
\myfrac{\partial v_1}{\partial t}&=-\myfrac{\eta}{\kappa1}\beta_{v}^{(1)}q_1,\\
\myfrac{\partial v_2}{\partial t}&=-\myfrac{\eta}{\kappa2}\beta_{v}^{(2)}q_2,\\
\myfrac{\partial v_2}{\partial t}&=-\myfrac{\eta}{\kappa3}\beta_{v}^{(3)}q_3,\\
\myfrac{\partial q_1}{\partial t}&=-\myfrac{\eta}{\kappa1}\beta_{q}^{(1)}q_1,\\
\myfrac{\partial q_2}{\partial t}&=-\myfrac{\eta}{\kappa2}\beta_{q}^{(2)}q_2,\\
\myfrac{\partial q_2}{\partial t}&=-\myfrac{\eta}{\kappa3}\beta_{q}^{(3)}q_3,
\end{align} 
where $\beta_{v}^{(i)}=\myfrac{\rho_f}{\rho_f^2-\rho m_i}$ and $\beta_{q}^{(i)}=\myfrac{\rho}{\rho_f^2-\rho m_i}$. The above system can be solved analytically, giving
\begin{align}
v_1^{*}&=v_1^n + \myfrac{\beta_v^{(1)}}{\beta_q^{(1)}}\left[\exp(\lambda_s^{(1)}dt)-1\right]q_1^n,\\
v_2^{*}&=v_2^n + \myfrac{\beta_v^{(2)}}{\beta_q^{(2)}}\left[\exp(\lambda_s^{(2)}dt)-1\right]q_2^n,\\
v_3^{*}&=v_3^n + \myfrac{\beta_v^{(3)}}{\beta_q^{(3)}}\left[\exp(\lambda_s^{(3)}dt)-1\right]q_3^n,\\
q_1^{*}&=\exp(\lambda_s^{(1)}dt)q_1^n,\\
q_2^{*}&=\exp(\lambda_s^{(2)}dt)q_2^n,\\
q_3^{*}&=\exp(\lambda_s^{(3)}dt)q_3^n,
\end{align}
where $\lambda_s^{(i)}=-(\eta/\kappa_i)\beta_q^{(i)}$. The state of intermediate vector $[v_1^{*},~v_2^{*},~v_3^{*},~q_1^{*},~q_2^{*},~q_3^{*}]$ is the input for the non-stiff being solved with DG scheme.

\end{document}